\documentclass[reqno,10pt,centertags,draft]{amsart}
\usepackage{amsmath,amsthm,amscd,amssymb,latexsym,upref}
\date{\today}

\newcommand{\bbN}{{\mathbb{N}}}
\newcommand{\bbR}{{\mathbb{R}}}

\newcommand{\bbZ}{{\mathbb{Z}}}
\newcommand{\bbC}{{\mathbb{C}}}
\newcommand{\bbQ}{{\mathbb{Q}}}

\newcommand{\cA}{{\mathcal A}}

\newcommand{\cC}{{\mathcal C}}

\newcommand{\cH}{{\mathcal H}}

\newcommand{\cK}{{\mathcal K}}

\newcommand{\cR}{{\mathcal R}}

\newcommand{\dott}{\,\cdot\,}
\newcommand{\no}{\notag}
\newcommand{\lb}{\label}
\newcommand{\f}{\frac}

\newcommand{\ol}{\overline}

\newcommand{\wti}{\widetilde}

\newcommand{\Oh}{O}

\newcommand{\loc}{\text{\rm{loc}}}

\newcommand{\ran}{\text{\rm{ran}}}

\newcommand{\Arg}{\text{\rm{Arg}}}

\newcommand{\dom}{\text{\rm{dom}}}

\newcommand{\supp}{\text{\rm{supp}}}

\newcommand{\bi}{\bibitem}

\newcommand{\beq}{\begin{equation}}
\newcommand{\eeq}{\end{equation}}
\newcommand{\ba}{\begin{align}}
\newcommand{\ea}{\end{align}}

\newcommand{\veps}{\varepsilon}

\renewcommand{\Re}{\text{\rm Re}}
\renewcommand{\Im}{\text{\rm Im}}
\renewcommand{\ln}{\text{\rm ln}}

\renewcommand{\le}{\leqslant}

\DeclareMathOperator*{\wlim}{w-lim}
\DeclareMathOperator*{\simlim}{\sim}
\DeclareMathOperator*{\eqlim}{=}
\DeclareMathOperator*{\arrow}{\rightarrow}

\DeclareMathOperator{\Tr}{Tr}

\allowdisplaybreaks
\numberwithin{equation}{section}

\newtheorem{theorem}{Theorem}[section]

\newtheorem{lemma}[theorem]{Lemma}
\newtheorem{corollary}[theorem]{Corollary}
\newtheorem{hypothesis}[theorem]{Hypothesis}
\theoremstyle{definition}
\newtheorem{definition}[theorem]{Definition}
\newtheorem{remark}[theorem]{Remark}

\newtheorem{example}[theorem]{Example}

\begin{document}

\title[Inverse Spectral Theory]{Inverse Spectral Theory \\ As Influenced By Barry Simon}
\author[F.\ Gesztesy]{Fritz Gesztesy}
\address{Department of Mathematics,
University of Missouri, Columbia, MO 65211, USA}
\email{fritz@math.missouri.edu}
\urladdr{http://www.math.missouri.edu/personnel/faculty/gesztesyf.html}
\dedicatory{Dedicated with great pleasure to Barry Simon, mentor and
friend, \\ on the occasion of his 60th birthday.}

\thanks{Appeared in {\it Spectral Theory and Mathematical Physics: A Festschrift
in Honor of Barry Simon's 60th Birthday. Ergodic Schr\"odinger Operators, Singular Spectrum, Orthogonal Polynomials, and Inverse Spectral Theory}, F.\ Gesztesy, P.\ Deift, C.\ Galvez,
P.\ Perry, and W.\ Schlag (eds.), Proceedings of Symposia in Pure
Mathematics, Amer. Math. Soc., Vol.\ 76/2, Providence, RI, 2007, pp.\ 741--820.}
\date{\today}
\subjclass[2000]{Primary 34A55, 34B20, 34B24; Secondary 34C10, 34L05.}
\keywords{Spectral deformations, oscillation theory, trace formulas, Weyl--Titchmarsh theory, uniqueness results, inverse spectral theory.}

\begin{abstract}
We survey Barry Simon's principal contributions to the field of inverse
spectral theory in connection with one-dimensional Schr\"odinger and Jacobi
operators. 
\end{abstract}

\maketitle

{\renewcommand{\baselinestretch}{.1} \scriptsize
\tableofcontents}

\section{Introduction}\label{s1}

This Festschrift contribution is devoted to a survey of Barry Simon's principal contributions to the area of inverse spectral theory for one-dimensional Schr\"odinger and Jacobi operators. We decided to put the emphasis on the following five groups of topics:

\smallskip
$\bullet$ The Dirichlet spectral deformation method
\smallskip

A general spectral deformation method applicable to Schr\"odinger, Jacobi, and Dirac-type operators in one dimension, which can be used to insert eigenvalues into spectral gaps of arbitrary background operators but is also an ideal technique to construct isospectral (in fact, unitarily equivalent) sets of operators starting from a given base operator.  

\smallskip
$\bullet$ Renormalized oscillation theory
\smallskip

Renormalized oscillation theory formulated in terms of Wronskians of appropriate solutions rather than solutions themselves, applies, in particular, to energies above the essential spectrum where real-valued solutions exhibit infinitely many zeros and traditional eigenvalue counting methods would naively lead to $\infty - \infty$. While not directly related to inverse spectral methods, we chose to include this topic because of its fundamental importance to the Dirichlet spectral deformation method. 

\smallskip
$\bullet$ The xi function and trace formulas for Schr\"odinger and Jacobi operators
\smallskip

The xi function, that is, essentially, the argument of the diagonal Green's function, which also takes on the role of a particular spectral shift function, is an ideal tool to derive a hierarchy of (higher-order) trace formulas for one-dimensional Schr\"odinger and Jacobi operators. The latter are the natural extensions of the well-known trace formulas for periodic and algebro-geometric finite-band potential coefficients to arbitrary coefficients. The xi function provides a tool for direct and inverse spectral theory.

\smallskip
$\bullet$ Uniqueness theorems in inverse spectral theorem
\smallskip

Starting from the Borg--Marchenko uniqueness theorem, the basic uniqueness result for Schr\"odinger and Jacobi operators in terms of the Weyl--Titchmarsh $m$-coefficient, a number of uniqueness results are discussed. The latter include the Borg-type two-spectra results as well as Hochstadt--Lieberman-type results with mixed prescribed data. In addition to these traditional inverse spectral problems, several new types of inverse spectral problems are addressed.

\smallskip
$\bullet$ Simon's new approach to inverse spectral theory
\smallskip

In some sense, Simon's new approach to inverse spectral theory for half-line problems, based on a particular representation of the Weyl--Titchmarsh $m$-function as a finite Laplace transform with control about the error term, can be viewed as a continuum analog of the continued fraction approach (based on the Riccati equation) to the inverse spectral problem for semi-infinite Jacobi matrices (the actual details, however, differ considerably). Among a variety of spectral-theoretic results, this leads to a formulation of the half-line inverse spectral problem alternative to that of Gel'fand and Levitan. In addition, it leads to a fundamental new result, the local Borg--Marchenko uniqueness theorem.

\medskip 

Each individual section focuses on a particular paper or on a group of papers to be surveyed, representing the five items just discussed. Since this survey is fairly long, it was our intention to write each section in such a manner that it can be read independently.

Only self-adjoint Schr\"odinger and Jacobi operators are considered in this survey. In particular, all potential coefficients $V$ and Jacobi matrix coefficients $a$ and $b$ are assumed to be real-valued throughout this survey (although, we occasionally remind the reader of this assumption).

To be sure, this is not a survey of the state of the art of inverse spectral 
theory for one-dimensional Schr\"odinger and Jacobi operators. Rather, it exclusively focuses on Barry Simon's contributions to and influence exerted on the field. Especially, the bibliography, although quite long, is far from complete and only reflects the particular purpose of this survey. The references included are typically of the following two kinds: First, background references that were used by Barry Simon and his coworkers in writing a particular paper. Such references are distributed throughout the particular survey of one of his papers. Second, at the end of each such survey we refer to more recent references which complement the results of the particular paper in question.   

\bigskip

It was 23 years ago in April of 1983 that Barry and I first met in person at Caltech and started our collaboration. Barry became a mentor and then a friend, and it is fair to say he has had a profound influence on my career since that time. Working with Barry has been exciting and most rewarding for me. Happy Birthday, Barry, and many more such anniversaries!

\section{The Dirichlet Spectral Deformation Method}
\label{s2}

In this section we describe some of the principal results of the paper: \\
\cite{GST96} F.\ Gesztesy, B.\ Simon, and G.\ Teschl, {\it Spectral
deformations of one-dimension- al Schr\"odinger operators}, J. Analyse
Math. {\bf 70}, 267--324 (1996).

\medskip\smallskip

Spectral deformations of Schr\"odinger operators in $L^2 (\bbR)$,
isospectral and certain
classes of non-isospectral ones, have attracted a lot of interest over the
past three decades
due to their prominent role in connection with a variety of topics,
including the
Korteweg-de Vries (KdV) hierarchy, inverse spectral problems,
supersymmetric quantum
mechanical models, level comparison theorems, etc. In fact, the
construction of $N$-soliton
solutions of the KdV hierarchy (and more generally, the construction of
solitons relative
to reflectionless backgrounds) is a typical example of a non-isospectral
deformation of
$H=-\frac{d^2}{dx^2}$ in $L^2 (\bbR)$ since the resulting deformation
$\tilde H = -\frac{d^2}{dx^2} + \tilde V$ acquires an additional point spectrum
$\{\lambda_1, \dots,\lambda_N\}\subset (-\infty, 0)$, $N\in\bbN$, such that
$$
\sigma(\tilde H) =\sigma (H)\cup \{\lambda_1, \dots, \lambda_N\}
$$
($\sigma(\,\cdot\,)$ abbreviating the spectrum). In the $N$-soliton context 
(ignoring the KdV time parameter for simplicity), $\wti V$ is of the explicit form 
\begin{equation}
\wti V(x) = - 2 \f{d^2}{dx^2} \ln[W(\Psi_1(x),\dots,\Psi_N(x))],  \quad x\in\bbR,  \lb{2.2.2}
\end{equation}
where $W(f_1,\dots,f_N)$ denotes the Wronskian of $f_1,\dots,f_N$ and 
the functions $\Psi_j$, $j=1,\dots,N$, are given by 
\begin{align*}
& \Psi_j(x)=(-1)^{j+1}e^{-\kappa_j x} +\alpha_j e^{\kappa_j x}, \quad x\in\bbR, \\
& 0<\kappa_1<\cdots<\kappa_N, \quad \alpha_j >0, \; j=1,\dots,N.   
\end{align*}
The Wronski-type formula in \eqref{2.2.2} is typical also for general background 
potentials and typical for the Crum--Darboux-type commutation approach \cite{Cr55}, 
\cite{DT79}  (cf.\ \cite{GSS91} and the references therein for general backgrounds) which underlies all standard spectral deformation methods for one-dimensional Schr\"odinger operators such as single and double commutation, and the Dirichlet deformation method presented in this section. 

On the other hand, the solution of the inverse periodic problem and the corresponding solution of the algebro-geometric quasi-periodic finite-band inverse problem for 
the KdV hierarchy (and certain almost-periodic limiting situations 
thereof) are intimately connected with isospectral (in fact, unitary)
deformations of a given base (background) operator
$H=-\frac{d^2}{dx^2}+V$.  Although not a complete bibliography on
applications of spectral deformations in mathematical physics,
the interested reader may consult \cite{Ba86}, \cite{BS89}, \cite{BGGSS87},
\cite{BF84}, \cite{Cr55}, \cite{De78}, \cite{DT79}, 
\cite[Sect.\ 4.3]{EK82a}, \cite{EK82}, 
\cite{EF85}, \cite{FIT87}, \cite{Fi89}, \cite{FM76}, 
\cite{GGKM74}, \cite{Ge91}, \cite{Ge92}, \cite{GS95}, 
\cite{Ge93}, \cite{GH00}, \cite[App.\ G]{GH03}, \cite{GSS91}, 
\cite{GST96}, \cite{GT96}, \cite{GW93}, \cite[Ch.\ 2, App.\ A]{GM97},
\cite{Iw87}, \cite{KM56}, \cite[Sect.\ 6.6]{Le87}, 
\cite{Mc85}--\cite{MM75}, \cite[Chs.\ 3, 4]{PT87}, \cite{RT88},
\cite{RS85}, \cite{Sc78}, \cite{Sc03}, \cite{Sh86}, \cite{Sh87}, and the
numerous references cited therein.

The main motivation for writing \cite{GST96} originated in our
interest in inverse spectral
problems. As pointed out later (see Remarks \ref{r2.4.5}--\ref{r2.4.8}),
spectral deformation methods can provide crucial insights into the
isospectral class of a given base potential $V$\!, and in some
cases can even determine the whole isospectral class 
of such potentials. A particularly interesting open problem in inverse
spectral theory concerns
the characterization of the isospectral class of potentials $V$ with purely
discrete spectra
(e.g., the harmonic oscillator $V(x)=x^2$, cf.\ \cite{Ch03}--\cite{CK05},
\cite{GS04}, \cite{Le88}, \cite{MT82}, \cite{PS05}).

The principal result in \cite{GST96}, reviewed in this section (cf.\
Theorem \ref{t2.4.4}\,$(i)$), provides a complete spectral characterization
of a new method of constructing isospectral (in fact, unitary)
deformations of general Schr\"odinger operators $H=-\frac{d^2}{dx^2}+V$ in
$L^2 (\bbR)$. The technique is connected to Dirichlet data, that is,to
the spectrum of the operator $H^D$ on $L^2 ((-\infty, x_0])\oplus L^2
([x_0, \infty))$ with a Dirichlet boundary condition at $x_0$. The
transformation moves a single eigenvalue of $H^D$ and perhaps flips the
half-line (i.e.,
$(-\infty,x_0)$ to $(x_0,\infty)$, or vice versa) to which the Dirichlet
eigenvalue belongs. On the remainder of the spectrum, the transformation
is realized by a unitary operator.

To describe the Dirichlet deformation method (DDM) as developed in
\cite{GST96} in some detail, we suppose that
$
V\in L^1_{\rm loc} (\bbR)
$
is real valued and introduce the differential expression
$\tau=-\frac{d^2}{dx^2}+V(x)$, $x\in\bbR$. Assuming $\tau$ to be in the
limit point case at $\pm\infty$ (for the general case we refer to
\cite{GST96}) one defines the self-adjoint base
(i.e., background) operator $H$ in $L^2 (\bbR)$ by
\begin{equation}
Hf =\tau f, \quad
f\in\dom(H) = \{g\in L^2 (\bbR) \,|\, g,g'\in
AC_{\rm loc}(\bbR);
\tau g\in L^2(\bbR)\}.  \lb{2.2.3}
\end{equation}
Here $W(f,g)(x)=f(x)g'(x)-f'(x)g(x)$ denotes the Wronskian of $f,g\in
AC_{\rm loc}
(\bbR)$ (the set of locally absolutely continuous functions on $\bbR$).
Given $H$ and a fixed reference point $x_0 \in\bbR$, we introduce
the associated Dirichlet operator $H^D_{x_0}$ in $L^2 (\bbR)$ by
\begin{align*}
H^D_{x_0} f = \tau f, \quad
f\in \dom(H^D_{x_0}) &= \{g\in L^2 (\bbR) \,|\, g\in
AC_{\rm loc}(\bbR), g'\in
AC_{\rm loc} (\bbR\backslash \{x_0\}); \\
& \quad \;\;\, \lim_{\epsilon\downarrow 0}
g(x_0 \pm\epsilon)=0; \, \tau g\in L^2 (\bbR)\}.
\end{align*}
Clearly, $H^D_{x_0}$ decomposes into
$
H^D_{x_0} = H^D_{-,x_0} \oplus H^D_{+,x_0}
$
with respect to the orthogonal decomposition
$
L^2 (\bbR) =L^2 ((-\infty, x_0]) \oplus L^2 ([x_0, \infty)).
$
Moreover, for any
$\mu\in\sigma_d (H^D_{x_0})\backslash \sigma(H)$  ($\sigma_d
(\,\cdot\,)$, the discrete spectrum,
$\sigma(\,\cdot\,)$
and $\sigma_{\rm ess}(\,\cdot\,)$, the spectrum and essential
spectrum, respectively),
we introduce the Dirichlet datum
\begin{equation}
(\mu,\sigma) \in\{\sigma_d (H^D_{x_0})\backslash\sigma_d (H)\} \times
\{-,+\},
\lb{2.2.7}
\end{equation}
which identifies $\mu$ as a discrete Dirichlet eigenvalue on the interval
$(x_0, \sigma
\infty)$, that is, $\mu\in\sigma_d (H^D_{\sigma,x_0})$, $\sigma\in
\{-,+\}$ (but excludes it
from being simultaneously a Dirichlet eigenvalue on $(x_0,
-\sigma\infty)$).

Next, we pick a fixed spectral gap $(E_0, E_1)$ of $H$, the endpoints of
which (without loss of generality) belong to the spectrum of $H$,
$$
(E_0, E_1) \subset \bbR\backslash\sigma(H), \quad
E_0, E_1 \in\sigma(H)
$$
and choose a discrete eigenvalue $\mu$ of $H^D_{x_0}$ in the closure of
that spectral gap,
\begin{equation}
\mu\in\sigma_d (H^D_{x_0})\cap [E_0, E_1] \lb{2.2.9}
\end{equation}
(we note there is at most one such $\mu$ since $(H^D_{x_0} -z)^{-1}$ is a
rank-one perturbation
of $(H-z)^{-1}$). According to \eqref{2.2.7}, this either gives rise to a
Dirichlet datum
\begin{equation}
(\mu,\sigma)\in (E_0, E_1) \times \{-,+\}, \lb{2.2.10}
\end{equation}
or else to a discrete eigenvalue of $H^D_{-,x_0}$ and $H^D_{+,x_0}$, that
is,
\begin{equation}
\mu\in\{E_0, E_1\}\cap\sigma_d (H)\cap\sigma_d(H^D_{-,x_0})\cap
\sigma_d (H^D_{+,x_0}) \lb{2.2.11}
\end{equation}
since the eigenfunction of $H$ associated with $\mu$ has a zero at $x_0$.
In particular,
since $(H^D_{x_0} -z)^{-1}$ is a rank-one perturbation of $(H-z)^{-1}$,
one infers
$$
\sigma_{\rm ess} (H^D_{x_0}) = \sigma_{\rm ess}(H),
$$
and thus, $\mu\in \{E_0, E_1\}\cap\sigma_{\rm ess}(H)$ is
excluded by assumption \eqref{2.2.9}. Hence, the case distinctions
\eqref{2.2.10} and \eqref{2.2.11} are exhaustive.

In addition to $\mu$ as in \eqref{2.2.9}--\eqref{2.2.11}, we also need to
introduce
$\tilde\mu\in
[E_0, E_1]$ and $\tilde\sigma\in \{-,+\}$ as follows: Either
$$
(\tilde\mu, \sigma)\in (E_0, E_1) \times \{-,+\},
$$
or else
$$
\tilde\mu\in \{E_0, E_1\}\cap\sigma_d (H).
$$
Given $H$, one introduces Weyl--Titchmarsh-type solutions $\psi_\pm (z,x)$
of $(\tau-z)\psi(z)=0$ by
\begin{align*}
& \psi_\pm (z,\,\cdot\,)\in L^2 ((R,\pm\infty)), \quad R\in \bbR, \\
& \lim_{x\to\pm\infty} W(\psi_\pm (z), g)(x)=0 \quad \text{for all
$g\in\dom(H)$}.
\end{align*}
If $\psi_\pm (z,x)$ exist, they are unique up to constant multiples. In
particular, $\psi_\pm
(z,x)$ exist for $z\in\bbC\backslash\sigma_{\rm ess}(H)$ and we
can (and
will) assume them to be holomorphic with respect to 
$z\in\bbC\backslash\sigma(H)$ and
real-valued for $z\in\bbR$ (cf.\ the discussion in connection with
\eqref{3.1.1}). 

Given $\psi_\sigma (\mu, x)$ and $\psi_{-\tilde\sigma}(\tilde\mu, x)$, one
defines
$$
W_{(\tilde\mu, \tilde\sigma)} (x) =
\begin{cases} (\tilde\mu -\mu)^{-1} W(\psi_\sigma (\mu),
\psi_{-\tilde\sigma}(\tilde\mu))(x),
& \mu,\tilde\mu \in [E_0, E_1], \  \tilde\mu \neq \mu, \\
-\sigma \int^x_{\sigma\infty} dx' \psi_\sigma (\mu, x')^2, &
(\tilde\mu, \tilde\sigma)
=(\mu, -\sigma), \ \mu\in (E_0, E_1), \
\end{cases}
$$
and the associated Dirichlet deformation
\begin{align}
& \tilde\tau_{(\tilde\mu, \tilde\sigma)} = -\frac{d^2}{dx^2} +\tilde
V_{(\tilde\mu,
\tilde\sigma)}(x), \no \\
& \tilde V_{(\tilde\mu, \tilde\sigma)}(x) = V(x)-2\{\ln [W_{(\tilde\mu,
\tilde\sigma)}
(x)]\}'', \quad x\in\bbR, \lb{2.2.20}  \\
& \mu, \tilde\mu \in [E_0, E_1], \  \mu\neq \tilde\mu
\text{ or } (\tilde\mu, \tilde\sigma) = (\mu, -\sigma), \
\mu\in (E_0, E_1).  \no
\end{align}
As discussed in Section \ref{s3}, $W_{(\tilde\mu, \tilde\sigma)}(x)\neq 0$,
$x\in\bbR$, and hence \eqref{2.2.20} yields a well-defined potential
$
\tilde V_{(\tilde\mu, \tilde\sigma)}\in L^1_{\rm loc} (\bbR).
$

In the remaining cases
$(\tilde\mu,
\tilde\sigma) =
(\mu, \sigma)$, $\mu\in [E_0, E_1]$, and $\mu=\tilde\mu \in\{E_0, E_1 \}
\cap \sigma_d (H)$, we define
$
\tilde V_{(\tilde\mu, \tilde\sigma)} =V
$
which represents the trivial deformation of $V$ (i.e., none at
all), and for notational simplicity these trivial cases are excluded in the
remainder of this section. For obvious
reasons we will allude to \eqref{2.2.20} as the Dirichlet
deformation method in the following.

If $\tilde\mu \in\sigma_d (H)$, then $\psi_- (\tilde\mu)=c\psi_+
(\tilde\mu)$ for some
$c\in\bbR\backslash \{0\}$, showing that $W_{(\tilde\mu,
\tilde\sigma)}(x)$ and hence,
$V_{(\tilde\mu, \tilde\sigma)}(x)$ in \eqref{2.2.20} becomes independent
of $\sigma$
or $\tilde\sigma$. In this case we shall occasionally use a more
appropriate notation
and write $\tilde V_{\tilde\mu}$ and $\tilde\tau_{\tilde\mu}$ (instead of
$\tilde V_{(\tilde\mu, \tilde\sigma)}$ and $\tilde\tau_{(\tilde\mu,
\tilde\sigma)}$).  

For later reference, we now summarize our basic assumptions on $V$,
$\mu$, and
$\tilde\mu$ in the following hypothesis.

\begin{hypothesis} \lb{h2.2.3}  Suppose $V\in L^1_{\rm loc}(\bbR)$ to
be real-valued. In addition, we assume
\begin{align*}
&(E_0, E_1)\subset\bbR\backslash\sigma(H), \quad E_0, E_1 \in\sigma(H), \\
&\mu\in\sigma_d (H^D_{x_0}), \quad (\mu,\sigma)\in (E_0, E_1)\times
\{-,+\}
\text{ or } \mu\in\{E_0, E_1\}\cap \sigma_d (H), \\
&(\tilde\mu,\tilde\sigma)\in (E_0, E_1)\times\{-,+\} \text{ or }
\tilde\mu \in \{E_0, E_1\}\cap\sigma_d (H), \\
&\mu,\tilde\mu \in [E_0, E_1],  \quad \mu\neq\tilde\mu \text{ or }
(\tilde\mu, \tilde\sigma) = (\mu,-\sigma), \quad \mu\in (E_0, E_1).
\end{align*}
\end{hypothesis}

Next, introducing the following solutions of
$(\tilde\tau_{(\tilde\mu,\tilde\sigma)}-z)
\tilde\psi (z)=0$,
\begin{align*}
&\tilde\psi_{-\sigma}(\mu,x) =\psi_{-\tilde\sigma}(\tilde\mu,x) \big/
W_{(\tilde\mu, \tilde\sigma)}(x), \\
&\tilde\psi_{\tilde\sigma}(\tilde\mu, x) = \psi_\sigma (\mu,x) \big/
W_{(\tilde\mu,\tilde\sigma)}(x), \quad \tilde\psi_{\tilde\sigma}(\tilde\mu,
x_0)=0,
\end{align*}
one 
infers
$$
(\tilde\tau_{(\tilde\mu,\tilde\sigma)}\tilde\psi_{-\sigma}(\mu))(x)=\mu
\tilde\psi_{-\sigma}(\mu,x), 
\quad 
(\tilde\tau_{(\tilde\mu,\tilde\sigma)}
\tilde\psi_{\tilde\sigma}(\tilde\mu))(x)=\tilde\mu 
\tilde\psi_{\tilde\sigma}
(\tilde\mu, x).
$$

The Dirichlet deformation operator $\tilde H_{(\tilde\mu,\tilde\sigma)}$
associated with
$\tilde\tau_{(\tilde\mu,\tilde\sigma)}$ in \eqref{2.2.20} is then defined
as  follows:
\begin{align}
& \tilde H_{(\tilde\mu,\tilde\sigma)}f
=\tilde\tau_{(\tilde\mu,\tilde\sigma)} f, \quad
f\in\dom(\tilde H_{(\tilde\mu,\tilde\sigma)}) = \{g\in L^2(\bbR) \,|\, 
g,g'\in AC_{\rm loc}(\bbR); \no  \\
& \hspace*{4.1cm}  g \text{ satisfies the b.c. in \eqref{2.2.33};} \,
\tilde\tau_{(\tilde\mu,\tilde\sigma)} g\in L^2 (\bbR)\}.  \lb{2.2.32}
\end{align}
The boundary conditions (b.c.'s) alluded to in \eqref{2.2.32} are chosen as
follows:
\begin{align}
\begin{split}
\lim_{x\to\tilde\sigma \infty}
W(\tilde\psi_{\tilde\sigma}(\tilde\mu), g)(x) &= 0
\text{ if $\tilde\tau_{(\tilde\mu,\tilde\sigma)}$ is l.c.~at $\tilde\sigma
\infty$}, \\
\lim_{x\to -\tilde\sigma\infty}
W(\tilde\psi_{-\sigma}(\mu), g)(x) &=0
\text{ if $\tilde\tau_{(\tilde\mu,\tilde\sigma)}$ is l.c.~at
$-\tilde\sigma\infty$}.   \lb{2.2.33}
\end{split}
\end{align}
Here we abbreviate the limit point and limit circle cases by l.p. and
l.c., respectively. As usual, the boundary condition at $\omega\infty$ in
\eqref{2.2.32} is omitted if
$\tilde\tau_{(\tilde\mu,\tilde\sigma)}$ is l.p.~at $\omega\infty$,
$\omega\in\{-,+\}$.

For future reference we note that in analogy to the Dirichlet operators
$H^D_{x_0}$, $H^D_{\pm,x_0}$ introduced in connection with the operator
$H$, one can also introduce the corresponding Dirichlet operators
$\tilde H^D_{(\tilde\mu,\tilde\sigma),x_0}$,
$\tilde H^D_{(\tilde\mu,\tilde\sigma),\pm,x_0}$ associated
with $\tilde H_{(\tilde\mu,\tilde\sigma)}$.

Next, we turn to the Weyl--Titchmarsh $m$-functions for the
Dirichlet deformation
operator $\tilde H_{(\tilde\mu,\tilde\sigma)}$ and relate them to those of
$H$.

Let $\phi(z,x),\theta(z,x)$ be the standard fundamental system of
solutions of
$(\tau -z)\psi(z)=0$, $z\in\bbC$ defined by
$\phi(z,x_0)=\theta'(z,x_0)=0$, $\phi'(z,x_0)=\theta(z,x_0)=1$,
$z\in\bbC$, and denote by
$\tilde\theta_{(\tilde\mu,\tilde\sigma)}(z,x),
\tilde\phi_{(\tilde\mu,\tilde\sigma)}(z,x)$
the analogously normalized fundamental system of solutions of
$(\tilde\tau_{(\tilde\mu,\tilde\sigma)}-z)
\tilde\psi(z)=0$, $z\in\bbC$, at $x_0$. One then has
$$
m_\sigma (z,x_0)=\psi'_\sigma (z,x_0)/\psi_\sigma (z,x_0) , \quad
\sigma\in\{-,+\}, \
z\in\bbC\backslash\bbR,
$$
where $m_\sigma(\cdot,x_0)$ denotes the Weyl--Titchmarsh $m$-function of
$H$ with respect to the half-line $(x_0, \sigma\infty)$,
$\sigma\in\{-,+\}$. For the corresponding half-line Weyl--Titchmarsh
$m$-functions of $\tilde H_{(\tilde\mu, \tilde\sigma)}$ in terms of those
of $H$ one then obtains the following result.

\begin{theorem} \lb{t2.3.2} Assume Hypothesis \ref{h2.2.3} and  
$z\in\bbC\backslash\bbR$. Let $H$ and $\tilde H_{(\tilde\mu,
\tilde\sigma)}$ be given by \eqref{2.2.3} and \eqref{2.2.32},
respectively, and denote by $m_\pm$
and $\tilde m_{(\tilde\mu, \tilde\sigma), \pm}$ the corresponding
$m$-functions
associated with the half-lines $(x_0, \pm\infty)$. Then,
\begin{align*}
& \tilde m_{(\tilde\mu, \tilde\sigma), \pm}(z,x_0) =
\frac{z-\mu}{z-\tilde\mu}\, m_\pm (z,x_0)
-\frac{\tilde\mu-\mu}{z-\tilde\mu}\, m_{-\tilde\sigma}(\tilde\mu,x_0),
\quad
\tilde\mu\neq\mu, \\
& \tilde m_{(\tilde\mu, \tilde\sigma), \pm}(z,x_0) = m_\pm (z,x_0) -
\biggl(\,\int^{x_0}_{\sigma
\infty}dx\, \phi(\mu,x)^2\biggr)^{-1} \frac{1}{z-\mu}, \quad (\tilde\mu,
\tilde\sigma) =
(\mu, -\sigma).
\end{align*}
\end{theorem}

Given the fundamental relation between $\tilde m_{(\tilde\mu,
\tilde\sigma),\pm}$
and $m_\pm$ in Theorem \ref{t2.3.2}, one can now readily derive the
ensuing relation between
the corresponding spectral functions $\tilde\rho_{(\tilde\mu,
\tilde\sigma),\pm}$
and $\rho_\pm$ associated with the half-line Dirichlet operators
$\tilde H^D_{(\tilde\mu, \tilde\sigma),\pm,x_0}$ and $H^D_{\pm,x_0}$. For
this and a complete spectral characterization of
$\tilde H^D_{(\tilde\mu,\tilde\sigma), \pm,x_0}$ in terms of
$H^D_{\pm,x_0}$ we refer to \cite{GST96}.

Next we turn to the principal results of \cite{GST96} including explicit
computations of the
Weyl--Titchmarsh and spectral matrices of $\tilde H_{(\tilde\mu,
\tilde\sigma)}$ in terms
of those of $H$ and a complete spectral
characterization of
$\tilde H_{(\tilde\mu, \tilde\sigma)}$ and $\tilde H^D_{(\tilde\mu,
  \tilde\sigma),x_0}$
in terms of $H$ and $H^D_{x_0}$.

We start with the Weyl--Titchmarsh matrices for $H$ and
$\tilde H_{(\tilde\mu,
\tilde\sigma)}$. To fix notation, we introduce the Weyl--Titchmarsh
$M$-matrix in
$\bbC^2$ associated with $H$ by
\begin{align*}
&M(z,x_0) = (M_{p, q}(z,x_0))_{1\leq p, q\leq 2} 
= [m_- (z,x_0) -m_+ (z,x_0)]^{-1} \\
& \quad \times \begin{pmatrix} m_- (z,x_0) m_+(z,x_0)
& [m_- (z,x_0)+m_+ (z,x_0)]/2 \\
[m_- (z,x_0)+m_+ (z,x_0)]/2 & 1 \end{pmatrix} , \quad
z\in\bbC\backslash\bbR,
\end{align*}
and similarly $\tilde M_{(\tilde\mu, \tilde\sigma)}$ in connection with
$\tilde H_{(\tilde\mu,\tilde\sigma)}$.
An application of Theorem \ref{t2.3.2} then yields

\begin{theorem} \lb{t2.4.1} Assume Hypothesis \ref{h2.2.3} and
$z\in\bbC\backslash\bbR$.Let  $H$ and $\tilde H_{(\tilde\mu,
\tilde\sigma)}$ be given by {\rm{(2.3)}} and {\rm{(2.32)}},
respectively. Then the corresponding  Weyl--Titchmarsh-matrices $M$ and
$\tilde M_{(\tilde\mu, \tilde\sigma)}$ are related by
\begin{align*}
\tilde M_{(\tilde\mu, \tilde\sigma), 1,1}(z,x_0) &=
\frac{z-\mu}{z-\tilde\mu}\, M_{1,1}(z,x_0)
- 2\frac{\tilde\mu-\mu}{z-\tilde\mu}\,
m_{-\tilde\sigma}(\tilde\mu,x_0)M_{1,2}(z,x_0) \\
&\quad + \frac {(\tilde\mu-\mu)^2}{(z-\mu)(z-\tilde\mu)}\,
m_{-\tilde\sigma} (\tilde\mu,x_0)^2 M_{2,2}(z,x_0),    \\
\tilde M_{(\tilde\mu, \tilde\sigma), 1,2}(z,x_0) &= M_{1,2}(z,x_0))
-\frac{\tilde\mu-\mu}{z-\mu}\,
m_{-\tilde\sigma}(\tilde\mu,x_0) M_{2,2}(z,x_0),   \\
\tilde M_{(\tilde\mu, \tilde\sigma),2,2}(z,x_0) &=
\frac{z-\tilde\mu}{z-\mu}\, M_{2,2}(z,x_0), \quad \tilde\mu\neq\mu.
\end{align*}
\end{theorem}

Given the basic connection between $\tilde M_{(\tilde\mu,
\tilde\sigma)}$ and $M$
in Theorem \ref{t2.4.1}, one can now proceed to derive the analogous
relations between the spectral
matrices $\tilde\rho_{(\tilde\mu, \tilde\sigma)}$ and
$\rho$ associated with
$\tilde H_{(\tilde\mu, \tilde\sigma)}$ and $H$, respectively (cf.\
\cite{GST96} for details).

The principal spectral deformation result of \cite{GST96} then reads as
follows.

\begin{theorem} \lb{t2.4.4} Assume Hypothesis \ref{h2.2.3}. Then, \\
  $(i)$ Suppose $\mu,\tilde\mu\in (E_0, E_1)$. Then $\tilde
H_{(\tilde\mu, \tilde\sigma)}$
and $H$ are unitarily equivalent. Moreover, $\tilde H^D_{(\tilde\mu,
\tilde\sigma),x_0}$ and
$H^D_{x_0}$, restricted to the orthogonal complements of the
one-dimensional eigenspaces
corresponding to $\tilde\mu$ and $\mu$, are unitarily equivalent. \\
$(ii)$ Assume $\mu\in\{E_0, E_1\}\cap\sigma_d (H)$, $\tilde\mu\in
(E_0, E_1)$. Then,
$$
\sigma_{(p)}(\tilde H_{(\tilde\mu, \tilde\sigma)})
=\sigma_{(p)}(H)\backslash
\{\mu\}, \quad
\sigma_{(p)}( \tilde H^D_{(\tilde\mu, \tilde\sigma),x_0})
= \{\sigma_{(p)}(H^D_{x_0})
\backslash \{\mu\}\}\cup\{\tilde\mu\}.  
$$
$(iii)$ Suppose $\mu\in 
(E_0, E_1)$, $\tilde\mu\in \{E_0, E_1\}\cap
\sigma_d (H)$. 
Then,
$$
\sigma_{(p)}(\tilde H_{\tilde\mu}) = 
\sigma_{(p)}(H)\backslash
\{\tilde\mu\}, \quad 
\sigma_{(p)}(\tilde 
H^D_{\tilde\mu,x_0}) =\sigma_{(p)}(H^D_{x_0})\backslash\{\mu\},
\quad 
\tilde\mu\notin\sigma(\tilde H^D_{\tilde\mu,x_0}).  
$$
$(iv)$ Assume 
$\mu, \tilde\mu\in \{E_0, E_1\}\cap
\sigma_d (H)$, $\mu 
\neq\tilde\mu$. Then,
$$
\sigma_{(p)}(\tilde H_{\tilde\mu}) = 
\sigma_{(p)}(H)\backslash\{E_0, E_1\},
\quad 
\sigma_{(p)}(\tilde 
H^D_{\tilde\mu,x_0}) = \sigma_{(p)}(H^D_{x_0})\backslash \{\mu\},
\quad 
\tilde\mu\notin\sigma (\tilde H^D_{\tilde\mu,x_0}).  
$$
In cases 
$(ii)$--$(iv)$, the corresponding pairs of operators,
restricted to 
the obvious orthogonal complements of the eigenspaces
corresponding 
to $\mu$ and/or $\tilde\mu$, are unitarily
equivalent. In 
particular,
$$
\sigma_{\rm ess, ac, sc}(\tilde H_{(\tilde\mu, 
\tilde\sigma)}) = \sigma_{\rm ess, ac, sc}(\tilde 
H^D_{(\tilde\mu,\tilde\sigma),x_0}) = \sigma_{\rm ess, ac, 
sc}(H^D_{x_0}) = \sigma_{\rm ess, ac, sc}(H). 
$$
\end{theorem}

\begin{remark} 
\lb{r2.4.5} $(i)$ Perhaps the most important consequence of
Theorem 
\ref{t2.4.4} (i),
from an inverse spectral point of view, is the fact 
that {\it{any}}
finite
number of
deformations of Dirichlet data 
within spectral gaps of $V$ yields a
potential $\tilde V$
in the 
isospectral class of $V$\!. No further constraints on 
$(\mu_j,
\sigma_j),
(\tilde\mu_j, \tilde\sigma_j)$, other than 
$(\mu_j, \sigma_j),
(\tilde\mu_j, \tilde\sigma_j)
\in (E_{j-1}, E_j) 
\times \{-,+\}$, $(E_{j-1}, E_j)\in\bbR\backslash
\sigma (H)$,
$j=1, 
\dots, N$, $N\in\bbN$, are involved. \\
$(ii)$ The isospectral 
property $(i)$ in Theorem
\ref{t2.4.4}, in the special case
of 
periodic potentials, was first proved by Finkel, Isaacson,
and 
Trubowitz \cite{FIT87}. Further results can be found in Buys 
and
Finkel \cite{BF84} and Iwasaki \cite{Iw87} (see also 
\cite{De78},
\cite{DT79}, \cite{Mc85},
\cite{Mc86}, \cite{Mc87}, 
\cite{Mc92}). Similar constructions for
Schr\"odinger operators on a 
compact interval can be found in P\"oschel and
Trubowitz \cite[Chs.\ 
3, 4]{PT87} and Ralston and Trubowitz \cite{RT88}. \\
$(iii)$ Let 
$\mu\in (E_0, E_1)$. Then the (isospectral) Dirichlet
deformation 
$(\mu,\sigma)
\to (\mu, -\sigma)$ is precisely the isospectral case 
of the double commutation method (cf.\ \cite{Ge93}, \cite{GH00}, 
\cite[App.\ B]{GST96},
\cite{GT96}). It simply flips the Dirichlet 
eigenvalue
$\mu$ on the half-line $(x_0, \sigma\infty)$ to the other 
half-line $(x_0,
-\sigma\infty)$.
In the special case where $V(x)$ is periodic, this procedure was first 
used by McKean and van Moerbeke \cite{MM75}. \\
$(iv)$ The topology of these Dirichlet data strongly depends on the nature of the 
endpoints $E_0, E_1$ of a particular spectral gap. For instance, in cases like the periodic one, different spectral gaps are separated by intervals of absolutely continuous spectrum and the two intervals $[E_0,E_1]$ together with $\sigma\in \{-,+\}$ can be identified with a circle (upon identifying the two intervals as two rims of a cut). Globally this then leads to a product of circles, that is, a torus. In particular, the Dirichlet eigenvalues in different spectral gaps can be prescribed independently of each other. The situation is entirely different if an endpoint, say $E_1$, belongs to the discrete spectrum of $H$. In this case there are two neighboring spectral gaps $(E_0,E_1)$ and 
$(E_1,E_2)$ and the two Dirichlet eigenvalues $\mu_j \in (E_{j-1},E_j)$, $j=1,2$, are not independent of each other. In fact, if one of $\mu_1$ or $\mu_2$ approaches $E_1$, then necessarily so does the other. The topology is then not a product of circles. For instance, a closer analysis of the case of $N$-soliton potentials in \eqref{2.2.2}  then illustrates that the appropriate coordinates parametrizing the $N$-soliton isospectral class are $\alpha_j \in (0,\infty)$ (compare also with positive norming constants), which results globally in a product of open half-lines.
\end{remark}

\begin{remark} \lb{r2.4.7} In certain cases where the base 
(background) potential $V$ is reflectionless (see, e.g., \cite{GY06}) 
and $H$ is bounded from below and has no singularly continuous
spectrum, the isospectral class $\text{Iso}(V)$ of $V$ (the set of all reflectionless 
$\tilde V$'s such that $\sigma(\tilde H)=\sigma(H)$) is completely
characterized by the distribution of Dirichlet
(initial) data $(\mu_{j+1}(x_0), \sigma_{j+1}(x_0))\in [E_j, E_{j+1}]
\times \{-,+\}$, $j\in J$, in nontrivial spectral gaps of $H$. Here
$x_0\in\bbR$ is a fixed reference point and $J=\{0,1,\dots, N-1\}$,
$N\in\bbN$, or $j\in J=\bbN_0$  ($=\bbN\cup\{0\}$) parametrizes these
nontrivial spectral gaps $(E_j, E_{j+1})$ of $H$ (the trivial one being
$(-\infty, \inf \sigma(H))$). Prime examples of this type are periodic
potentials, algebro-geometric quasi-periodic finite-band potentials (cf.\ 
\cite[Ch.\ 3]{BBEIM94}, \cite[Ch.\ 1]{GH03}, \cite[Chs.\ 8--12]{Le87}, 
\cite[Ch.\ 4]{Ma86}, \cite[Ch.\ II]{NMPZ84}), and certain limiting cases
thereof (e.g., soliton potentials). In these cases, an iteration of the
Dirichlet deformation method, in the sense that 
$(\mu_{j+1}(x_0), \sigma_{j+1} (x_0))\to (\tilde\mu_{j+1}(x_0),
\tilde\sigma_{j+1}(x_0))$ within $[E_j, E_{j+1}]\times \{-, +\}$ for each
$j\in J$, independently of each other, yields an explicit realization of
the underlying isospectral class $\text{Iso}(V)$ of reflectionless potentials with base
$V$\!.  In the periodic case, this was first proved by Finkel,
Isaacson, and Trubowitz \cite{FIT87} (see also \cite{BF84}, \cite{Iw87}).
More precisely, the inclusion of limiting cases 
$\mu_{j+1}(x_0)\in \{E _j, E_{j+1}\} \cap \sigma_{\rm ess}(H)$ requires a
special argument (since it is excluded by Hypothesis \ref{h2.2.3})
but this can be provided in the special cases at hand.
\end{remark}


\begin{remark} \lb{r2.4.8} Another case of primary interest 
concerns
potentials
$V$ with
purely discrete spectra bounded from 
below, that is,
$$
 \sigma(H)=\sigma_d(H) = \{E_j\}_{j\in\bbN_0}, 
\quad -\infty< E_0,  \
E_j <E _{j+1}, \  j\in\bbN_0, \quad 
\sigma_{\rm ess}(H)=\emptyset.
$$
(For simplicity, one may think in terms of the harmonic oscillator
$V(x)=x^2$, \cite{Ch03}--\cite{CK05}, \cite{GS04}, \cite{Le88}, \cite{MT82},
\cite{PS05}.) In this case, either
$$
(\mu_{j+1}(x_0), \sigma_{j+1}(x_0)) \in (E_j, E_{j+1}) \times \{-, +\}
\, \text{ or } \,
\mu_{j+1}(x_0) =E _{j+1} =\mu_{j+2} (x_0),
$$
that is, Dirichlet eigenvalues necessarily meet in pairs whenever they hit
an eigenvalue
of $H$. The following trace formula for $V$ in terms of $\sigma(H)=
\{E _j\}_{j\in\bbN_0}$ and $\sigma (H^D_x)=\{\mu_j (x)\}_{j\in\bbN}$
(with $H^D_y$ the Dirichlet operator associated with
$\tau=-\frac{d^2}{dx^2} +V(x)$ and a Dirichlet boundary condition at
$x=y$), proved in \cite{GS96} (cf.\ Section \ref{s4}),
\begin{equation}
V(x)=E _0 +\lim_{\alpha\downarrow 0} \alpha^{-1} \sum^{\infty}_{j=1}
\Big(2e^{-\alpha\mu_j (x)} - e^{-\alpha E _j} - e^{-\alpha E _{j+1}}\Big),
\lb{2.4.23}
\end{equation}
then shows one crucial difference to the periodic-type cases mentioned
previously. Unlike
in the periodic case, though, the initial Dirichlet eigenvalues
$\mu_{j+1}(x_0)$ {\it{cannot}}
be prescribed arbitrarily in the spectral gaps $(E _j, E_{j+1})$ of $H$.
Indeed, the fact
that the Abelian regularization in the trace formula \eqref{2.4.23} for
$V(x)$ converges to a limit
restricts the asymptotic distribution of $\mu_{j+1}(x)\in [E_j, E_{j+1}]$
as $j\to\infty$. For instance, consider $V(x)=x^2-1$, then $E_j=2j$, $j\in\bbN_0$. 
The choice $\mu_j(x_0)=2j-\gamma$, $\gamma\in (0,1)$, then yields for the 
Abelian regularization on the right-hand side of \eqref{2.4.23},
$$
\lim_{\alpha\downarrow 0} \alpha^{-1} \sum_{j=1}^\infty 
\Big(2 e^{-\alpha(2j-\gamma)-e^{-\alpha 2j}-e^{-\alpha(2j-2)}} \Big)
= \lim_{\alpha\downarrow 0} 2[(\gamma -1)+\Oh(\alpha)] \f{1}{1-e^{-2\alpha}} =\infty.
$$
Put differently, our choice of $\mu_j(x_0)=2j-\gamma$, $\gamma\in (0,1)$,  
was not an admissible choice of Dirichlet eigenvalues for the (shifted) harmonic oscillator potential  $V(x)=x^2-1$. 
However, as stressed in Remark \ref{r2.4.5}\,$(i)$, one of the fundamental
consequences of \cite{GST96} 
concerns the fact that there is no such restriction for any finite number
of spectral gaps of $H$. In other words, only the tail end of the Dirichlet
eigenvalues
$\mu_{j+1}(x_0)$ as $j\to\infty$ is restricted (the precise nature of this
restriction being
unknown at this point), any finite number of them can be placed arbitrarily
in the spectral
gaps $(E _j, E_{j+1})$ (with the obvious ``crossing'' restrictions at the
common boundary
$E _{j+1}$ of $(E _j, E_{j+1})$ and $(E _{j+1}, E_{j+2})$). The only other
known restriction to date on Dirichlet initial data $(\mu_j
(x_0), \sigma_j (x_0))$ is that $\sigma_j (x_0) =-$ and $\sigma_j (x_0)=+$
infinitely often, that is,
both half-lines $(-\infty, x_0)$ and $(x_0, \infty)$ support (naturally)
infinitely many Dirichlet eigenvalues.
\end{remark}

For various extensions of the results presented, including a careful
discussion of limit point/limit circle properties of the Dirichlet
deformation operators, iterations of DDM to insert finitely many
eigenvalues in spectral gaps, applications to reflectionless
Schr\"odinger operators, general
Sturm--Liouville operators in a weighted $L^2$-space, applications to
short-range scattering theory, and a concise summary of single and double
commutation methods, we refer to \cite{GST96}.

\medskip 

{\it More recent references:} An interesting refinement of Theorem \ref{t2.4.4}\,$(i)$, in which a unitary operator relating $\tilde H_{(\tilde\mu, \tilde\sigma)}$ and $H$ is
explicitly characterized, is due to Schmincke \cite{Sc03}. DDM  
for one-dimensional Jacobi and Dirac-type operators has been worked out by Teschl \cite{Te97}, \cite[Ch.\ 11]{Te00} (see also \cite{GT96a}, \cite{Te99}), \cite{Te98a}.

\section{Renormalized Oscillation Theory} \label{s3}

In this section we summarize some of the principal results of
the paper: \\  
\cite{GST96a} F.\ Gesztesy, B.\ Simon, and G.\ Teschl, {\it Zeros of the
Wronskian and renormalized oscillation theory}, Amer. J. Math. {\bf 118},
571--594 (1996).

\medskip\smallskip

For over a hundred and fifty years, oscillation theorems for
second-order differential equations have fascinated mathematicians.
Originating with Sturm's celebrated memoir \cite{St36}, extended in a
variety of ways by B\^ocher \cite{Bo17} and others, a large body of
material has been accumulated since then (thorough treatments can be found,
e.g., in \cite{Co71}, \cite{Kr73}, \cite{Re80}, \cite{Sw68}, and the
references therein). In \cite{GST96a} a new wrinkle to oscillation theory
was added by showing that zeros of Wronskians can be used to count
eigenvalues in situations where a naive use of oscillation theory would
give $\infty -\infty$ (i.e., Wronskians lead to renormalized oscillation
theory). In a nutshell, we will show the following result for general
Sturm--Liouville operators $H$ in $L^2((a,b);rdx)$ with separated boundary
conditions at $a$ and $b$ in this section: If $E_{1,2}\in\bbR$ and if
$u_{1,2}$ solve the differential equation $Hu_j=E_j u_j$, $j=1,2$ and
respectively satisfy the boundary condition on the left/right, then the
dimension of the spectral projection $P_{(E_1, E_2)}(H)$ of $H$ equals the
number of zeros of the Wronskian of $u_1$ and $u_2$.

The main motivation in writing \cite{GST96a} had its origins in attempts
to provide a general construction of isospectral potentials for
one-dimensional Schr\"odinger operators following previous work
by Finkel, Isaacson, and Trubowitz \cite{FIT87} (see also \cite{BF84}) in
the special case of periodic potentials. In fact, in the case of periodic
Schr\"odinger operators $H$, the nonvanishing of $W(u_1, u_2)(x)$ for
Floquet solutions $u_1 = \psi_{\varepsilon_1}(E_1)$, 
$u_2 = \psi_{\varepsilon_2}(E_2)$, \;
$\varepsilon_{1,2} \in \{+,-\}$ of $H_p$, for $E_1$ and  $E_2$ in the
same spectral gap of $H$, is proved in \cite{FIT87}. The extension of these
ideas to general one-dimensional Schr\"odinger operators was done in
\cite{GST96} and is reviewed in Section \ref{s2} of this survey. So while 
\cite{GST96a} is not directly related to the overarching inverse spectral theory 
topic of this survey, we  decided to include it because of its relevance in
connection with Section \ref{s2}.

To set the stage, we consider Sturm--Liouville differential expressions
of the form
$$
(\tau u)(x)=r(x)^{-1}[-(p(x)u'(x))' + q(x)u(x)], \quad x\in(a,b), \quad
-\infty\leq a<b\leq
\infty
$$
where
\begin{equation*}
r, p^{-1}, q \in L^1_{\loc}((a,b)) \text{ are real-valued
and $r,p >0$ a.e.~on $(a,b)$}. 
\end{equation*}
We shall use $\tau$ to describe the formal differentiation expression and
$H$ for the operator in $L^{2}((a,b); r\,dx)$ given by $\tau$ with separated
boundary conditions at $a$ and/or $b$.

If $a$ (resp., $b$) is finite and $q,p^{-1},r$ are in addition integrable
near $a$ (resp., $b$), $a$ (resp., $b$) is called a {\it{regular}}
end point. $\tau$, respectively $H$, is called {\it{regular}} if
both $a$ and $b$ are regular. As is usual (\cite[Sect.\ XIII.2]{DS88}, 
\cite[Sect.\ 17]{Na68}, \cite[Ch.\ 3]{We87}), we consider the local
domain
\begin{equation*}
D_{\loc}=\{u \in AC_{\loc}((a,b)) \,|\, pu'\in
AC_{\loc}((a,b)),\ \tau u \in L^{2}_{\loc}((a,b); r\,dx)\}, 
\end{equation*}
where $AC_{\loc}((a,b))$ is the set of locally
absolutely continuous functions on $(a,b)$. General ODE
theory shows that for any $E\in\bbC$, $x_0\in (a,b)$, and $(\alpha,
\beta)\in\bbC^2$, there is a unique $u\in D_{\loc}$
such that $-(pu')'+qu - Eru=0$ for a.e.~$x \in (a,b)$ and $(u(x_0),
(pu')(x_0))=(\alpha,\beta)$.

The maximal and minimal operators are defined by taking
$$
\dom(T_{\rm max})=\{u\in L^{2}((a,b); r\,dx)\cap
D_{\loc}\,|\, \tau u\in L^{2}((a,b); r\,dx)\},
$$
with
\begin{equation*}
T_{\rm max}u=\tau u. 
\end{equation*}
$T_{\rm min}$ is the operator closure of
$T_{\rm max}\restriction D_{\loc}\cap\{u
\text{ has compact support in $(a,b)$}\}$. Then $T_{\rm min}$
is symmetric and $T^{*}_{\rm min}=T_{\rm max}$.

According to Weyl's theory of self-adjoint extensions (\cite[Sect.\ 
XIII.6]{DS88}, \cite[Sect.\ 18]{Na68}, \cite[App.\ to X.1]{RS75},  
\cite[Section 8.4]{We80}, \cite[Chs.\ 4, 5]{We87}), the deficiency indices
of $T_{\rm min}$ are $(0,0)$ or $(1,1)$ or $(2,2)$ depending
on whether it is limit point at both, one, or neither endpoint.
Moreover, the self-adjoint extensions can be described in terms of
Wronskians (\cite[Sect.\ XIII.2]{DS88}, \cite[Sects.\ 17, 18]{Na68},  
\cite[Sect.\ 8.4]{We80}, \cite[Ch.\ 3]{We87}). Define
\begin{equation*}
W(u_1,u_2)(x)=u_1(x)(pu_2')(x)-(pu_1')(x)u_2(x). 
\end{equation*}
Then if $T_{\rm min}$ is limit point at both ends,
$T_{\rm min}=T_{\rm max}=H$. If $T_{\rm min}$
is limit point at $b$ but not at $a$, for $H$ any self-adjoint extension
of $T_{\rm min}$, if $\varphi_-$ is any function in $\dom(H)
\backslash \dom(T_{\rm min})$, then
$$
\dom(H)=\{u\in \dom(T_{\rm max})\,|\, W(u,\varphi_-)(x)\to 0
\text{ as $x\downarrow a$}\}.
$$
Finally, if $u_1$ is limit circle at both ends, the operators $H$ with
separated boundary conditions are those for which we can find
$\varphi_\pm \in \dom(H)$, $\varphi_+= 0$ near $a$, $\varphi_-
= 0$ near $b$, and $\varphi_\pm \in \dom(H)\backslash
\dom(T_{\rm min})$. In that case,
$$
\dom(H)=\{u\in D(T_{\rm max}) \,|\,  W(u,\varphi_-)(x)\to 0
\text{ as $x\downarrow a$}, \,
W(u,\varphi_+)(x)\to 0 \text{ as $x\uparrow b$}\}.
$$
Of course, if $H$ is regular, we can just specify the boundary
conditions by taking values at $a,b$ since by regularity any $u\in
\dom(T_{\rm max})$ has $u,pu'$ continuous on $[a,b]$.
It follows from this analysis that
$$
\text{if $u_{1,2}\in \dom(H)$, then $W(u_1,u_2)(x)
\to 0$ as $x\to a$ or $b$.}
$$

Such operators will be called SL operators (for Sturm--Liouville, but SL
includes separated boundary conditions (if necessary)) and denoted by $H$.

It will be convenient to write $\ell_-=a$, $\ell_+=b$.

Throughout this section we will denote by $\psi_{\pm}(z,x) \in
D_{\loc}$ solutions of $\tau \psi = z \psi$ so that
$\psi_{\pm}(z,\,.\,)$ is $L^2$ at $\ell_{\pm}$ and $\psi_{\pm}
(z,\,.\,)$ satisfies the appropriate boundary condition at $\ell_\pm$
in the sense that for any $u \in \dom(H)$, $\lim_{x\to\ell_{\pm}}
W(\psi_{\pm}(z),u)(x)=0$. If $\psi_{\pm}(z,\,.\,)$ exist, they are
unique up to constant multiples. In particular, $\psi_{\pm}(z,\,.\,)$
exist for $z$ not in the essential spectrum of $H$ and we can assume
them to be holomorphic with respect to $z$ in ${\bbC} \backslash
\sigma (H)$ and real for $z\in\bbR$. One can choose
\begin{equation}
\psi_\pm(z,x) = ((H-z)^{-1} \chi_{(c,d)})(x) \,
\text{ for } \, x<c \, \text{ and } x> d, \quad a<c<d<b  \lb{3.1.1}
\end{equation}
and uniquely continue $\psi_\pm(z,x)$  for $x > c$ and $x<d$.
Here $(H-z)^{-1}$ denotes the resolvent of $H$ and $\chi_\Omega$ the
characteristic function of the set $\Omega  \subseteq \bbR$. Clearly
we can include a finite number of isolated eigenvalues in the domain
of holomorphy of $\psi_\pm$ by removing the corresponding poles. Moreover,
to simplify notations, all solutions $u$ of $\tau u = Eu$ are understood
to be not identically vanishing and solutions associated with real values
of the spectral parameter $E$ are assumed to be real-valued in this paper.
Thus if $E$ is real and in the resolvent set for $H$ or an isolated
eigenvalue, we are guaranteed there are solutions that satisfy the boundary
conditions at $a$ or $b$. If $E$ is in the essential
spectrum,  it can happen  that such solutions do not exist or it may happen that they do.
In Theorems \ref{t3.1.3}, \ref{t3.1.4} below, we shall explicitly assume
such solutions exist for the energies of interest. If these energies are
not in the essential spectrum, that is automatically fulfilled.

The key idea in \cite{GST96a} is to look at zeros of the Wronskian. That
zeros of the Wronskian are related to oscillation theory is indicated by
an old paper of Leighton \cite{Le52}, who noted that if $u_j,pu_j' \in
AC_{\loc}((a,b))$, $j=1,2$ and $u_1$ and $u_2$ have a
nonvanishing Wronskian $W(u_1,u_2)$ in $(a,b)$, then their zeros must
intertwine each other. (In fact, $pu_1'$ must have opposite signs at
consecutive zeros of $u_1$, so by nonvanishing of $W$, $u_2$ must have
opposite signs at consecutive zeros of $u_1$ as well. Interchanging the
role of $u_1$ and $u_2$ yields strict interlacing of their zeros.)
Moreover, let $E_1<E_2$ and $\tau u_j = E_j u_j$, $j=1,2$. If $x_0,x_1$
are two consecutive zeros of $u_1$, then the number of zeros of $u_2$
inside $(x_0,x_1)$ is equal to the number of zeros of the Wronskian
$W(u_1,u_2)$ plus one (cf.\ Theorem \ref{t3.7.4}). Hence the Wronskian
comes with a built-in renormalization counting the additional zeros of
$u_2$ in comparison to $u_1$.  

We let $W_0(u_1,u_2)$ be the number of zeros of the Wronskian in the
open interval $(a,b)$ not counting multiplicities of zeros. Given
$E_1<E_2$, we let $N_0(E_1, E_2)=\dim(
\ran( P_{(E_1, E_2)}(H)))$ be the dimension of the spectral
projection $P_{(E_1, E_2)}(H)$ of $H$.

We begin by presenting two aspects of zeros of the Wronskian which are
critical for the two halves of our proofs (i.e., for showing $N_0\geq
W_0$ and that $N_0\leq W_0$). First, the vanishing of the Wronskian lets
us patch solutions together:

\begin{lemma} \lb{l3.3.1}  Suppose that $\psi_{+,j}, \psi_-\in
D_{\loc}$ and that $\psi_{+,j}$ and $\tau\psi_{+,j}$,
$j=1,2$ are in $L^{2}((c,b))$ and that $\psi_-$ and $\tau\psi_-$ are
in $L^{2}((a,c))$ for all $c\in (a,b)$.  Suppose, in addition, that
$\psi_{+,j}$, $j=1,2$ satisfy the boundary condition defining $H$ at
$b$ $($i.e., $W(u,\psi_{+,j})(c)\to 0$ as $c\uparrow b$ for all $u\in
\dom(H)$$)$ and similarly, that $\psi_-$ satisfies the boundary condition
at $a$. Then, \\
$(i)$ If $W(\psi_{+,1}, \psi_{+,2})(c)=0$ and $(\psi_{+,2}(c),
(p\psi'_{+,2})(c)) \neq (0,0)$, then there exists a $\gamma$ such that
$$
\eta=\chi_{[c,b)}(\psi_{+,1}-\gamma\psi_{+,2})\in \dom(H)
$$
and
\begin{equation*}
H\eta=\chi_{[c,b)}(\tau\psi_{+,1}-\gamma\tau\psi_{+,2}). 
\end{equation*}
$(ii)$ If $W(\psi_{+,1}, \psi_-)(c)=0$ and $(\psi_-(c),
(p\psi_-')(c))\neq (0,0)$, then there is a $\gamma$ such that
$$
\eta=\gamma\chi_{(a,c]}\psi_- + \chi_{(c,b)}\psi_{+,1}\in \dom(H)
$$
and
\begin{equation*}
H\eta =\gamma\chi_{(a,c]}\tau\psi_- + \chi_{(c,b)}\tau\psi_{+,1}.
\end{equation*}
\end{lemma}

The second aspect connects zeros of the Wronskian to Pr\"ufer
variables $\rho_u, \theta_u$ (for $u,pu'$ continuous) defined by
$$
u(x)=\rho_u(x)\sin(\theta_u(x)), \quad (pu')(x)=\rho_u(x)
\cos(\theta_u(x)).
$$
If $(u(x), (pu')(x))$ is never $(0,0)$, then $\rho_u$ can be chosen
positive and $\theta_u$ is uniquely determined once a value of $\theta_u
(x_0)$ is chosen subject to the requirement that $\theta_u$ be continuous in $x$.

Notice that
$$
W(u_1,u_2)(x)= \rho_{u_1}(x)\rho_{u_2}(x)\sin(\theta_{u_1}(x)
- \theta_{u_2}(x)).
$$
Thus, one obtains the following results.

\begin{lemma} \lb{l3.3.2} Suppose $(u_j,pu_j')$, $j=1,2$ are never
$(0,0)$. Then $W(u_1,u_2)(x_0)$ is zero if and only if $\theta_{u_1}
(x_0)= \theta_{u_2}(x_0) \, (\text{mod} \; \pi)$.
\end{lemma}
In linking Pr\"ufer variables to rotation numbers, an important role
is played by the observation that because of
$$
u(x) = \int_{x_0}^x dt \, \frac{\rho_u(t) \cos(\theta_u(t))}{p(t)},
$$
$\theta_u(x_0)= 0 \, (\text{mod} \; \pi)$ implies $[\theta_u(x)-
\theta_u(x_0)]\big/ (x-x_0) >0$ for $0 <|x - x_0|$ sufficiently small
and hence for all $0 <|x - x_0|$ if $(u,pu') \neq (0,0)$. (In fact,
suppose
$x_1 \ne x_0$ is the closest $x$ such that $\theta_u(x_1)=\theta_u
(x_0)$ then apply the local result at $x_1$ to obtain a contradiction.)
We summarize:

\begin{lemma} \lb{l3.3.3} $(i)$ If $(u,pu') \neq (0,0)$ then
$\theta_u(x_0)= 0 \, (\text{mod} \; \pi)$ implies
$$
[\theta_u(x)-\theta_u(x_0)]\big/ (x-x_0) >0
$$
for $x \ne x_0$. In particular, if $\theta_u(c)\in [0,\pi)$ and $u$ has
$n$ zeros in $(c,d)$, then $\theta_u(d-\epsilon)\in (n\pi, (n+1)\pi)$
for sufficiently small $\epsilon > 0$. \\
$(ii)$ Let $E_1<E_2$ and assume that $u_{1,2}$
solve $\tau u_j = E_j u_j$, $j=1,2$. Let $\Delta(x)=\theta_{u_2}(x)-
\theta_{u_1}(x)$. Then $\Delta(x_0)= 0 \, (\text{mod} \; \pi)$ implies
$(\Delta(x)-\Delta(x_0))/ (x-x_0)>0$ for $0<|x -x_0|$.
\end{lemma}

\begin{remark} \lb{r3.3.4} $(i)$ Suppose $r,p$ are continuous on $(a,b)$.
If $\theta_{u_1}(x_0)= 0 \, (\text{mod} \; \pi)$ then $\theta_{u_1}(x) -
\theta_{u_1}(x_0) = c_0(x-x_0) + o(x-x_0)$ with $c_0>0$. If $\Delta
(x_0)= 0 \, (\text{mod} \; \pi)$ and $\theta_{u_1}(x_0)\neq 0 \, (\text{mod} \; \pi)$, then $\Delta(x)-\Delta (x_0)=c_1(x-x_0)+o(x-x_0)$ with
$c_{1}>0$. If $\theta_{u_1} (x_0)= 0=\Delta(x_0) \, (\text{mod} \; \pi)$,
then $\Delta(x) - \Delta(x_0)=c_2(x-x_0)^{3}+o(x-x_0)^{3})$ with $c_2>0$.
Either way, $\Delta$ increases through $x_0$. (In fact, $c_0=p(x_0)^{-1}$,
$c_{1}= (E_2-E_1)r(x_0)\sin^2 (\theta_{u_1}(x_0))$ and $c_{2}=\frac{1}{3}
r(x_0)p(x_0)^{-2}(E_2-E_1))$. \\
$(ii)$ In other words, Lemma \ref{l3.3.3} implies that the integer parts
of $\theta_u/ \pi$ and $\Delta_{u,v}/ \pi$ are increasing with respect
to $x\in(a,b)$ (even though $\theta_u$ and $\Delta_{u,v}$ themselves
are not necessarily monotone in $x$). \\
$(iii)$ Let $E \in [E_1,E_2]$ and assume $[E_1,E_2]$ to be outside the
essential spectrum of $H$. Then, for $x \in (a,b)$ fixed,
\begin{equation*}
\frac{d\theta_{\psi_\pm}}{dE}\,(E,x) =
-\frac{\int^{\ell_\pm}_x dt \, \psi_\pm(E,t)^2}
{\rho_{\psi_\pm}(E,x)} 
\end{equation*}
proves that $\mp\theta_{\psi_\pm}(E,x)$ is strictly increasing with
respect to $E$.
\end{remark}

We continue with some preparatory results in the regular case.

\begin{lemma} \lb{l3.4.2} Assume $H$ to be a regular SL operator. \\
$(i)$ Let $u_{1,2}$ be eigenfunctions of $H$ with eigenvalues $E_1<E_2$
and let $\ell$ be the number of eigenvalues of $H$ in $(E_1, E_2)$. Then
$W(u_1,u_2)(x)$ has exactly $\ell$ zeros in $(a,b)$. \\
$(ii)$ Let $E_1\leq E_2$ be eigenvalues of $H$ and suppose $[E_1, E_2]$ has
$\ell$ eigenvalues. Then for $\epsilon \geq 0$ sufficiently small,
$W_0(\psi_-(E_1 - \epsilon), \psi_+(E_2+\epsilon)) =\ell$. \\
$(iii)$ Let $E_3<E_4<E$ and $u$ be any solution of $\tau u=E u$. Then,
\begin{equation}
W_0(\psi_-(E_3), u)\geq W_0(\psi_-(E_4), u). \lb{3.4.3}
\end{equation}
Similarly, if $E_3>E_4>E$ and $u$ is any solution of $\tau u=E u$,
then \eqref{3.4.3} holds. \\
$(iv)$ Item $(iii)$ remains true if every $\psi_-$ is
replaced by a $\psi_+$.
\end{lemma}

\begin{remark} \lb{r3.4.3a} $(i)$ Since $(E_1, E_2)$ has $\ell-2$
eigenvalues, Lemma \ref{l3.4.2}\,$(i)$ implies that the Wronskian
$W(\psi_-(E_1), \psi_+(E_2))(x)$ has $\ell-2$ zeros in
$(a,b)$ and clearly it has zeros at $a$ and $b$. Essentially, Lemma
\ref{l3.4.2}\,$(ii)$ implies that replacing $E_1$ by $E_1-\epsilon$ and
$E_2$ by $E_2+\epsilon$ moves the zeros at $a,b$ inside $(a,b)$ to give
$\ell-2+2=\ell$ zeros. \\
$(ii)$ Lemma \ref{l3.4.2}\,$(iv)$ follows from Lemma \ref{l3.4.2}\,$(iii)$
upon reflecting at some point $c \in (a,b)$, implying an interchange of
$\psi_+$ and $\psi_-$.
\end{remark}

Lemma \ref{l3.4.2} then implies the following  result.

\begin{lemma} \lb{l3.4.1} Let $H$ be a regular SL operator and suppose
$E_1<E_2$. Then,
\begin{equation*}
W_0(\psi_-(E_1), \psi_+(E_2))\geq N_0(E_1, E_2). 
\end{equation*}
\end{lemma}

Using the approach of Weidmann (\cite[Ch.\ 14]{We87}) to control some
limits, one can remove the assumption that $H$ is regular in Lemma
\ref{l3.4.1} as follows.

Fix functions $u_1,u_2 \in D_{\loc}$. Pick $c_{n}\downarrow a$,
$d_{n}\uparrow b$. Define $\tilde{H}_n$ on $L^{2}((c_{n}, d_{n}); r\,dx)$
by imposing the following boundary conditions on
$\eta\in \dom(\tilde{H}_{n})$
\begin{equation*}
W(u_1,\eta)(c_{n})=0=W(u_2,\eta)(d_{n}).  
\end{equation*}
On $L^{2}((a,b);r\,dx)= L^{2}((a,c_{n});r\,dx)\oplus L^{2}
((c_{n}, d_{n});r\,dx)\oplus L^{2}((d_{n}, b);r\,dx)$ take $H_{n}
=\alpha I\oplus\tilde{H}_{n}\oplus\alpha I$ with $\alpha$ a
fixed real constant. Then Weidmann proves:

\begin{lemma} \lb{l3.5.1} Suppose that either $H$ is limit point at
$a$ or that $u_1$ is a $\psi_-(E,x)$ for some $E$ and similarly, that
either $H$ is limit point at $b$ or $u_2$ is a $\psi_+(E',x)$ for
some $E'$. Then $H_n$ converges to $H$ in strong resolvent sense as
$n\to\infty$.
\end{lemma}

The idea of Weidmann's proof is that it suffices to find a core $D_0$ of
$H$ such that for every $\eta \in D_0$ there exists an $n_0 \in \bbN$
with $\eta\in D_0$ for $n \geq n_0$ and $H_{n}\eta\to H\eta$ as $n$
tends to infinity (see \cite[Theorem 9.16\,(i)]{We80}). If $H$ is limit
point at both ends, take $\eta\in D_0=\{u\in 
D_{\loc}\,|\,
\text{supp}(u)\text{ compact in }(a,b)\}$. Otherwise, 
pick
$\tilde{u}_1,\tilde{u}_2\in \dom(H)$ with $\tilde{u}_2 =u_2$ 
near $b$ and
$\tilde{u}_2=0$ near $a$ and with $\tilde{u}_1=u_1$ near 
$a$ and
$\tilde{u}_1=0$ near $b$. Then pick $\eta\in D_0+\text{span}[\tilde{u}_1,
\tilde{u}_2]$ which one can show is a core for $H$ (\cite[Ch.\ 14]{We87}).

Secondly one uses the following fact:

\begin{lemma} \lb{l3.5.2} Let $A_{n}\to A$ in strong resolvent sense as $n
\to \infty$. Then,
\begin{equation*}
\dim(\ran(P_{(E_1,E_2)}(A)))\leq
\varliminf_{n \to \infty} \dim(\ran(
P_{(E_1,E_2)}(A_{n}))).  
\end{equation*}
\end{lemma}

Combining Lemmas \ref{l3.4.1}--\ref{l3.5.2} then yields the following
result.

\begin{theorem} \lb{t3.1.5} Let $E_1 < E_2$. If $u_1=\psi_-(E_1)$ and
either $u_2=\psi_+(E_2)$ or $\tau u_2 = E_2 u_2$ and $H$ is limit point
at $b$. Then,
$$
W_0(u_1,u_2)\geq\dim(\ran((P_{(E_1, E_2)}(H))).
$$
\end{theorem}

Next, we indicate how the following result can be proved:

\begin{theorem} \lb{t3.1.6} Let $E_1 < E_2$. Let either $u_1=\psi_+(E_1)$
or
$u_1=\psi_-(E_1)$ and either $u_2=\psi_+(E_2)$ or $u_2=\psi_-(E_2)$.
Then,
\begin{equation}
W_0(u_1,u_2)\leq\dim(\ran(P_{(E_1, E_2)}(H))). \lb{3.1.8}
\end{equation}
\end{theorem}

Let $E_1<E_2$. Suppose
first that $u_1=\psi_-(E_1)$ and $u_2=\psi_+(E_2)$. Let $x_1,\dots,
x_{m}$ be zeros of $W(u_1,u_2)(x)$. Suppose we can prove that $\dim
P_{(E_1, E_2)}(H) \geq m$. If $W_0(u_1,u_2)=m$, this proves \eqref{3.1.8}.
If
$W_0=\infty$, we can take $m$ arbitrarily large, and again \eqref{3.1.8}
holds. Define
$$
\eta_j(x) = \begin{cases} u_1(x), & x\leq x_j, \\
\gamma_j u_2(x), & x\geq x_j, \end{cases} \quad 1 \le j \le m,
$$
where $\gamma_j$ is defined such that $\eta_j\in \dom(H)$ by Lemma
\ref{l3.3.1}. Let
$$
\tilde{\eta}_j(x) = \begin{cases} u_1(x), & x\leq x_j, \\
-\gamma_j u_2(x), & x>x_j, \end{cases} \quad 1 \le j \le m.
$$
If $E_2$ is an eigenvalue of $H$, we define in addition $\eta_0 = u_2
= -\tilde{\eta}_0$, $x_0=a$ and if $E_1$ is an eigenvalue of $H$,
$\eta_{m+1} = u_1 = \tilde{\eta}_{m+1}$, $x_{m+1}=b$.

\begin{lemma} \lb{l3.6.1} $\langle\eta_j, \eta_{k}\rangle =\langle
\tilde{\eta}_j, \tilde{\eta}_{k}\rangle$ for all $j,k$, where
$\langle \,\cdot\, ,\,\cdot\,\rangle$ is the $L^{2}((a,b);r\,dx)$
inner product.
\end{lemma}

Notice that by (3.2),
\begin{equation*}
\biggl(H-\frac{E_2+E_1}{2}\biggr)\, \eta_j=
\biggl(\frac{E_1-E_2}{2}\biggr)\, \tilde{\eta}_j.  
\end{equation*}
This result and Lemma \ref{l3.6.1} imply the following lemma.

\begin{lemma} \lb{l3.6.2} If $\eta$ is in the span of the $\eta_j$, then
$$
\left\|\biggl(H-\frac{E_2+E_1}{2}\biggr)\, \eta\right\| =
\frac{|E_2-E_1|}{2}\, \|\eta\|.
$$
\end{lemma}

Thus, $\dim(\ran( P_{[E_1, E_2]}(H)))\geq
\dim(\text{span}(\{\eta_j\}))$. But $u_1$ and $u_2$ are independent
on each interval (since their Wronskian is nonconstant) and so the
$\eta_j$ are linearly independent. This proves Theorem \ref{t3.1.6} in the
$\psi_-(E_1), \psi_+(E_2)$ case. The case $u_1=\psi_-(E_1)$,
$u_2=\psi_-(E_2)$ is proved similarly. The cases $u_1=\psi_+(E_1)$,
$u_2=\psi_\pm(E_2)$ can be obtained by reflection.

Next one proves the following result.

\begin{theorem} \lb{t3.7.1} Let $E_1\neq E_2$. Let $\tau u_j=E_j u_j$,
$j=1,2$, $\tau v_2=E_2 v_2$ with $u_2$ linearly independent of $v_2$.
Then the zeros of $W(u_1,u_2)$ interlace the zeros of $W(u_1,v_2)$ and
vice versa $($in the sense that there is exactly one zero of one
function in between two zeros of the other$)$. In particular,
$|W_0(u_1,u_2)-W_0(u_1,v_2)| \leq 1$.
\end{theorem}

Theorems \ref{t3.1.5}, \ref{t3.1.6}, and \ref{t3.7.1} then yield
the following two theorems, the principal results of \cite{GST96a}:

\begin{theorem} \lb{t3.1.3}  Suppose $E_1 < E_2$. Let $u_1=\psi_-(E_1)$ and
$u_2=\psi_+(E_2)$. Then,
$$
W_0(u_1,u_2)=N_0(E_1, E_2).
$$
\end{theorem}

\begin{theorem} \lb{t3.1.4} Suppose $E_1 < E_2$. Let $u_1=\psi_-(E_1)$ and
$u_2=\psi_-(E_2)$. Then either,
\begin{equation}
W_0(u_1,u_2)=N_0(E_1, E_2),  \lb{3.1.6}
\end{equation}
or,
\begin{equation}
W_0(u_1,u_2)=N_0(E_1, E_2) -1. \lb{3.1.7}
\end{equation}
If either $N_0=0$ or $H$ is limit point at $b$, then \eqref{3.1.6} holds.
\end{theorem}

One infers that if $b$ is a regular point and $E_2 > e > E_1$ with $e$
an eigenvalue and $|E_2-E_1|$ is small, then \eqref{3.1.7} holds rather
than \eqref{3.1.6}. One also sees that if $u_{1,2}$ are arbitrary
solutions of $\tau u_j=E_j u_j$, $j=1,2$, then, in general, $|W_0-N_0|\leq
2$ (this means that if one of the quantities is infinite, the other is as
well) and we note that any of $0, \pm 1,\pm 2$ can occur for $W_0-N_0$.
Especially, if either $E_1$ or $E_2$ is in the interior of the essential
spectrum of $H$ (or $\dim(\ran(P_{(E_1, E_2)}(H)))=\infty$), then $W_0
(u_1,u_2)=\infty$ for any $u_1$ and $u_2$ satisfying $\tau u_j=E_j
u_j$, $j=1,2$ (cf.\ Theorem \ref{t3.7.3}).

\begin{remark} \lb{r3.1.4a} Of course, by reflecting about a point
$c \in (a,b)$, Theorems \ref{t3.1.5}, \ref{t3.1.3}, and \ref{t3.1.4} hold
for $u_1= \psi_+(E_1)$ and $u_2 =
\psi_-(E_2)$ (and either $N_0 = 0$ or $H$ is limit point at $a$ in the
corresponding analog of Theorem \ref{t3.1.4} yields \eqref{3.1.6}) and
similarly,
$\tau u_2 = E_2 u_2$ and $H$ is limit point at $a$ yields the conclusion
in the corresponding analog of Theorem \ref{t3.1.5}.
\end{remark}

We add a few more results proved in \cite{GST96a}.

By applying Theorem \ref{t3.7.1} twice, one concludes

\begin{theorem} \lb{t3.7.2} Let $E_1\neq E_2$. Let $u_1, u_2, v_1, v_2$ be
the linearly independent functions with $\tau u_j=E_j u_j$ and $\tau v_j
=E_j v_j$. Then,
$$
|W_0(u_1, u_2)-W_0(v_1, v_2)|\leq 2.
$$
\end{theorem}

\begin{theorem} \lb{t3.7.3} If $\dim(\ran(P_{(E_1,E_2)}(H)))
=\infty$, then $W_0(u_1,u_2)=\infty$ for any $u_1$ and $u_2$ satisfying
$\tau u_j = E_j u_j$, $j=1,2$.
\end{theorem}

\begin{theorem} \lb{t3.7.4} Let $E_1<E_2$. Let $\tau u_j = E_j u_j$,
$j=1,2$. If $a < x_0 < x_1 < b$ are zeros of $u_1$ or of $W(u_1, u_2)
(\,.\,)$, then the number of zeros of $u_2$ inside $(x_0,x_1)$ equals
the number of zeros of $W(u_1,u_2)(\,.\,)$ inside $(x_0,x_1)$ plus the
number of zeros of $u_1$ inside $(x_0,x_1)$ plus one.
\end{theorem}

The following result is of special interest in connection with
the problem of whether the total number of eigenvalues of $H$ in
one of its essential spectral gaps is finite or infinite. In
particular, the energies $E_1, E_2$ in Theorem 7.5 below may lie
in the essential spectrum of $H$. For this purpose we consider an
auxiliary Dirichlet operator $H_{x_0}^D$, $x_0\in(a,b)$ associated
with $H$. $H_{x_0}^D$ is obtained by taking the direct sum of the
restrictions $H_{x_0,\pm}^D$ of $H$ to $(a,x_0)$, respectively
$(x_0,b)$, with a Dirichlet boundary condition at $x_0$. We
emphasize that the Dirichlet boundary conditions can be replaced
by boundary conditions of the type $\lim_{\epsilon\downarrow 0}
[u'(x_0\pm\epsilon)+ \beta u(x_0\pm\epsilon)] =0$, $\beta\in\bbR$.

\begin{theorem} \lb{t3.7.5}
Let $E_1<E_2$. Let $\tau u_j = E_j u_j$, $\tau s_j = E_j s_j$,
and $s_j(E_j,x_0)=0$, $j=1,2$. Then,
\begin{align*}
& \dim(\ran(P_{(E_1,E_2)}(H)))<\infty
\, \text { if and only if } \, W_0(u_1,u_2)<\infty, \\
& \dim(\ran(P_{(E_1,E_2)}(H)))-1
\le \dim(\ran(P_{(E_1,E_2)}(H^D_{x_0}))) \\
& \hspace*{4.13cm} \le
\dim(\ran(P_{(E_1,E_2)}(H)))+2,  \\
& W_0(s_1,s_2)-1 \le \dim(\ran(
   P_{(E_1,E_2)}(H^D_{x_0}))) \le W_0(s_1,s_2)+1.
\end{align*}
\end{theorem}

For an application of this circle of ideas to the notion of the density of
states we refer to \cite{GST96a}.

\medskip

{\it More recent references:} Oscillation and renormalized oscillation theory was also put in perspective by Simon's contribution \cite{Si05} to the the Festschrift \cite{AHP05} in honor of Sturm and Liouville.  
Renormalized oscillation theory for one-dimensional Jacobi and Dirac-type operators was
developed by Teschl \cite{Te96} (see also \cite[Sect.\ 4.3]{Te00}) and \cite{Te98b}. For an
interesting application of some of the results in \cite{GST96a} to the stability theory of complete minimal surfaces we refer to a paper by Schmidt \cite{Sc01}. For additional results on oscillation theory, critical coupling constants, and eigenvalue asymptotics we refer to Schmidt \cite{Sc00}.

\section{Trace Formulas for Schr\"odinger and Jacobi Operators: The xi
 Function} \label{s4}

In this section we summarize some of the principal results of
the following papers:  \\
\cite{GHS95} F.~Gesztesy, H.~Holden, and B.~Simon, {\it Absolute
summability of the trace relation for certain Schr\"odinger operators},
Commun.  Math. Phys. {\bf 168}, 137--161 (1995). \\
\cite{GHSZ93} F.~Gesztesy, H.~Holden, B.~Simon, and Z.~Zhao, {\it Trace
formulae and inverse spectral theory for Schr\"odinger operators},
Bull. Amer. Math. Soc. {\bf 29}, 250--255 (1993). \\
\cite{GHSZ95} F.~Gesztesy, H.~Holden, B.~Simon, and Z.~Zhao, {\it Higher
order trace relations for Schr\"odinger operators}, Rev. Math.
Phys. {\bf 7}, 893--922 (1995).\\
\cite{GHSZ96} F.~Gesztesy, H.~Holden, B.~Simon, and Z.~Zhao, {\it A trace
formula for multidimensional Schr\"odinger operators}, J. Funct. Anal.
{\bf 141}, 449--465 (1996). \\
\cite{GS96} F.~Gesztesy and B.~Simon, {\it The xi function}, Acta Math.
{\bf 176}, 49--71 (1996).

\medskip\smallskip
 
We start with \cite{GS96}. One of the principal goals in \cite{GS96} was
to introduce a special function
$\xi(\cdot,\cdot)$ on $\bbR\times\bbR$ associated with one-dimensional
Schr\"odinger operators $H$ (and Jacobi operators $h$) which led to a 
generalization of the known trace formula for periodic Schr\"odinger
operators for general potentials $V$ and established
$\xi$ as a new tool in the spectral theory of one-dimensional
Schr\"odinger operators and (multi-dimensional) Jacobi operators.
 
To illustrate this point we recall the well-known trace formula for
periodic potentials $V$ of period $a>0$: Then, by Floquet theory (see,
e.g., \cite{Ea73}, \cite{MW79}, \cite{RS78})
$$
\sigma (H)=[E_{0}, E_{1}]\cup [E_{2}, E_{3}]\cup\dots
$$
a set of bands.  If $V$ is $C^1(\bbR)$, one can show that the sum of the
gap sizes is finite, that is,
\begin{equation}
\sum^{\infty}_{n=1}|E_{2n}-E_{2n-1}|<\infty.    \lb{4.1.1}
\end{equation}
 
For fixed $y$, let $H_y$ be the operator $-\frac{d^2}{dx^2}+V$
in $L^{2}([y, y+a])$ with Dirichlet boundary conditions $u(y)=u(y+a)=0$ . 
Its spectrum is discrete, that is, there are eigenvalues
$\{\mu_{n}(y)\}^{\infty}_{n=1}$ with
\begin{equation}
E_{2n-1}\leq\mu_{n}(y)\leq E_{2n}.   \lb{4.1.2}
\end{equation}
The trace formula for $V$ then reads
\begin{equation}
V(y)=E_{0}+\sum^{\infty}_{n=1} [E_{2n}+E_{2n-1}-2\mu_{n}(y)].   \lb{4.1.3}
\end{equation}
By \eqref{4.1.2},
$
|E_{2n}+E_{2n-1}-2\mu_{n}(y)|\leq |E_{2n}-E_{2n-1}|
$
so \eqref{4.1.1} implies the convergence of the sum in \eqref{4.1.3}. An 
elegant direct proof of \eqref{4.1.3} can be found, for instance, in 
\cite[Sect.\ 26]{Si79}.
 
The earliest trace formula for Schr\"odinger operators was found on a
finite interval in 1953 by Gel'fand and Levitan \cite{GL87} with later
contributions by Dikii \cite{Di61}, Gel'fand \cite{Ge87},
Halberg and Kramer \cite{HK60}, and Gilbert and Kramer
\cite{GK64}.  The first trace formula for periodic $V$ was obtained in
1965 by Hochstadt \cite{Ho65}, who showed that for finite-band 
potentials
$$
V(x)-V(0)=2\sum^{g}_{n=1}[\mu_{n}(0)-\mu_{n}(x)].
$$
Dubrovin [9] then proved \eqref{4.1.3} for finite-band potentials.  The 
general formula \eqref{4.1.3} under the hypothesis that $V$ is periodic
and in $C^{\infty}(\bbR)$ was proved in 1975 by Flaschka \cite{Fl75}
and McKean and van Moerbeke \cite{MM75}, and later for general $C^3(\bbR)$
potentials by Trubowitz \cite{Tr77}.  Formula \eqref{4.1.3} is a key
element of the solution of inverse spectral problems for periodic
potentials \cite{Du75}, \cite{Fl75}, \cite{Ho65}, \cite[Ch.\ 11]{Le87},
\cite[Sect.\ 4.3]{Ma86}, \cite{MM75}, \cite{MT76},, \cite{Tr77}.  

There have been two classes of potentials for which \eqref{4.1.3} has been
extended.  Certain almost periodic potentials are studied in Craig
\cite{Cr89}, Levitan \cite{Le85}, \cite[Ch.\ 11]{Le87}, and Kotani-Krishna
\cite{KK88}.
 
In 1979, Deift and Trubowitz \cite{DT79} proved that if $V(x)$ decays
sufficiently rapidly at infinity and $-\frac{d^2}{dx^2}+V$ has no negative
eigenvalues, then
\begin{equation}
V(x)=\frac{2i}{\pi}\int^{\infty}_{-\infty}\,dk\,k\,
\ln\biggl[1+R(k)\frac{f_{+}(x, k)}{f_{-}(x, k)}\biggr]     \lb{4.1.4}
\end{equation}
(where $f_{\pm}(x, k)$ are the Jost functions at energy $E=k^{2}$ and
$R(k)$ is a reflection coefficient) which can be shown to be an analog
of \eqref{4.1.3}.  Previously, Venakides \cite{Ve88} studied a trace
formula for $V$, a positive smooth potential of compact support, by
writing \eqref{4.1.3} for the periodic potential
$
V_{L}(x)=\sum^{\infty}_{n=-\infty}V(x+nL)
$
and then taking $L$ to $\infty$.  He found an integral formula which 
is precisely \eqref{4.1.4} (although, this was not identified as such in
\cite{Ve88}). 
 
The basic definition of $\xi$ depends on the theory of the Lifshits--Krein
spectral shift function \cite{Kr62}.  If $A$ and $B$ are self-adjoint
operators bounded from below, that is, $A\geq\eta$, $B\geq\eta$ for some
real
$\eta$, and so that
$[(A+i)^{-1}- (B+i)^{-1}]$ is trace class, then there exists a measurable
function
$\xi(\lambda)$ associated with the pair $(B, A)$ so that
\begin{equation}
\Tr[f(A)-f(B)]=-\int_{\bbR} d\lambda \, f'(\lambda)\xi(\lambda) 
\lb{4.1.5}
\end{equation}
for a class of functions $f$ which are sufficiently smooth and which
decay sufficiently rapidly at infinity, and, in particular, for
$f(\lambda)=e^{-t\lambda}$ for any $t>0$; and so that
\begin{equation}
\xi(\lambda)=0 \, \text{ for } \, \lambda<\eta.   \lb{4.1.6}
\end{equation}
 
Moreover, \eqref{4.1.5}, \eqref{4.1.6} uniquely determine $\xi(\lambda)$ 
for a.e.~$\lambda$. Moreover, if $[(A+i)^{-1}-(B+i)^{-1}]$ is rank $n$,
then
$
|\xi(\lambda)|\leq n
$
and if $B\geq A$, then $\xi(\lambda)\geq 0$.
 
For the rank-one case of importance in this paper, an extensive study
of $\xi$ can be found in \cite{Si95}.
 
Let $V$ be a continuous function on $\bbR$ which is bounded from below.
Let $H=-\frac{d^2}{dx^2}+V$ in $L^2(\bbR)$ which is essentially self-adjoint 
on $C^{\infty}_{0} (\bbR)$ and let $H_{x}^D$ be the operator on
$L^{2}((-\infty, x))\oplus L^{2}((x, \infty))$ with $u(x)=0$ Dirichlet
boundary conditions.  Then $[(H_{x}^D+i)^{-1}-(H+i)^{-1}]$ is rank one,
so there results a spectral shift function $\xi(\lambda,x)$ for the pair
$(H_{x}^D, H)$ which, in particular, satisfies,
\begin{equation}
\Tr[e^{-tH}-e^{-tH_{x}^D}]=t\int^{\infty}_{0} d\lambda \, 
e^{-t\lambda}\xi(\lambda,x).     \lb{4.1.7}
\end{equation}
 
While $\xi$ is defined in terms of $H$ and $H_{x}^D$, there is a
formula that only involves $H$, or more precisely, the Green's
function $G(z,x,y)$ of $H$ defined by
$$
((H-z)^{-1}f)(x)=\int_{\bbR} dy\,G(z,x, y)f(y), \quad \Im(z)\neq 0. 
$$
Then by general principles, $\lim_
{\epsilon\downarrow 0}\, G(\lambda+i\epsilon,x, y)$ exists for
a.e.~$\lambda\in\bbR$, and
\begin{equation}
\xi(\lambda,x)=\frac{1}{\pi}\,\Arg \bigl(\lim_
{\epsilon\downarrow 0}\, G(\lambda+i\epsilon, x, x)\bigr). \lb{4.1.8a}
\end{equation}
 
This is {\it formally} equivalent to formulas that Krein \cite{Kr62} 
has for
$\xi$ but in a singular setting (i.e., corresponding to an infinite 
coupling constant).  With
this definition out of the way, we can state the general trace formula
derived in \cite{GS96}:
\begin{theorem} \lb{t4.3.1}  Suppose $V$ is a continuous function bounded
from below on $\bbR$.  Let $\xi(\lambda,x)$ be the spectral
shift function for the pair $(H_{x}^D, H)$ with $H_{x}^D$ the operator on
$L^{2}((-\infty, x))\oplus L^{2}((x, \infty))$ obtained from
$H=-\frac{d^2} {dx^2}+V$ in $L^2(\bbR)$ with a Dirichlet boundary 
condition at $x$.  Let $E_{0}\leq \inf\,\sigma(H)$. Then
\begin{equation}
V(x)=\lim_{\alpha\downarrow 0} \biggl[E_{0}+\int^{\infty}
_{E_0} d\lambda \, e^{-\alpha\lambda} [1-2\xi(\lambda,x)]\,
\biggr].
\lb{4.3.1}
\end{equation}
\end{theorem} 

In particular, if
$\int^{\infty}_{E_0} d\lambda \, |1-2\xi(\lambda,x)| <\infty$,
then
$$
V(x)=E_{0}+\int^{\infty}_{E_0} d\lambda \, [1-2\xi(\lambda,x)].
$$

We note that the trace formula extends to real-valued potentials $V\in
L^1_{\loc}(\bbR)$ as long as $H$ stays bounded from below (it then is in
the limit point case at $\pm\infty$). Equation \eqref{4.3.1} then holds at
all Lebesgue points of $V$ and hence for a.e.\ $x\in\bbR$.
 
For certain almost periodic potentials, Craig \cite{Cr89} used a
regularization similar to the $\alpha$-regularization in \eqref{4.3.1}.
 
Basically, \eqref{4.3.1} follows from \eqref{4.1.7} and an asymptotic
formula,
\begin{equation}
\Tr[e^{-tH}-e^{-tH_{x}^D}]=\frac{1}{2}[1-tV(x)+o(t)].  \lb{4.3.2}
\end{equation}
We present a few examples next:
\begin{example} \lb{e4.3.2} Pick a constant $C\in\bbR$ such that $V(x)=C$ 
for all $x\in\bbR$. Then 
$G(\lambda,x,x)=(C-\lambda)^{-1/2}/2$ and hence one infers that 
$\Arg(G(\lambda,x,x))=0$ (resp., $\pi/2$) if
$\lambda<C$ (resp., $\lambda>C$).  Thus, by \eqref{4.1.8a},
$\xi(\lambda,x)=1/2$ on $(C, \infty)$ and $\xi(\lambda,x)=0$ on $(-\infty,C)$.  
When $\xi=1/2$ on a subset 
of $\sigma (H)$, that set does not contribute to the integral in 
\eqref{4.3.1} and one verifies for $E_0 \leq C$, 
$$
V(x)=E_0+\int_{E_0}^C d\lambda = E_0 + (C-E_0)=C, \quad x\in\bbR.
$$
\end{example}

\begin{example} \lb{e4.3.3} Suppose that $V(x)\to\infty$ as $|x|\to\infty$.
Then $H$ has eigenvalues $E_{0}<E_{1}<E_{2}<\cdots$ and $H^{D}_{x}$ has
eigenvalues $\{\mu_{j}(x)\}^{\infty}_{j=1}$ with 
$E_{j-1}\leq\mu_{j}(x)\leq E_{j}$.  We have
$$
\xi(\lambda,x) = \begin{cases} 1, & E_{j-1}<\lambda <\mu_{j}(x), \\
0, & \lambda<E_{0} \text{ or } \mu_{j}(x)<\lambda<E_{j}.
\end{cases}
$$
Thus \eqref{4.3.1} becomes:
\begin{equation}
V(x)=E_{0}+\lim_{\alpha\downarrow 0}\, \biggl[\sum^{\infty}_{j=1}
\bigl(\left.2e^{-\alpha\mu_{j}(x)}-e^{-\alpha E_{j}}-e^{-\alpha E_{j-
1}}\bigr)\right/\alpha\biggr].    \lb{4.3.3}
\end{equation}
If we could take $\alpha$ to zero inside the sum, we would get
\begin{equation}
V(x)=E_{0}+\sum^{\infty}_{j=1}\,[E_{j}+E_{j-1}-2\mu_{j}(x)]
\;\; \text{ (formal!)}    \lb{4.3.4}
\end{equation}
which is just a limit of the periodic formula \eqref{4.1.3} in the limit
of vanishing band widths.  \eqref{4.3.3} is just a kind of Abelianized
summation procedure applied to \eqref{4.3.4}.
 
As a special case of this example, consider $V(x)=x^{2}-1$.  Then
$E_{j}=2j$ and $\{\mu_{j}(0)\}$ is the set $\{2, 2, 6, 6, 10, 10, 14,
14, \dots \}$ of $j$ odd eigenvalues, each doubled.  Thus \eqref{4.3.4} is the
formal sum
$$
-1=-2+2-2+2\dots \, \text{ (formal!)}
$$
with \eqref{4.3.3}
\begin{align*}
-1 &=\lim_{\alpha\downarrow 0}\,[(e^{-2\alpha}-1)/\alpha]\,[1-e^{-2\alpha}
+e^{-4\alpha}\dots ] \\
&=\lim_{\alpha\downarrow 0}\,[(e^{-2\alpha}-1)/\alpha]\,
[1/(1+e^{-2\alpha})]
\end{align*}
with a true abelian summation.
\end{example}

\begin{example} \lb{e4.3.4}  Suppose $V$ is periodic.  Let 
$E_{j}, \mu_{j}(x)$ be the band edges and Dirichlet eigenvalues as in
\eqref{4.1.2}, \eqref{4.1.3}.  Then it follows from the fact that the two 
Floquet solutions are complex conjugates of each other on the spectrum 
of $H$, and the Wronskian is antisymmetric in its argument ($W(f,g)=-W(g,f)$), 
that $g(\lambda,x)$ is purely imaginary on
$\sigma (H)$; that is,
$\xi(\lambda,x)=1/2$ there, so
$$
\xi(\lambda,x)= \begin{cases} 1/2, & E_{2n}<\lambda<E_{2n+1},
\\ 1, & E_{2n-1}<\lambda<\mu_{n}(x), \\
0, & \mu_{n}(x)<\lambda<E_{2n}.
\end{cases}
$$
(This fact has also been used by Deift and Simon
\cite{DS83} and Kotani
\cite{Ko84} and also follows from the fact that $g(\lambda,x)=G(\lambda+i0,x,x)
=-[m_{+}(\lambda,x)+m_{-} (\lambda,x)]^{-1}$ in terms of the Weyl--Titchmarsh 
$m$-functions.) Thus, 
$$
\int^{\infty}_{E_0} d\lambda \, |1-2\xi(\lambda,x)|=\sum^{\infty}_{n=1}
|E_{2n}-E_{2n-1}|
$$
is finite if \eqref{4.1.1} holds.  In that case one can take the limit 
inside the integral in \eqref{4.3.1} and so recover \eqref{4.1.3}.
\end{example}

\begin{example} \lb{e4.3.5}  In \cite{GHS95} it is proved that if $V$ is
short-range, that is, $V\in H^{2, 1}(\bbR)$, then, $\int^{\infty}_{E_0}
d\lambda \, |1-2\xi(\lambda,x)| <\infty$ and one can take the
limit in \eqref{4.3.1} inside the integral.  This recovers Venakides'
result \cite{Ve88} with an explicit form for $\xi$ in terms of the Green's
function (see Theorem \ref{t4.3.1}). Similarly, one can treat short-range
perturbations $W$ of periodic background potentials $V$ (modeling
scattering off defects or impurities, described by $W$, in
one-dimensional solids) and ``cascading'' potentials, that is, potentials
approaching different spatial asymptotes sufficiently fast (cf.\
\cite{GHS95} for details). 
\end{example}

Next we mention a striking inverse spectral application of the trace
formula \eqref{4.3.1} to a celebrated two-spectra inverse spectral 
theorem due to Borg \cite{Bo46}:

\begin{theorem} \lb{t4.3.2}
Let $V\in L^1_{\loc} (\bbR)$ be real-valued and periodic. Let
$H=-\f{d^2}{dx^2}+V$ be the associated self-adjoint 
Schr\"odinger operator in $L^2(\bbR)$ and suppose that
$\sigma(H)=[E_0,\infty)$ for some $E_0\in\bbR$.  Then
$$
V(x)=E_0 \, \text{ for a.e.\ $x\in\bbR$}. 
$$
\end{theorem}
Given the trace formula \eqref{4.3.1} (observing the a.e.\ extension
noted after Theorem \ref{t4.3.1}) and using the fact that for all
$x\in\bbR$ and a.e.\ $\lambda > E_0$, $\xi(\lambda,x)=1/2$ (cf.\ Example
\ref{4.3.4}), the proof of Borg's Theorem \ref{t4.3.2} is effectively
reduced to just a one-line argument (as was observed in \cite{CGHL00}).
In addition, the new proof permits one to replace periodic by
reflectionless potentials and hence applies to algebro-geometric
quasiperiodic (KdV) potentials and certain classes of almost periodic
potentials.

Now we turn to an analog for Theorem \ref{t4.3.1} for
Jacobi operators. This turns out to be a special case of the following
result.
\begin{theorem} \lb{t4.4.1}  Let $A$ be a bounded self-adjoint operator
in some complex separable Hilbert space $\cH$ with $\alpha=\inf\,
\sigma(A)$, $\beta=\sup\,\sigma(A)$.  Let $\varphi\in\cH$,
$\|\varphi\|_{\cH}=1$ be an arbitrary unit vector in $\cH$ and let
$\xi(\lambda)$ be the spectral shift function for the pair $(A_{\infty},A)$,
where
$A_\infty$ is defined by
$$
(A_\infty-z)^{-1}=(A-z)^{-1}-(\varphi,(A-z)^{-1}\varphi)^{-1}
((A-\ol z)^{-1}\varphi,\cdot)(A-z)^{-1}\varphi.
$$
Then for any $E_{-}\leq\alpha$ and $E_{+}\geq\beta$:
\begin{align*}
(\varphi, A\varphi) &=E_{-}+\int^{E_{+}}_{E_{-}} d\lambda \, 
[1-\xi(\lambda)]   
= E_{+}-\int^{E_{+}}_{E_{-}} d\lambda \, \xi(\lambda)    \\
&=\frac12\, (E_{+}+E_{-})+\frac12 \int^{E_{+}}_{E_{-}} d\lambda \, 
[1-2\xi(\lambda)].  
\end{align*}
\end{theorem} 

\begin{corollary} \lb{c4.4.2}  Let $H$ be a Jacobi matrix on 
$\ell^{2}(\bbZ^{\nu})$, that is, for a bounded function $V$ on
$\bbZ^{\nu}$,
\begin{equation}
(Hu)(n)=\sum_{|n-m|=1} u(m)+V(n)u(n), \quad n\in\bbZ^{\nu}.   \lb{4.4.4}
\end{equation}
For $r\in\bbZ^{\nu}$, let $H^{D}_{r}$ be the operator on 
$L^{2}(\bbZ^{\nu}\backslash\{r\})$ given by \eqref{4.4.4} with $u(r)=0$
boundary conditions.  Let $\xi(\lambda,r)$ be the spectral shift function for
the pair
$(H^{D}_{r}, H)$.  Then
\begin{align}
\begin{split}
V(r) &= E_{-}+\int^{E_{+}}_{E_{-}} d\lambda \, [1-\xi(\lambda,r)]  
= E_{+}-\int^{E_{+}}_{E_{-}} d\lambda \, \xi(\lambda,r)  \\
&= \frac12\, (E_{+}+E_{-})+\frac12
\int^{E_{+}}_{E_{-}}  d\lambda \, [1-2\xi(\lambda,r)]   \lb{4.4.5}
\end{split}  
\end{align}
for any $E_{-}\leq\inf\,\sigma(H)$,
$E_{+}\geq\sup\,\sigma(H)$.
\end{corollary}
Only when $\nu=1$ does this have an interpretation as a formula using
Dirichlet problems on the half-line. 

Next, we look at some applications to absolutely continuous spectra. In
particular, we will point out that $\xi(\lambda,x)$ for a single fixed
$x\in\bbR$ determines the absolutely continuous spectrum of a
one-dimensional Schr\"odinger or Jacobi operator.  We begin with a
result that holds for higher-dimensional Jacobi operators as well:

\begin{lemma} \lb{l4.5.1} $(i)$ For an arbitrary Jacobi matrix,
$H$, on $\bbZ^{\nu}$, $\operatornamewithlimits{\cup}_
{j\in\bbZ^{\nu}}\{\lambda\in\bbR \,|\, 0<\xi(\lambda,j)<1\}$ is an
essential support for the absolutely continuous spectrum of $H$. \\
$(ii)$ For a one-dimensional Schr\"odinger operator,
$H=-\frac{d^2}{dx^2}+V$ bounded from below,
$\operatornamewithlimits{\cup}_{x\in\bbQ}\{\lambda
\in\bbR\,|\, 0<\xi(\lambda,x)<1\}$ is an essential support for the
absolutely continuous spectrum of $H$.
\end{lemma}

\begin{remark} \lb{r4.5.2} We recall that every absolutely continuous
measure, $d\mu$, has the form $f(E)\,dE$.  $S=\{E\in\bbR\,|\, f(E)\neq
0\}$ is called an essential support.  Any Borel set which differs from
$S$ by sets of zero Lebesgue measure is also called an essential support.
If $A$ is a self-adjoint operator on $\cH$ and $\varphi_{n}$, an
orthonormal basis for $\cH$, and $d\mu_{n}$, the spectral measure for
the pair, $A, \varphi_{n}$ (i.e., $(\varphi_{n}, e^{isA}\varphi_{n})=
\int_{\bbR} e^{isE}\,d\mu_{n}(E)$) and if $d\mu^{\rm ac}_{n}$
is the absolutely continuous component of $d\mu_{n}$ with $S_n$ its
essential support, then $\operatornamewithlimits{\cup}_{n}S_n$ is the
essential support of the absolutely continuous spectrum for $A$.
\end{remark}
  
In one dimension though, a single $x$ suffices:

\begin{theorem} \lb{t4.5.2}  For one-dimensional Schr\"odinger 
$($resp., Jacobi$)$ operators, $\{\lambda\in\bbR\,|\,
0<\xi(\lambda,x)<1\}$ is an essential support for the absolutely
continuous measure for \rm{any } fixed $x\in\bbR$ $($resp., $\bbZ$$)$.
\end{theorem}
 
These results are of particular interest because of their implications
for a special kind of semi-continuity of the spectrum.  
 
\begin{definition} \lb{d4.5.2} Let $\{V_{n}\},V$ be continuous potentials 
on $\bbR$ (resp., on $\bbZ$).  We say that $V_n$ converges to $V$ locally
as $n\to\infty$ if and only if \\
$(i)$  $\inf_{(n, x)\in\bbN\times\bbR}V_{n}(x)>-\infty$
($\bbR$ case) or $\sup_{(n, j)\in\bbN\times\bbZ}
|V_{n}(j)|<\infty$ ($\bbZ$ case). \\
$(ii)$ For each $R<\infty$, $\sup_{|x|\leq R} |V_{n}(x)
-V(x)|\to 0$ as $n\to\infty$.
\end{definition}

\begin{lemma} \lb{l4.5.3}  If $V_{n}\to V$ locally as $n\to\infty$ and
$H_n, H$ are the corresponding Schr\"odinger operators $($resp., Jacobi
matrices$)$, then $(H_{n}-z)^{-1}\to (H-z)^{-1}$ strongly for
$\text{Im}\,z\neq 0$ as $n\to\infty$.
\end{lemma} 
 
\begin{theorem} \lb{t4.5.4}  If $V_{n}\to V$ locally as $n\to\infty$ and
$\xi_{n}(\lambda,x)$, $\xi(\lambda,x)$ are the corresponding xi
functions for fixed $x$, then $\xi_{n}(\lambda,x)\, d\lambda$ converges
to
$\xi(\lambda,x)\, d\lambda$ weakly in the sense that for any $f\in
L^{1}(\bbR; d\lambda)$,
$$
\int_{\bbR} d\lambda \, f(\lambda) \, \xi_{n}(\lambda,x)\, \to
\int_{\bbR} d\lambda \, f(\lambda) \, \xi(\lambda,x)\,  \,
\text{ as } \, n\to\infty.    
$$
\end{theorem}

\begin{definition}  \lb{d4.5.4} For any $H$, let $|S_{\rm ac}(H)|$
denote the Lebesgue measure of the essential support of the absolutely
continuous spectrum of $H$.
\end{definition}

\begin{theorem} \lb{t4.5.5} $($For one-dimensional Schr\"odinger
or Jacobi operators$)$  Suppose $V_{n}\to V$ locally as
$n\to\infty$ and each $V_n$ is periodic.  Then for any interval
$(a, b)\subset\bbR$,
$$
|(a, b)\cap S_{\rm ac}|\geq \varlimsup_{n\to\infty}
|(a, b)\cap S_{\rm ac} (H_{n})|.
$$
\end{theorem}
We note that the periods of $V_{n}$ need {\it not} be fixed;
indeed, almost periodic potentials are allowed.
 
\begin{example} \lb{e4.5.6}  Let $\alpha_n$ be a sequence of rationals and
$\alpha=\lim_{n\to\infty}\,\alpha_n$.  Let $H_n$ be the Jacobi
matrix with potential $\lambda\cos (2\pi\alpha_{n}+\theta)$ for $\lambda,
\theta$ fixed.  Then [2] have shown for $|\lambda|\leq 2$, $|S_{n}|
\geq 4-2|\lambda|$.  It follows from Theorem \ref{t4.5.5} that $|S|\geq
4-2|\lambda|$, providing a new proof (and a strengthening) of a result of
Last \cite{La93}. At present much more is known about this example and the interested reader may want to consult the survey by Last \cite{La05} for additional results. 
\end{example}

\begin{example} \lb{e4.5.7}  Let $\{a_{m}\}_{m\in\bbN}$ be a sequence with
$s=\sum^{\infty}_{m=1}2^{m}|a_{m}|<2$. Let
$V(n)=\sum^{\infty}_{m=1}a_{m}\cos(2\pi n/2^{m})$, a limit periodic
potential on $\bbZ$.  Let $h$ be the corresponding Jacobi matrix, then
one can show that
$
|\sigma_{\rm ac}(h)|\geq 2(2-s).    
$
\end{example}

\smallskip
\begin{center}
$\bigstar$  \qquad $\bigstar$ \qquad $\bigstar$
\end{center}
\smallskip

Next we very briefly turn to higher-order trace formulas derived in \cite{GHSZ95} obtained by
higher-order expansions in \eqref{4.3.2} as $t\downarrow 0$. For
simplicity we now assume that $V\in C^\infty(\bbR)$ is bounded from
below. Then \eqref{4.3.2} can be extended to 
$$
\Tr[e^{tH_{x}^D}-e^{-tH}]
\operatornamewithlimits{\sim}\limits_{t\downarrow 0}
\sum^{\infty}_{j=0} s_{j}(x)t^{j}, \quad x\in\bbR.  
$$
Similarly, one has,
\begin{align*}
& \Tr[(H_{x}^D-z)^{-1}-(H-z)^{-1}]\operatornamewithlimits{\sim}
\limits_{z\downarrow -\infty} \sum^{\infty}_{j=0} r_{j}(x)z^{-j-1}, \\
& r_{0}(x)=1/2, \quad r_{1}(x)=V(x)/2, \quad x\in\bbR 
\end{align*}
and one can show that 
$$
s_{j}(x)=(-1)^{j+1}(j!)^{-1}r_{j}(x), \quad j\in\bbN\cup \{0\}.
$$
In particular, $r_j(x)$ and $s_{j}(x)$ are the celebrated KdV invariants
(up to inessential numerical factors). They can be computed recursively
(see, e.g., \cite{GHSZ95}). The higher-order analogs of \eqref{4.3.1}
then read
\begin{align*}
s_{0}(x) &=-\frac12 , \\
s_{j}(x) &=\frac{(-1)^{j+1}}{j!} \left\{\frac{E^{j}_{0}}{2}\right. + j \,
\lim_{t\downarrow 0} \int^{\infty}_{E_0}d\lambda\,
e^{-t\lambda}\lambda^{j-1} \biggl[\frac12 -\xi(\lambda, x)\biggr]
\quad j\in\bbN,\, x\in\bbR.
\end{align*}
and similarly using a resolvent rather than a heat kernel regularization,
\begin{align*}
&r_{1}(x) =\frac12 \,V(x) 
=\frac{E_0}{2}+\lim_{z\to i\infty}\int^{\infty}_{E_0}
d\lambda\, \frac{z^2}{(\lambda -z)^{2}}\biggl[\frac12 -\xi(\lambda,
x)\biggr], \\
&r_{j}(x)=\frac{E^{j}_{0}}{2}+\lim_{z\to
i\infty}\int^{\infty}_{E_0}
d\lambda\,\frac{z^{j+1}}{(\lambda-z)^{j+1}}\,j(-\lambda)^{j-1}
\biggl[\frac12 -\xi(\lambda, x)\biggr], \quad j\in\bbN,\,
x\in\bbR.
\end{align*}
In the special periodic case, the corresponding extension of
\eqref{4.1.3} then reads
$$
2(-1)^{j+1}j!\, s_{j}(x)=2r_{j}(x)=E^{j}_{0}+\sum^{\infty}_{n=1}
[E^{j}_{2n-1}+E^{j}_{2n}-2\mu_{n}(x)^{j}], \quad j\in\bbN,\, x\in\bbR. 
$$
The latter formulas were originally found in \cite{Fl75} and \cite{MM75}.

We also note that the use of the Dirichlet boundary boundary
condition $u(x)=0$ and hence the choice of the Dirichlet operator
$H_x^D$ in connection with \eqref{4.1.7} can be replaced by any
self-adjoint boundary condition of the type 
$u'(x_{\pm})+\beta u(x_{\pm})=0$, $\beta\in\bbR$, and the corresponding
Schr\"odinger operator $H_x^\beta$ in $L^{2}((-\infty, x))\oplus 
L^{2}((x,\infty))$. This is worked out in detail in \cite{GHSZ95}. 

Additional results on trace formulas for Schr\"odinger operators were presented in 
\cite{Ge95}, \cite{GH95}, \cite{GH97}, \cite{GM00}, \cite{GRT96}. 

\medskip

{\it More recent references:} Important extensions of the trace formula \eqref{4.3.1},
including the case of Schr\"odinger operators unbounded from below, were
discussed by Rybkin \cite{Ry01}, \cite{Ry01a}. Further discussions of the 
trace formula \eqref{4.4.5} for Jacobi operators can be found in
\cite{GS97a} and Teschl \cite{Te98}, \cite[Ch.\ 6]{Te00}. An extension of 
Corollary \ref{c4.4.2} to Schr\"odinger operators on a countable set was 
discussed by Shirai \cite{Sh98}.

Removal of the resolvent regularization limit in the above trace formula for 
$r_1$ (resp., $V$) under optimal conditions on $V$ has been studied by 
Rybkin \cite{Ry01a}, \cite{Ry03} (the latter reference offers necessary 
and sufficient conditions for absolute summability of the trace formula).

A certain multi-dimensional variant of these trace formulas, inspired by
work of Lax \cite{La94}, was established in \cite{GHSZ96} (see also
\cite{GH95}).

Matrix-valued extensions of the trace formula for Schr\"odinger,
Dirac-type, and Jacobi matrices, as well as Borg and Hochstadt-type
theorems were studied in \cite{BGMS03}, \cite{CG01}, \cite{CG02},
\cite{CGHL00}, \cite{CGR05}, \cite{GH97}, and \cite{GS03}. 

Trace formulas and an ensuing general Borg-type theorem for CMV operators
(i.e., in connection with orthogonal polynomials on the unit circle, cf.\
\cite{Si05}) appeared in \cite{GZ05}.

An application of $\xi$-function ideas to obtain Weyl-type
asymptotics using $\zeta$-function regularizations of determinants of
certain operators on complete Riemannian manifolds can be found in
Carron \cite{Ca02}. 

Theorem \ref{t4.5.2} was used in \cite{GKT96} to solve an inverse
spectral problem for Jacobi matrices and most recently in \cite{GY06} in
connection with proving purely absolutely continuous spectrum of a class
of reflectionless Schr\"odinger operators with homogeneous spectrum. It
has also recently been discussed in \cite[Sect.\ 1.5]{DK05}.

\section{Various Uniqueness Theorems in Inverse Spectral Theory} 
\label{s5}

In this section we summarize some of the principal results of
the following papers:  \\
\cite{DGS97} R.~del Rio, F.~Gesztesy, and B.~Simon, {\it Inverse
spectral analysis with partial information on the potential, III.
Updating boundary conditions}, Intl. Math. Research Notices {\bf 
1997,  No.\ 15}, 751--758. \\
\cite{DGS99} R.~del Rio, F.~Gesztesy, and B.~Simon, {\it Corrections and
Addendum to ``Inverse spectral analysis with partial information on the
potential, III. Updating boundary conditions''}, Intl. Math. Research
Notices {\bf 1999, No.\ 11}, 623--625.  \\
\cite{GS96a} F.~Gesztesy and B.~Simon, {\it Uniqueness theorems in inverse
spectral theory for one-dimensional Schr\"odinger operators},  
Trans. Amer. Math. Soc. {\bf 348}, 349--373 (1996). \\
\cite{GS97} F.~Gesztesy and B.~Simon, {\it Inverse spectral
analysis with partial information on the potential, I. The
case of an a.c.\ component in the spectrum}, Helv. Phys. Acta 
{\bf 70}, 66--71, 1997. \\
\cite{GS97a} F.~Gesztesy and B.~Simon, {\it $m$-functions and inverse
spectral analysis for finite and semi-infinite Jacobi matrices},
J. Analyse Math. {\bf 73}, 267--297 (1997). \\
\cite{GS99} F.~Gesztesy and B.~Simon, {\it On the determination of a
potential from three spectra}, in {\it Differential Operators and Spectral Theory}, 
V.\ Buslaev, M.\ Solomyak, and D.\ Yafaev (eds.), Amer. Math. Soc.
Transl. Ser.\ 2, {\bf 189}, 85--92 (1999). \\
\cite{GS00} F.~Gesztesy and B.~Simon, {\it Inverse spectral analysis with
partial information on the potential, II. The case of discrete spectrum},
Trans. Amer. Math. Soc. {\bf 352}, 2765--2787 (2000). 

\medskip\smallskip

One can argue that inverse spectral theory, especially, the case of uniqueness theorems in inverse spectral theory, started with the paper by Ambarzumian 
\cite{Am29} in 1929 and was turned into a full fledged discipline by the seminal 1946 paper by Borg \cite{Bo46}. Ambarzumian proved the special uniqueness theorem that if the eigenvalues of a Schr\"odinger operator $-d^2/dx^2 + V$ in $L^2([0,\pi])$ with Neumann boundary conditions at the endpoints $x=0$ and $x=\pi$ coincide with the sequence of numbers $n^2$, $n=0,1,2,\dots$, then $V=0$ a.e.\ on $[0,\pi]$. This result is very special. Indeed, Borg showed that for more general boundary conditions, one set of eigenvalues, in general (i.e., in the absence of symmetries of $V$), is insufficient to determine $V$ uniquely. Moreover, he described in great detail when two spectra  guarentee unique determination of the potential $V$. In this section we will discuss Borg's celebrated two-spectra uniqueness result and many of its extension due to Gasymov, Hald, Hochstadt, Levitan, Lieberman, Marchenko, and Simon and collaborators.  

We start with paper \cite{GS96a}. It contains a variety of
new uniqueness theorems for potentials $V$ in one-dimensional
Schr\"odinger operators $-\frac{d^2}{dx^2}+V$ on $\bbR$ and on 
the half-line $\bbR_{+}=[0,\infty)$ in terms of appropriate 
spectral shift functions introduced in a series of papers
describing new trace formulas for $V$ on $\bbR$
\cite{GHS95}, \cite{GHSZ93}, \cite{GHSZ95}, \cite{GS96} and on
$\bbR_{+}$ \cite{GH95}. In particular, it contains a generalization
of a well-known uniqueness theorem of Borg and Marchenko for Schr\"odinger
operators on the half-line with purely discrete spectra to arbitrary
spectral types and a new uniqueness result for Schr\"odinger operators
with confining potentials on the entire real line.
 
Turning to the half-line case first, we recall one of the principal 
uniqueness results proved in \cite{GS96a}, which extends a well-known
theorem of Borg \cite{Bo52} and Marchenko \cite{Ma73} in the special
case of purely discrete spectra to arbitrary spectral types. We suppose
\begin{equation}
V\in L^{1}([0,R]) \text{ for all $R>0$}, \quad V\text{ real-valued,} 
\lb{5.2.1}
\end{equation} 
and introduce the differential expression
$$
\tau_{+}=-\frac{d^2}{dx^2}+V(x), \quad x\geq 0,  
$$
for simplicity assuming that $\tau$ is in the limit point case at
$\infty$. (We refer to \cite{GS96a} for a general treatment that includes
the limit circle case.) Associated with
$\tau_+$ one introduces the following self-adjoint operator $H_{+,\alpha}$
in $L^{2}(\bbR_{+})$. 
\begin{align}
&H_{+,\alpha}f=\tau_{+}f, \quad \alpha\in [0,\pi), \no \\
&f\in\dom(H_{+,\alpha})=\{g\in L^{2}(\bbR_{+}) \,|\, g,g'\in
AC([0,R])\text{ for all } R>0;   \lb{5.2.2} \\
&  \hspace*{3.02cm}
\sin(\alpha)g'(0_{+})+\cos(\alpha)g(0_{+})=0; \, \tau_{+}g\in
L^{2}(\bbR_{+})\}.  \no
\end{align}
Then $H_{+,\alpha}$ has uniform spectral multiplicity one.
 
Next we introduce the fundamental system $\phi_{\alpha}(z,x)$,
$\theta_{\alpha}(z,x)$, $z\in\bbC$, of solutions of
$\tau_{+}\psi(z,x)=z\psi(z,x)$, $x\geq 0$, satisfying
$$
\phi_{\alpha}(z,0)=-\theta'_{\alpha}(z,0)=-\sin(\alpha), \quad
\phi'_{\alpha}(z,0)=\theta_{\alpha}(z,0)=\cos(\alpha)  
$$
such that $W(\theta_{\alpha}(z), \phi_{\alpha}(z))=1$. 
Furthermore, let
$\psi_{+,\alpha}(z,x)$, $z\in\bbC\backslash\bbR$ be the unique
solution of $\tau\psi(z)=\psi(z)$ which satisfies
$$
\psi_{+,\alpha}(z,\,\cdot\,)\in L^{2}(\bbR_{+}), \quad
\sin(\alpha)\psi'_{+,\alpha}(z,0_{+}) + \cos(\alpha) \psi_{+,\alpha}
(z,0_{+})=1.
$$ 
$\psi_{+,\alpha}$ is of the form
$$
\psi_{+,\alpha}(z,x)=\theta_{\alpha}(z,x)+m_{+,\alpha}(z)
\phi_{\alpha}(z,x)  
$$
with $m_{+,\alpha}(z)$ the half-line Weyl--Titchmarsh $m$-function. Being 
a  Herglotz function (i.e., an analytic function in the open upper-half
plane that maps the latter to itself)),
$m_{+,\alpha}(z)$ has the following representation in terms of a positive
measure $d\rho_{+,\alpha}$, 
$$
m_{+,\alpha} = \begin{cases} a_{+,\alpha} 
+\int_{\bbR} \biggl[\frac{1}{\lambda-z}-\frac{\lambda}
{1+\lambda^2}\biggr]\,d\rho_{+,\alpha}
(\lambda), \quad \alpha\in [0,\pi),   \\
\cot(\alpha)+\int_{\bbR}(\lambda-z)^{-1} d\rho_{+,\alpha}
(\lambda), \quad \alpha\in (0,\pi).  \end{cases}
$$
 
The basic uniqueness criterion for Schr\"odinger operators on the
half-line $[0,\infty)$, due to Marchenko \cite{Ma73}, that we shall rely
on repeatedly in the following, can be stated as follows.

\begin{theorem} \lb{t5.2.1} Suppose
$\alpha_{1}, \alpha_{2}\in [0,\pi)$, $\alpha_{1}\neq\alpha_{2}$ and
define $H_{+, j, \alpha_{j}}$, $m_{+, j, \alpha_{j}}$,
$\rho_{+,j,\alpha_{j}}$ associated with the differential expressions
$\tau_{j}=-\frac{d^2}{dx^2}+V_{j}(x)$, $x\geq 0$, where $V_{j}, j=1,2$
satisfy assumption \eqref{5.2.1}. Then the following three assertions are
equivalent:
\\
$(i)$ $m_{+,1,\alpha_{1}}(z)=m_{+,2,\alpha_{2}}(z)$, $z\in\bbC_{+}$. \\
$(ii)$ $\rho_{+,1,\alpha_{1}}((-\infty,\lambda])=
\rho_{+,2,\alpha_{2}} ((-\infty,\lambda])$, $\lambda\in\bbR$. \\
$(iii)$ $\alpha_{1}=\alpha_{2}$ and $V_{1}(x)=V_{2}(x)$ for 
a.e.~$x\geq 0$.
\end{theorem}

Next we relate Green's functions for different boundary conditions
at $x=0$. 

\begin{lemma} \lb{l5.2.2} Let $\alpha_{j}\in [0,\pi)$, $j=1,2$, $x,x'\in
\bbR_{+}$, and $z\in\bbC\backslash\{\sigma(H_{+,\alpha_{1}})\cup
\sigma(H_{+,\alpha_{2}})\}$. Then, 
\begin{align*}
&G_{+,\alpha_{2}}(z,x,x')-G_{+,\alpha_{1}}(z,x,x')=-
\frac{\psi_{+,\alpha_{1}}(z,x)\psi_{+,\alpha_{1}}(z,x')}
{\cot(\alpha_{2}-\alpha_{1})+m_{+, \alpha_{1}}(z)}, \\
&\frac{G_{+,\alpha_{2}}(z,0,0)}{G_{+,\alpha_{1}}(z,0,0)} =\frac{1}
{(\beta_{1}-\beta_{2})\sin^{2}(\alpha_{1})[\cot(\alpha_{2}-\alpha_{1})
+m_{+,\alpha_{1}}(z)]} \\
&\quad = (\beta_{1}-
\beta_{2})\sin^{2}(\alpha_{2})[\cot(\alpha_{2}-\alpha_{1})-
m_{+,\alpha_{2}}(z)],\quad \beta_{j}=\cot(\alpha_{j}), j=1,2  \\
& \Tr [(H_{+,\alpha_{2}}-z)^{-1}
-(H_{+,\alpha_{1}}-z)^{-1}] =
-\frac{d}{dz}\,\ln[\cot(\alpha_{2}-\alpha_{1})+m_{+,\alpha_{1}}(z)] \\
&\hspace*{5.42cm} =\frac{d}{dz}\,
\ln [\cot(\alpha_{2}-\alpha_{1})-m_{+,\alpha_{2}}(z)].
\end{align*}
\end{lemma}
 
Since $m_{+,\alpha}(z)$ is a Herglotz function, we may now introduce
spectral shift function [27] $\xi_{\alpha_{1},
\alpha_{2}}(\lambda)$ for the pair $(H_{+,\alpha_{2}}, H_{+,
\alpha_{1}})$ via the exponential Herglotz representation of
$m_{+,\alpha}(z)$ (cf.\ \cite{AD56})
\begin{align*}
& \cot(\alpha_{2}-\alpha_{1})+m_{+,\alpha_{1}}(z) 
=\exp\biggl\{\text{Re}[\ln
(\cot(\alpha_{2}-\alpha_{1})+m_{+,\alpha_{1}}(i))] \\
&\quad  + \int_{\bbR} \biggl[\frac{1}{\lambda-z}-
\frac{\lambda}{1+\lambda^2}\biggr]\xi_{\alpha_{1}, \alpha_{2}}
(\lambda)\,d\lambda\biggr\}, \quad 0\leq\alpha_{1} <\alpha_{2} <\pi,
z\in\bbC\backslash\bbR.
\end{align*}
This is extended to all $\alpha_{1},\alpha_{2}\in [0,\pi)$ by
$$
\xi_{\alpha,\alpha} (\lambda)=0, \quad \xi_{\alpha_{2}, \alpha_{1}}
(\lambda)=-\xi_{\alpha_{1}, \alpha_{2}}(\lambda)
\text{ for a.e.~$\lambda\in\bbR$}.  
$$
 
Next we summarize a few properties of $\xi_{\alpha_{1},\alpha_{2}}
(\lambda)$.

\begin{lemma} \lb{l5.2.3} $(i)$ Suppose $0\leq\alpha_{1} <\alpha_{2}
<\pi$. Then for a.e.~$\lambda\in\bbR$,
$$
\xi_{\alpha_{1},\alpha_{2}}(\lambda) = \begin{cases}
\lim\limits_{\epsilon\downarrow 0}\,\pi^{-1} \Im
\{\ln[\cot(\alpha_{2}-\alpha_{1})+m_{+,\alpha_{1}}(\lambda+
i\epsilon)]\}, \\
-\lim\limits_{\epsilon\downarrow 0}\,\pi^{-1}\Im \{\ln
[\cot(\alpha_{2}-\alpha_{1})-m_{+,\alpha_{2}}(\lambda+
i\epsilon)]\}, \\
\lim\limits_{\epsilon\downarrow 0}\,\pi^{-1}\Im \bigl\{\ln
\bigl[\tfrac{1}
{\sin(\alpha_{1})}\,\tfrac{G_{+,\alpha_{1}}(\lambda +
i\epsilon, 0, 0)}{G_{+,\alpha_{2}}(\lambda + i\epsilon, 0,0)}\bigr]
\bigr\}.  
\end{cases}
$$
$($For $\alpha_{1}=0$, 
$G_{+,\alpha_{1}}(\lambda+i\epsilon,0,0)/\sin
(\alpha_{1})$ has to be replaced by $-1$ in the last expression.$)$
Moreover,
$$
0\leq \xi_{\alpha_{1},\alpha_{2}}(\lambda)\leq 1 \, \text{ a.e.}
$$
$(ii)$ Let $\alpha_{j}\in [0,\pi)$, $1\leq j\leq 3$. Then the
``chain rule''
$$
\xi_{\alpha_{1},\alpha_{3}}(\lambda)=\xi_{\alpha_{1},
\alpha_{2}}(\lambda)
+\xi_{\alpha_{2},\alpha_{3}}(\lambda)  
$$
holds for a.e.~$\lambda\in\bbR$. \\
$(iii)$ For all $\alpha_{1},\alpha_{2}\in [0,\pi)$,
$$
\xi_{\alpha_{1},\alpha_{2}}\in L^{1}(\bbR; (1+\lambda^{2})^{-1}\, 
d\lambda).  
$$
$(iv)$ Assume $\alpha_{1}, \alpha_{2}\in [0,\pi)$,
$\alpha_{1}\neq \alpha_{2}$. Then,
$$
\xi_{\alpha_{1},\alpha_{2}}\in L^{1}(\bbR; (1+|\lambda|)^{-1}
d\lambda) \, \text{ if and only if } \, \alpha_{1},\alpha_{2}\in
(0,\pi).  
$$
$(v)$ For all $\alpha_{1},\alpha_{2}\in [0,\pi)$,
$$
\Tr [(H_{+,\alpha_{2}}-z)^{-1} -(H_{+,\alpha_{1}}-z)^{-1}]
=-\int_{\bbR} \f{d\lambda \, 
\xi_{\alpha_{1},\alpha_{2}} (\lambda)}{(\lambda-z)^{-2}}. 
$$
\end{lemma}
 
We note that $\xi_{\alpha_{1}, \alpha_{2}}(\lambda)$ (for $\alpha_{1},
\alpha_{2}\in (0,\pi)$) has been introduced by Javrjan \cite{Ja66},
\cite{Ja71}. In particular, he proved Lemma \ref{l5.2.2}\,$(iii)$ and
Lemma \ref{l5.2.3}\,$(v)$ in the non-Dirichlet cases where
$0<\alpha_{1},\alpha_{2}<\pi$. 
 
Given these preliminaries, we are now able to state the main
uniqueness result for half-line Schr\"odinger operators of \cite{GS96a}.

\begin{theorem} \lb{t5.2.4} Suppose $V_j$ satisfy assumption
\eqref{5.2.1} and define $H_{+,j,\alpha_{j,\ell}}$, 
$j,\ell=1,2$, associated with the differential expressions
$\tau_{j}=-\frac{d^2} {dx^2}+V_{j}(x)$, $x\geq 0$, $j=1,2$, where 
$\alpha_{j,\ell}\in [0,\pi)$,
$\ell=1,2$, and we suppose $0\leq \alpha_{1,1} <\alpha_{1,2} <\pi$, $0\leq
\alpha_{2,1} <\alpha_{2,2} <\pi$. In
addition, let $\xi_{j,\alpha_{j,1},\alpha_{j,2}}$, $j=1,2$ be the
spectral shift function for the pair $(H_{+,j,\alpha_{j,1}},
H_{+,j,\alpha_{j,2}})$. Then the following two assertions are equivalent: 
\\
$(i)$ $\xi_{1,\alpha_{1,1},\alpha_{1,2}}(\lambda) =
\xi_{2,\alpha_{2,1},\alpha_{2,2}}(\lambda)$ for a.e.~$\lambda\in\bbR$. \\
$(ii)$ $\alpha_{1,1}=\alpha_{2,1}$,
$\alpha_{1,2}=\alpha_{2,2}$, and $V_{1}(x)=V_{2}(x)$ for a.e.~$x\geq
0$.
\end{theorem}

As a corollary, one obtains a well-known uniqueness result 
originally due to Borg \cite[Theorem\ 1]{Bo52} and Marchenko 
\cite[Theorem\ 2.3.2]{Ma73} (see also \cite{LG64}). 

\begin{corollary} \lb{c5.2.5} Define
$\tau_j$, $H_{+,j,\alpha}$, $\alpha\in [0,\pi)$ as in Theorem
\ref{t5.2.4}.  Assume
in addition that $H_{+,1,\alpha_{1}}$ and $H_{+,2,\alpha_{2}}$ 
have purely discrete spectra for some $($and hence for all$)$ 
$\alpha_{j}\in [0,\pi)$, that is,
$$
\sigma_{\rm ess} (H_{+,j,\alpha_{j}})=\emptyset
\, \text{ for some } \, \alpha_{j}\in [0,\pi), j=1,2.  
$$
Then the following two assertions are equivalent: \\
$(i)$ $\sigma(H_{+,1,\alpha_{1,1}})=\sigma
(H_{+,2,\alpha_{2,1}})$, $\sigma(H_{+,1,\alpha_{1,2}})=
\sigma (H_{+,2,\alpha_{2,2}})$, 
$\alpha_{j,\ell}\in [0,\pi)$, $j,\ell=1,2$,  $\sin(\alpha_{1,1}-
\alpha_{1,2})\neq 0$. \\
$(ii)$ $\alpha_{1,1}=\alpha_{2,1}$,
$\alpha_{1,2}=\alpha_{2,2}$, and $V_{1}(x)=V_{2}(x)$ for 
a.e.~$x\geq 0$.
\end{corollary}
 
Roughly speaking, Corollary \ref{c5.2.5} implies that two sets of purely 
discrete spectra $\sigma(H_{+,\alpha_{1}}), \sigma(H_{+,\alpha_{2}})$
associated with distinct boundary conditions at $x=0$ (but a fixed
boundary condition (if any) at $+\infty$), that is, 
$\sin(\alpha_{2}-\alpha_{1})\neq 0$, uniquely determine $V$ a.e. The
first main  result in \cite{GS96a}, Theorem \ref{t5.2.4}, removes all a
priori spectral hypotheses and  shows that the spectral shift function
$\xi_{\alpha_{1}, \alpha_{2}}(\lambda)$ for the pair $(H_{+,\alpha_{2}},
H_{+,\alpha_{1}})$  with distinct boundary conditions at $x=0$, 
$\sin(\alpha_{2}-\alpha_{1})\neq 0$, uniquely determines $V$ a.e. This
illustrates that  Theorem \ref{t5.2.4} is the natural generalization of
Borg's and Marchenko's theorems  from the discrete spectrum case to
arbitrary spectral types.

Now we turn to uniqueness results for Schr\"odinger operators
on the whole real line. We shall rely on the notation $\tau$, 
$\phi_{\alpha}$, $\theta_\alpha$, $\psi_{\pm,\alpha}$, $m_{\pm,\alpha}$, 
$d\rho_{\pm,\alpha}$, which are defined in complete analogy to the
half-line case (with $x\in\bbR$), and we shall assume 
\begin{equation}
V\in L^{1}_{\text{loc}}(\bbR), 
\quad V\text{ real-valued}. \lb{5.3.1}
\end{equation}
Following \cite{GHSZ95}, we introduce, in
addition, the following family of self-adjoint operators $H^{\beta}_{y}$
in
$L^{2}(\bbR)$,
\begin{align*}
& H^{\beta}_{y}f=\tau f, \quad \beta\in\bbR\cup\{\infty\}, \quad
y\in\bbR, \\
& \dom (H^{\beta}_{y})
=\{g\in L^{2}(\bbR\,|\, g,g'\in AC([y,\pm
R])\text{ for all }R>0;\, g'(y_{\pm})+\beta g(y_{\pm})=0; \\
& \hspace*{3.7cm} \lim\limits_{R\to\pm\infty}\,W(f_{\pm}(z_{\pm}),g)(R)=0;
\, \tau g\in L^{2}(\bbR)\}.
\end{align*}
Thus $H^{D}_{y}:=H^{\infty}_{y}$ (resp., $(H^{N}_{y}:=H^{0}_{y})$
corresponds to the Schr\"odinger operator with an additional Dirichlet
(resp., Neumann) boundary condition at $y$. In obvious notation,
$H^{\beta}_{y}$ decomposes into the direct sum of half-line operators
$$
H^{\beta}_{y}=H^{\beta}_{-,y}\oplus H^{\beta}_{+,y}  
$$
with respect to
$ 
L^{2}(\bbR)=L^{2}((-\infty, y]) \oplus L^{2}([y,\infty)) 
$.
In particular, $H^{\beta}_{+,y}$ equals $H_{+,\alpha}$ for $\beta
=\cot(\alpha)$ (and $y=0$)) in our notation \eqref{5.2.2} and, as 
done in our previous Sections \ref{s2} and \ref{s4}, the reference
point $y$ will be added as a subscript to obtain $\theta_{\alpha,
y}(z,x)$, $\phi_{\alpha, y}(z,x)$, $\psi_{\pm,\alpha,y}(z,x)$, 
$m_{\pm,\alpha,y}(z)$, $M_{\alpha, y}(z)$, etc.  
 
Next, we recall a few results from \cite{GHSZ95}. With $G(z,x,x')$ and
$G^{\beta}_{y}(z,x,x')$ the Green's functions of $H$ and
$H^{\beta}_{y}$, one obtains (for $z\in\bbC\backslash
\{\sigma(H^{\infty}_{y})\cup\sigma(H)\}$)
\begin{align*}
G^{\beta}_{y}(z,x,x') &=G(z,x,x') 
-\frac{(\beta+\partial_{2})G(z,x,y)
(\beta+\partial_{1}) G(z,y,x')}{(\beta+\partial_{1})
(\beta+\partial_{2})G(z,y,y)}, \\
&\hspace*{3.4cm} \beta\in\bbR, z\in\bbC\backslash
\{\sigma(H^{\beta}_{y})\cup\sigma(H)\},  \\
G^{\infty}_{y}(z,x,x') &= G(z,x,x')-G(z,y,y)^{-1} G(z,x,y)
G(z,y,x'). 
\end{align*}
Here we abbreviated
\begin{align*}
& \partial_{1}G(z,y,x')
=\partial_{x}G(z,x,x') |_{x=y}, \quad \partial_{2}(G,z,x,y)
=\partial_{x'}G(z,x,x') |_{x'=y}, \\
&  \partial_{1}\partial_{2} G(z,y,y)
=\partial_{x}\partial_{x'}
G(z,x,x') |_{x=y=x'}, \, \text{ etc.} \\
& \partial_{1}G(z,y,x)=\partial_{2}G(z,x,y), \quad x\neq y.  
\end{align*}
As a consequence,
\begin{align*}
&\Tr [(H^{\beta}_{y}-z)^{-1}-(H-z)^{-1}] 
=-\frac{d}{dz}\,\ln [(\beta+\partial_{1})(\beta+\partial_{2})G(z,y,y)], 
\quad \beta\in\bbR,  \\
&\Tr [(H^{\infty}_{y}-z)^{-1}-(H-z)^{-1}] =-\frac{d}{dz}\,\ln
[G(z,y,y)].  
\end{align*}

In analogy to the Herglotz property of $G(z,y,y)$, 
$(\beta+\partial_{1})(\beta+\partial_{2}) G(z,y,y)$ is also Herglotz
for each $y\in\bbR$. Hence, both admit exponential representations of the
form
\begin{align*}
& G(z,y,y)=\exp\biggr\{c_{\infty}+\int_{\bbR} d\lambda\, 
\biggl[\frac
{1}{\lambda-z}-\frac{\lambda}{1+\lambda^2}\biggl] \xi^{\infty}
(\lambda, y)\biggr\}, \\
& c_{\infty}\in\bbR, \quad 0\leq\xi^{\infty}(\lambda, y)\leq 1
\text{ a.e.}, \\
& \xi^{\infty}(\lambda, y)
=\lim\limits_{\epsilon\downarrow 0}\, \pi^{-1}
\Im \{\ln [G(\lambda+i\epsilon, y, y)]\}
\text{ for a.e.~$\lambda\in\bbR$}, \\
& (\beta+\partial_{1})(\beta+\partial_{2})G(z,y,y)
=\exp\biggl\{c_{\beta}
+\int_{\bbR} d\lambda \, 
\biggl[\frac{1}{\lambda-z}-\frac{\lambda}
{1+\lambda^2}\biggr] [\xi^{\beta}(\lambda, y)+1]\biggr\}, \\
& c_{\beta}\in\bbR, \quad -1\leq\xi^{\beta}(\lambda, y)\leq 0
\text{ a.e.}, \quad \beta\in\bbR, \\
& \xi^{\beta}(\lambda, y)=
 \lim\limits_{\epsilon\downarrow 0}\, \pi^{-1}
\Im \{\ln [(\beta+\partial_{1})(\beta+\partial_{2})
G(\lambda+i\epsilon, y, y)]\}-1, \quad \beta\in\bbR 
\end{align*}
for each $y\in\bbR$. Moreover,
$$
\Tr [(H^{\beta}_{y}-z)^{-1}-(H-z)^{-1}]=-\int_{\bbR} d\lambda \, 
(\lambda-z)^{-2}\xi^{\beta}(\lambda, y), \quad
\beta\in\bbR\cup\{\infty\}. 
$$
  
Applying the basic uniqueness criterion for Schr\"odinger operators to
both half-lines $(-\infty,y]$ and $[y,\infty)$ then yields the following
principal characterization result for Schr\"odinger operators on
$\bbR$ first proved in \cite{GS96a}.
 
\begin{theorem} \lb{t5.3.2} Let $\beta_{1},\beta_{2}\in\bbR\cup
\{\infty\}$, $\beta_{1}\neq\beta_{2}$, and $x_{0}\in\bbR$. Then the following 
assertions hold: \\
$(i)$ $\xi^{\beta_1}(\lambda, x_{0})$ and $\xi^{\beta_2}
(\lambda, x_{0})$ for a.e.~$\lambda\in\bbR$ uniquely determine
$V(x)$ for a.e.~$x\in\bbR$ if the pair $(\beta_{1},\beta_{2})$
differs from $(0,\infty)$, $(\infty,0)$.  \\
$(ii)$ If $(\beta_{1},\beta_{2})=(0,\infty)$ or $(\infty,
0)$, assume in addition that $\tau$ is in the limit point case at
$+\infty$ and $-\infty$. Then $\xi^{\infty}(\lambda, x_{0})$ and
$\xi^{0}(\lambda, x_{0})$ for a.e.~$\lambda\in\bbR$ uniquely
determine $V$ a.e.~up to reflection symmetry with respect to $x_0$;
that is, both $V(x)$, $V(2x_{0}-x)$ for a.e.~$x\in\bbR$ correspond
to $\xi^{\infty}(\lambda, x_{0})$ and $\xi^{0}(\lambda, x_{0})$ for
a.e.~$\lambda\in\bbR$.
\end{theorem}

\begin{corollary} \lb{c5.3.3} Suppose $\tau$ is in the limit point case at
$+\infty$ and $-\infty$ and let $\beta\in\bbR\cup\{\infty\}$ and
$x_{0}\in\bbR$. Then $\xi^{\beta}(\lambda, x_{0})$ for
a.e.~$\lambda\in\bbR$ uniquely determines $V(x)$ for a.e.~$x\in\bbR$ 
if and only if $V$ is reflection symmetric with respect to $x_0$,
that is, $V(2x_{0}-x)=V(x)$ a.e.
\end{corollary}
 
In view of Corollary \ref{c5.2.5}, it seems appropriate to formulate
Theorem \ref{t5.3.2} in the special case of purely discrete spectra.

\begin{corollary} \lb{c5.3.4} Suppose $H$ $($and hence $H^{\beta}_{y}$
for all $y\in\bbR$, $\beta\in\bbR\cup\{\infty\}$$)$ has purely
discrete spectrum, that is, $\sigma_{\text{\rm{ess}}}(H)=\emptyset$
and let $\beta_{1},\beta_{2}\in\bbR\cup\{\infty\}$, $\beta_{1}\neq
\beta_{2}$, and $x_{0}\in\bbR$.  \\
$(i)$ $\sigma(H)$, $\sigma(H^{\beta_j}_{x_0})$, $j=1,2$
uniquely determine $V$ a.e.~if the pair $(\beta_{1},\beta_{2})$
differs from $(0,\infty)$ and $(\infty, 0)$.  \\
$(ii)$ If $(\beta_{1},\beta_{2})=(0,\infty)$ or $(\infty,
0)$, assume in addition that $\tau$ is in the limit point case at
$+\infty$ and $-\infty$. Then $\sigma(H)$, $\sigma(H^{\infty}_{x_0})$,
and $\sigma(H^{0}_{x_0})$ uniquely determine $V$ a.e.~up to reflection
symmetry with respect to $x_0$, that is, both $V(x)$ and $\widehat{V}
(x)=V(2x_{0}-x)$ for a.e.~$x\in\bbR$ correspond to $\sigma(H)=\sigma
(\widehat{H})$, $\sigma(H^{\infty}_{x_0})=\sigma(\widehat{H}^{\infty}
_{x_0})$, and $\sigma(H^{0}_{x_0})=\sigma(\widehat{H}^{0}_{x_0})$.
Here, in obvious notation, $\widehat H$, $\widehat{H}^{\infty}_{x_0}$,
$\widehat{H}^{0}_{x_0}$ correspond to $\widehat{\tau}=-\frac{d^2}
{dx^2}+\widehat{V}(x)$, $x\in\bbR$.   \\
$(iii)$ Suppose $\tau$ is in the limit point case at
$+\infty$ and $-\infty$ and let $\beta\in\bbR\cup\{\infty\}$. Then
$\sigma(H)$ and $\sigma(H^{\beta}_{x_0})$ uniquely determine $V$
a.e.~if and only if $V$ is reflection symmetric with respect to $x_0$.  \\
$(iv)$ Suppose that $V$ is reflection symmetric with
respect to $x_0$ and that $\tau$ is nonoscillatory at $+\infty$ and $-
\infty$. Then $V$ is uniquely determined a.e.~by $\sigma(H)$ in the
sense that $V$ is the only potential symmetric with respect to $x_0$
with spectrum $\sigma(H)$.
\end{corollary}
  
Of course, Corollary \ref{c5.3.4}\,$(iii)$ is implied
by the result of Borg [5] and Marchenko [32] (see Corollary \ref{c5.2.5} 
with $\alpha_{1}=0$, $\alpha_{2}=\pi/2$).
 
Thus far, we exclusively dealt with $\xi$-functions and spectra in
connection with uniqueness theorems. A variety of further uniqueness
results can be obtained by invoking alternative information such as
the left/right distribution of $\lambda^{\beta}_{n}(x_{0})$ (i.e.,
whether $\lambda^{\beta}_{n}(x_{0})$ is an eigenvalue of $H^{\beta}_
{-,x_{0}}$ in $L^{2}((-\infty, x_{0}])$ or of $H^{\beta}_{+,x_{0}}$ in
$L^{2}([x_{0},\infty))$) and/or associated norming constants. For details
we refer to the discussion in \cite{GS96a}. 

\medskip

{\it More recent references:} Uniqueness theorems related to Theorem \ref{t5.2.4} in the short-range case with spectral shift data replaced by scattering data were studied by Aktosun and Weder \cite{AW06}, \cite{AW07}. Analogs of Corollaries \ref{c5.2.5} and \ref{c5.3.4}\,$(i)$ for Jacobi operators were derived by Teschl \cite{Te98}.

\smallskip
\begin{center}
$\bigstar$  \qquad $\bigstar$ \qquad $\bigstar$
\end{center}
\smallskip

Next we focus on \cite{GS00}, which discussed results where the discrete
spectrum (or partial information on the discrete spectrum) and partial
information on the potential $V$ of a one-dimensional Schr\"odinger
operator $H=-\frac{d^2}{dx^2}+q$ determines the potential completely.
Included are theorems for finite intervals and for the whole line. In
particular, a new type of  inverse spectral problem involving fractions
of the eigenvalues of $H$ on a finite interval and knowledge of $V$ over
a corresponding fraction of the interval was posed and solved in
\cite{GS00}. The methods employed in \cite{GS00} rest on Weyl--Titchmarsh 
$m$-function techniques starting with the basic Borg--Marchenko uniqueness 
result (cf.\ Theorems \ref{t6.1.5} and \ref{t6.1.6}) and densities of zeros of a class of entire functions since the $m$-functions are meromorhic functions in this context.

In 1978, Hochstadt and Lieberman \cite{HL78} proved the following
remarkable theorem:

\begin{theorem} \lb{t6.1.1} Let $h_0 \in \bbR$, $h_1 \in \bbR \cup
\{\infty\}$ and assume $V_1, V_2 \in L^1 ((0,1))$ to be real-valued.
Consider the Schr\"odinger operators $H_1, H_2$ in $L^2 ([0,1])$ given by
$$
H_j = -\frac{d^2}{dx^2}+ V_j, \quad j=1,2,
$$
with the boundary conditions
\begin{align}
\begin{split}
u'(0) + h_0 u(0) & = 0, \lb{6.1.1} \\
u'(1) + h_1 u(1) & = 0.  
\end{split}
\end{align}
Let $\sigma(H_j)=\{\lambda_{j,n}\}$ be the $($necessarily simple$)$
spectra of $H_j, j=1,2$. Suppose that $V_1 = V_2$ a.e.\ on
$[0,1/2]$ and that $\lambda_{1,n} = \lambda_{2,n}$ for all $n$. Then
$V_1 = V_2$ a.e.\ on $[0,1]$.
\end{theorem}
Here, in obvious notation, $h_1 =\infty$ in \eqref{6.1.1} singles out
the Dirichlet boundary condition $u(1)=0$.

For each $\varepsilon >0$, there are simple examples  where $V_1 = V_2$ on
$[0,(1/2) - \varepsilon]$ and $\sigma (H_1)$ = $\sigma (H_2)$ but $V_1
\neq V_2$. (Choose $h_0=-h_1$, $V_1 (x)=0$ for $x \in (0,(1/2) -
\varepsilon] \cup [1/2, 1]$ and nonzero on $((1/2) - \varepsilon,
1/2)$, and $V_2 (x) = V_1 (1-x)$. See also Theorem I$^\prime$ in the
appendix of \cite{Su86}.)

Later refinements of Theorem \ref{t6.1.1} in \cite{Ha80}, \cite{Su86} 
(see also the summary in \cite{Su82}) showed that the boundary
condition for $H_1$ and
$H_2$ at $x=1$  need not be assumed a priori to be the same, and that if
$V$ is continuous,  then one only needs $\lambda_{1,n} = \lambda_{2,m(n)}$
for all values of $n$ but one.  The same boundary condition for $H_1$
and  $H_2$ at $x=0$, however, is crucial for Theorem~1.1 to hold (see
\cite{Ha80}, \cite{De90}).

Moreover, analogs of Theorem \ref{t6.1.1} for certain Schr\"odinger
operators are considered in \cite{Kh84} (see also \cite[Ch.\ 4]{PT87}).
Reconstruction techniques for $V$ in this context are discussed in
\cite{RS92}.  

Our purpose in \cite{GS00} was to provide a new approach to Theorem
\ref{t6.1.1} that we felt was more transparent and, moreover, capable of
vast generalizations. To state our generalizations, we will introduce a
shorthand notation to paraphrase Theorem \ref{t6.1.1} by saying ``$V$ on
$[0, 1/2]$ and the eigenvalues of $H$ uniquely determine $V$.'' This
is just a shorthand notation for saying $V_1 =V_2$ a.e.\ if the  obvious
conditions hold.

Unless explicitly stated otherwise, all potentials $V, V_1$, and $V_2$ 
will be real-valued and in $L^1 ([0,1])$ for the remainder of this paper.
Moreover, to avoid too many case distinctions we shall assume $h_0 ,
h_1 \in \bbR$ in \eqref{6.1.1}. In particular, for
$h_0, h_1 \in \bbR$ we index the corresponding eigenvalues $\lambda_n$
of $H$ by $n \in \bbN_0 = \bbN \cup \{0\}$. The case of Dirichlet
boundary conditions, where $h_0 = \infty$ and/or $h_1 = \infty$ has been 
dealt with in detail in \cite[Appendix\ A]{GS00}.

Here is a summary of the generalizations proved for
Schr\"odinger operators on $[0,1]$ in \cite{GS00}:

\begin{theorem} \lb{t6.1.2} Let $H=-\frac{d^2}{dx^2}+ V$ in $L^2 ([0,1])$
with boundary conditions \eqref{6.1.1} and $h_0,h_1 \in \bbR$. Suppose
$V$ is $C^{2k}(((1/2) - \varepsilon, (1/2) + \varepsilon))$ for some
$k=0,1,\dots$\, and for some $\varepsilon > 0$. Then $V$ on
$[0,1/2]$,
$h_0$, and all the eigenvalues of $H$ except for $(k+1)$ uniquely
determine $h_1$ and $V$ on all of $[0,1]$.
\end{theorem}

\begin{remark} $(i)$. The case $k=0$ in Theorem \ref{t6.1.2} is due to
Hald \cite{Ha80}. \\
$(ii)$ In the non-shorthand form of this theorem (cf.\ the paragraphs preceding Theorem 
\ref{t6.1.2}), we mean that both $V_1$
and $V_2$ are $C^{2k}$ near $x=1/2$. \\
$(iii)$ One need not know which eigenvalues are missing. Since the
eigenvalues asymptotically satisfy
$$
\lambda_n = (\pi n)^2 + 2(h_1 - h_0) + \int_0^1 dx\, V(x) +o(1) \,
\text{ as } \, n \to \infty,  
$$
given a set of candidates for the spectrum, one can tell how many are 
missing. \\
$(iv)$ For the sake of completeness we mention the precise
definition of $H$ in
$L^2 ([0,1])$ for real-valued $V \in L^1([0,1])$ and boundary condition
parameters
$h_0, h_1 \in \bbR \cup \{\infty \}$ in (1.1):
\begin{align}
H=&-\frac{d^2}{dx^2} + V, \no \\
\dom (H)=&\{g \in L^2([0,1]) \, | \, g,g' \in AC([0,1]); \, (-g'' + Vg)
\in L^2([0,1]);  \lb{6.1.3} \\
& \hspace*{3.1cm} g'(0) + h_0g(0)=0, \, g'(1) + h_1g(1) = 0 \},
\no
\end{align}
where $AC([0,1])$ denotes the set of absolutely continuous functions on
$[0,1]$
and $h_{x_0}=\infty$ represents the Dirichlet boundary condition 
$g(x_0)=0$ for $x_0 \in \{0,1\}$ in \eqref{6.1.3}.
\end{remark}

By means of explicit examples, it has been shown in Section\ 3 of
\cite{GS00}, that Theorem \ref{t6.1.2} is optimal in the sense that if $V$
is only assumed to be $C^{2k-1}$ near
$x=1/2$ for some $k \geq 1$, then it is not uniquely determined by
$V\restriction [0,1/2]$ and all the
eigenvalues but $(k+1)$.

Theorem \ref{t6.1.2} works because the condition that $V$ is $C^{2k}$ near
$x=1/2$ gives us partial information about $V$ on $[1/2, 1]$; indeed, we
know the values of 
$$
V(1/2), V'(1/2), \dots, V^{(2k)}(1/2)
$$ 
computed on $[1/2, 1]$ since one can compute them on $[0,1/2]$. This
suggests that knowing $V$ on more than $[0,1/2]$ should let one dispense
with a finite density of eigenvalues. That this is indeed the case is the
content of the following theorem. (We denote by $\#\{\cdots\}$ the cardinality of 
the set $\{\cdots\}$.)

\begin{theorem} \lb{t6.1.3} Let $H = -\frac{d^2}{dx^2}+V$ in $L^2 ([0,1])$
with boundary conditions \eqref{6.1.1} and $h_0,h_1 \in \bbR$. Then $V$
on $[0,(1+\alpha)/2]$ for some $\alpha \in (0,1)$, $h_0$, and a 
subset $S \subseteq \sigma(H)$ of all the eigenvalues $\sigma(H)$ of $H$
satisfying
\begin{equation}
\#\{\lambda\in S \,|\, \lambda\leq \lambda_0 \} \geq (1-\alpha)
\# \{\lambda\in \sigma(H) \,|\, \lambda\leq \lambda_0 \} + (\alpha/2) 
\lb{6.1.4}
\end{equation}
for all sufficiently large $\lambda_0 \in \bbR$, uniquely determine 
$h_1$ and $V$ on all of $[0,1]$.
\end{theorem}

\begin{remark} $(i)$ As a typical example, knowing slightly more than half
the eigenvalues and knowing $V$ on $[0,\frac34]$ determines $V$ uniquely
on all of $[0,1]$. To the best of our knowledge, Theorem \ref{t6.1.3}
introduced and solved a new type of inverse
spectral problem. \\
$(ii)$ As in the case $\alpha =0$, one has an extension of the same type
as Theorem \ref{t6.1.2}. Explicitly, if $V$ is assumed to be $C^{2k}$ near
$x=(1 +\alpha)/2$, we only need
$$
\# \{\lambda\in S \,|\, \lambda\leq \lambda_0 \} \geq (1-\alpha)
\# \{\lambda\in \sigma(H) \,|\, \lambda\leq \lambda_0 \} + (\alpha/2) -
(k+1)  
$$
instead of \eqref{6.1.4}.
\end{remark}

One can also derive results about problems on all of $\bbR$. 

\begin{theorem} \lb{t6.1.4} Suppose that $V\in L^1_{\loc}
(\bbR)$ satisfies the following two conditions: \\
$(i)$ $V(x) \geq C|x|^{2+\varepsilon} -D$ for some
$C,\varepsilon, D>0$.  \\
$(ii)$ $V(-x) \geq V(x) \quad x\geq 0$. \\
Then $V$ on $[0,\infty)$ and the spectrum of $H = -\frac{d^2}{dx^2} + V$ 
in $L^2 (\bbR)$ uniquely determine $V$ on all of $\bbR$.
\end{theorem}

Hochstadt-Lieberman \cite{HL78} used the details of the inverse
spectral theory in their proof. In a sense, we only used in \cite{GS00}
the main uniqueness theorem of that theory due to Marchenko
\cite{Ma73}, which we now describe. For $V \in L^1([a,b])$
real-valued, $-\infty <a <b <\infty$, consider $-u'' + Vu = zu$ 
with the boundary condition 
\begin{equation}
u'(b) + h_b u(b) = 0  \lb{6.1.7}
\end{equation}
at $x=b$. Let $u_+ (z,x)$ denote the solution of this equation, normalized,
say, by
$u_+ (z,b) = 1$. The $m_+$-function is then defined by
\begin{equation}
m_+ (z,a) = \frac{u'_+ (z,a)}{u_+ (z,a)}\, . \lb{6.1.8}
\end{equation}
Similarly, given a boundary condition at $x=a$, 
\begin{equation}
u' (a) + h_a u(a) = 0,  \lb{6.1.9}
\end{equation} 
we define the solution $u_- (z,x)$ of $-u'' + Vu = zu$ normalized by
$u_- (z,a)=1$  and then define
\begin{equation}
m_- (z,b) = \frac{u'_- (z,b)}{u_- (z,b)}\, . \lb{6.1.10}
\end{equation}
In our present context where $-\infty < a < b < \infty$, $m_{\pm}$ are
even meromorphic on $\bbC$). Moreover,
$$
\text{Im}\, (z) > 0 \, \text{ implies } \, \Im (m_- (z,b)) < 0,
\, \Im (m_+ (z,a)) > 0.  
$$

Marchenko's \cite{Ma73} fundamental uniqueness theorem of inverse
spectral theory then reads as follows:

\begin{theorem} \lb{t6.1.5} $m_+ (z,a)$ uniquely determines $h_b$ as well
as $V$ a.e.\ on $[a,b]$.
\end{theorem}

If $V \in L^1_{\loc}([a,\infty))$ is real-valued (with
$|a|<\infty$) and
$-\frac{d^2}{dx^2} + V$ is in the limit point case at infinity, one can
still define a unique
$m_+ (z,a)$ function but now for $Im (z) \neq 0$ rather than all
$z \in \bbC$. For such $z$, there is a unique function $u_+(z,
\,\cdot\,)$ which is
$L^2$ at infinity (unique up to an overall scale factor which drops out of
$m_+ (z,a)$ defined by \eqref{6.1.8}). Again, one has the following
uniqueness result independently proved by Borg \cite{Bo52} and
Marchenko \cite{Ma73}.

\begin{theorem} \lb{t6.1.6} $m_+ (z,a)$ uniquely determines $V$
a.e.\ on $[a,\infty)$.
\end{theorem}

It is useful to have $m_-(z,b)$ because of the following basic fact:

\begin{theorem} \lb{t6.1.7} Let $H = -\frac{d^2}{dx^2} + V$ be a
Schr\"odinger operator in $L^2 ([a,b])$ with boundary conditions
\eqref{6.1.7}, \eqref{6.1.9} and let $G(z,x,y)$ be the integral kernel 
of $(H-z)^{-1}$. Suppose $c \in (a,b)$ and let $m_+ (z,c)$ be the
corresponding $m_+$-function for $[c,b]$ and $m_-(z,c)$ the
$m_-$-function for $[a,c]$. Then
\begin{equation}
G(z,c,c) = \frac{1}{m_- (z,c) - m_+ (z,c)}\, . \lb{6.1.12}
\end{equation}
\end{theorem}

Theorems \ref{t6.1.5} and \ref{t6.1.6} are deep facts; Theorem
\ref{t6.1.7} is an elementary calculation following from the explicit
formula for the integral kernel of $(H-z)^{-1}$,
$$
G(z,x,y) = \frac{u_- (z,\min (x,y)) u_+ (z, \max (x,y))}{W (u_- (z), u_+
(z))(x)},
$$
where as usual 
$
W(f,g)(x) = f' (x) g(x) - f(x)g'(x)
$
denotes the Wronskian of $f$ and $g$. An analog of Theorem \ref{t6.1.7}
holds in case $[a,b]$ is replaced by $(-\infty,\infty)$.

We can now describe the strategy of our proofs of Theorems
\ref{t6.1.1}--\ref{t6.1.4}. $G(z,c,c)$ has poles at the eigenvalues of $H$
(this is not quite true; see below), so by \eqref{6.1.12}, at eigenvalues
$\lambda_n$ of $H$:
\begin{equation}
m_+ (\lambda_n,c) = m_- (\lambda_n,c). \lb{6.1.13}
\end{equation}
If we know $V$ on a left partial interval $[a,c]$ and we know some 
eigenvalue $\lambda_n$, then we know $m_- (z,c)$ exactly; so by
\eqref{6.1.13}, we know the value of $m_+ (\lambda_n,c)$ at the point
$\lambda_n$. Below we indicate when knowing the values of
$f(\lambda_n)$ of an analytic function of the type of the
$m$-functions uniquely determines $f(z)$. If $m_+(z,c)$ is determined, 
then by Theorem \ref{t6.1.5}, $V$ is determined on $[a,b]$ and so is
$h_b$.

So the logic of the argument for a theorem like Theorem \ref{t6.1.1} is
the following: \\
$(i)$ $V$ on $[0,1/2]$ and $h_0$ determine $m_- (z,1/2)$
by direct spectral theory.  \\
$(ii)$ The $\lambda_n$ and \eqref{6.1.13} determine $m_+ (\lambda_n,1/2)$, and
then by suitable theorems in complex analysis, $m_+ (z, 1/2)$ is uniquely
determined for all $z$.  \\
$(iii)$ $m_+ (z,1/2)$ uniquely determines $V$ (a.e.) on
$[1/2,1]$ and $h_1$ by inverse spectral theory.

It is clear from this approach why $h_0$ is required and $h_1$ is free in 
the context of Theorem \ref{t6.1.1} (see \cite{De90} for examples where
$h_1$ and $V \restriction [0,1/2]$ do not determine $V$); without
$h_0$ we cannot compute $m_-(z, 1/2)$ and so start the process.

As indicated before \eqref{6.1.13}, $G(z,c,c)$ may not have a pole at an
eigenvalue $\lambda_n$ of $H$. It will if $u_n (c) \neq 0$, but if $u_n
(c) =0$, then $G(z,c,c) = 0$ rather than $\infty$. Here $u_n$ denotes the
eigenfunction of $H$ associated with the (necessarily simple) eigenvalue
$\lambda_n$. Nevertheless, \eqref{6.1.13} holds at points where $u_n (c)
=0$ since then $u_- (c)=u_+ (c) = 0$, and so both sides of \eqref{6.1.13}
are infinite. (In spite of \eqref{6.1.13}, $m_- - m_+$ is also infinite at
$z = \lambda_n$ and so $G(\lambda_n, c, c) = 0$.) We summarize this 
discussion next:

\begin{theorem} \lb{t6.1.8} For any $c\in (a,b)$, \eqref{6.1.13} holds at
any eigenvalue $\lambda_n$ of $H_{[a,b]}$ $($with the possibility of
both sides of \eqref{6.1.13} being infinite$)$.
\end{theorem}

\medskip

{\it More recent references:} A new inverse nodal problem was reduced to Theorem \ref{t6.1.3} by Yang \cite{Ya01}. A substantial generalization of Theorem \ref{t6.1.4}, replacing condition $(i)$ by $H$ being bounded from below with purely discrete spectrum, was proved by Khodakovsky \cite{Kh99}, \cite{Kh00}. He also found other variants of Theorem \ref{t6.1.4}. 

\smallskip
\begin{center}
$\bigstar$  \qquad $\bigstar$ \qquad $\bigstar$
\end{center}
\smallskip

We end our survey of \cite{GS00} by briefly indicating the
uniqueness theorems for entire functions needed in the proofs of Theorems
\ref{t6.1.1}--\ref{t6.1.4}. In discussing extensions of Hochstadt's
discrete (finite matrix) version \cite{Ho79} of the
Hochstadt--Lieberman theorem in \cite{GS97a}, we made use of the
following simple lemma which is an elementary consequence of the fact
that any polynomial of degree $d$ with $d+1$ zeros must be the zero
polynomial:

\begin{lemma} \lb{t6.B.1} Suppose $f_1 = \frac{P_1}{Q_1}$ and $f_2 =
\frac{P_2}{Q_2}$ are two rational fractions where the polynomials satisfy
$\deg (P_1) = \deg (P_2)$ and $\deg (Q_1)=\deg (Q_2)$. Suppose that $d =
\deg (P_1) + \deg (Q_1)$ and that $f_1 (z_n) = f_2 (z_n)$ for $d+1$ distinct
points $\{z_n\}^{d+1}_{n=1}\in \bbC$. Then $f_1 = f_2$.
\end{lemma}

In the context of \cite{GS00}, one is interested in entire functions of the
form
\begin{equation}
f(z) = C \prod^\infty_{n=0} \biggl( 1- \frac{z}{x_n}\biggr), \lb{6.B.1}
\end{equation}
where $0<x_0 < x_1 < \cdots$ is a suitable sequence of positive numbers which
are the zeros of $f$ and $C$ is some complex constant.

Given a sequence $\{x_n\}_{n=0}^{\infty}$ of positive reals, we define
$$
N(t) = \#\{n \in \bbN \cup \{0\}\,|\, x_n < t\}. 
$$
We recall the following basic theorem (see, e.g., 
\cite[Ch.\ I]{Le80}, \cite[Sects.\ II.48, II.49]{Ma85}):

\begin{theorem} \lb{t6.B.2} Fix $0 < \rho_0 < 1$. Then: \\
$(i)$ If $\{x_n\}^\infty_{n=0}$ is a sequence of positive reals with
\begin{equation}
\sum_{n=0}^{\infty} x^{-\rho}_n < \infty \, \text{ for all 
$\rho > \rho_0$,}
\lb{6.B.3}
\end{equation}
then the product in \eqref{6.B.1} defines an entire function $f$ with
\begin{equation}
|f(z)| \leq C_1 \exp (C_2 |z|^\rho ) \, \text{ for all $\rho > \rho_0$}.
\lb{6.B.4}
\end{equation}
$(ii)$ Conversely, if $f$ is an entire function satisfying \eqref{6.B.4}
with all its $($complex$)$ zeros on $(0, \infty)$, then its zeros
$\{x_n \}^\infty_{n=0}$ satisfy \eqref{6.B.3}, and $f$ has the canonical
product expansion \eqref{6.B.1}. \\

Moreover, \eqref{6.B.3} holds if and only if
\begin{equation}
N(t) \leq C |t|^\rho \, \text{ for all $\rho > \rho_0$}. \lb{6.B.5}
\end{equation}
\end{theorem}

Given this theorem, we single out the following definition.

\begin{definition} \lb{d6.B.3} A function $f$ is called of $m$-type if and
only if $f$ is an entire function satisfying \eqref{6.B.4} (of order $0<\rho<1$
in the usual definition) with all the zeros of $f$ on $(0, \infty)$.
\end{definition}

Our choice of ``$m$-type'' in Definition \ref{d6.B.3} comes from the fact
that in many cases we discuss in this paper, the $m$-function is a ratio
of functions of $m$-type. By Theorem \ref{t6.B.2}, $f$ in Definition
\ref{d6.B.3} has the form \eqref{6.B.1} and $N(t)$, which we will denote
as $N_f (t)$, satisfies \eqref{6.B.5}. 

\begin{lemma} \lb{l6.B.5} Let $f$ be a function of $m$-type. Then there
exists
a $0 < \rho < 1$ and a sequence $\{R_k\}_{k=1}^{\infty}$, $R_k \to \infty$ as
$k \to \infty$, so that
$$
\inf \{ |f(z)| \,|\, |z| = R_k\} \geq C_1 \exp (-C_2 R^{\rho}_k).
$$
\end{lemma}

\begin{lemma} \lb{l6.B.6} Let $F$ be an entire function that satisfies
the following two conditions: \\
$(i)$ $\sup_{|z|=R_k} |F(z)| \leq C_1 \exp (C_2 R^\rho_k)$ for
some $0 \leq \rho <1$, $C_1, C_2 >0$, and some sequence $R_k \to \infty$ as
$k \to \infty$. \\
$(ii)$ $\lim_{|x|\to\infty; x\in\bbR} |F(ix)| =0$. \\
Then $F \equiv 0$.
\end{lemma}

Lemmas \ref{l6.B.5} and \ref{l6.B.6} finally yield the following result. 

\begin{theorem} \lb{t6.B.4} Let $f_1 ,f_2 ,g$ be three functions of 
$m$-type so that the following two conditions hold: \\
$(i)$ $f_1 (z) = f_2 (z)$ at any point $z$ with $g(z) =0$.
$(ii)$ For all sufficiently large $t$,
$$
\max (N_{f_1} (t), N_{f_2} (t)) \leq N_g (t) -1.
$$
Then, $f_1 =f_2$.
\end{theorem}

\smallskip
\begin{center}
$\bigstar$  \qquad $\bigstar$ \qquad $\bigstar$
\end{center}
\smallskip

Refinements of the results of \cite{GS00} can be found in
\cite{DGS97}, \cite{DGS99}. Here we just mention the following facts.

\begin{theorem} \lb{6.3.1a} Let $H_1 (h_0), H_2(h_0)$ be associated with
two potentials $V_1, V_2$ on $[0,1]$ and two potentially distinct
boundary conditions $h^{(1)}_1, h^{(2)}_1 \in \bbR$ at $x=1$.
Suppose that $\{(\lambda_n, h^{(n)}_0)\}_{n\in\bbN_0}$ is a sequence 
of pairs with $\lambda_0 <\lambda_1 <\cdots \to \infty$ and $h^{(n)}_0
\in\bbR \cup\{\infty\}$ so that both $H_1 (h^{(n)}_0)$ and 
$H_2 (h^{(n)}_0)$ have eigenvalues at $\lambda_n$. Suppose that
$$
\sum^\infty_{n=0} \frac{(\lambda_n - 
\frac{1}{4} \pi^2 n^2)_+}{n^2}
<\infty 
$$
holds. Then $V_1 = V_2$ a.e.\ on $[0,1]$ and  $h^{(1)}_1 = h^{(2)}_1$.
\end{theorem}

This implies Borg's celebrated two-spectra uniqueness result \cite{Bo46}
(see also, \cite{Le68}, \cite{LG64}, \cite[Ch.\ 3]{Le87}, \cite{Ma73}):

\begin{corollary} \lb{c6.3.2a}  Fix $h^{(1)}_0, h^{(2)}_0 \in \bbR$.
Then all the eigenvalues of $H(h^{(1)}_0)$ and all the 
eigenvalues of $H(h^{(2)}_0)$, save one, uniquely determine 
$V$ a.e.\ on $[0,1]$.
\end{corollary}

It also implies the following amusing result: 

\begin{corollary} \lb{c6.3.3a} Let $h^{(1)}_0, h^{(2)}_0, h^{(3)}_0
\in \bbR$ and denote by $\sigma_j=\sigma(H(h^{(j)}_0))$ the spectra
of $H(h^{(j)}_0)$, $j=1,2,3$. Assume $S_j \subseteq \sigma_j$, $j=1,2,3$
and suppose that for all sufficiently large $\lambda_0 >0$ one has
$$
\#\{\lambda \in \{S_1 \cup S_2 \cup S_3\} \, \text{ with } \, 
\lambda\leq \lambda_0 \} \geq \tfrac23 \#\{\lambda \in
\{\sigma_1 \cup \sigma_2 \cup \sigma_3 \}  \, 
\text{ with } \, \lambda\leq \lambda_0 \} -1.
$$
Then $V$ is uniquely determined a.e.\ on $[0,1]$.
\end{corollary}

In particular, two-thirds of three spectra determine $V$. 

\medskip

{\it More recent references:} Further refinements of Corollary \ref{c6.3.3a}, involving $N$ spectra, were proved by 
Horv\'ath \cite{Ho01} (he also studies the corresponding analog for a Dirac-type operator). Optimal and nearly optimal conditions for a set of eigenvalues to determine the potential in terms of closedness properties of the exponential system corresponding to the known eigenvalues (implying Theorem \ref{6.3.1a} and a generalization thereof) were also derived by Horv\'ath \cite{Ho05}. For an interesting half-line problem related to this circle of ideas we also refer to Horv\'ath \cite{Ho06}. A variant of Theorem 
\ref{6.3.1a} was discussed by Ramm \cite{Ra99}, \cite{Ra00}. 
Hochstadt--Lieberman-type problems for Schr\"odinger operators including a reconstruction algorithm has been presented by L.\ Sakhnovich. 
The analog of the two-spectra result, Corollary \ref{c6.3.2a}, including a reconstruction algorithm, for a class of singular potentials has been discussed by Hryniv and Mykytyuk \cite{HM03a}, \cite{HM04} (see also \cite{HM05}). They also studied Hochstadt--Lieberman-type results for such a class of singular potentials in \cite{HM04a}. 
Hochstadt--Lieberman-type results for a class of Dirac-type operators relevant to the AKNS system were published by del Rio and Gr\'ebert \cite{DG01}. Borg- and Hochstadt--Lieberman-type inverse problems for systems including matrix-valued 
Schr\"odinger and Dirac-type equations, were studied in depth by M.\ Malamud 
\cite{Ma99}, \cite{Ma99a}, \cite{Ma05} \cite{Ma06}. He also studied Borg-type theorems for $n$th-order scalar equations \cite{Ma06a}. Borg- and Hochstadt--Lieberman-type inverse problems for matrix-valued Schr\"odinger operators were also studied by Shen \cite{Sh01}. He also considered Borg-type inverse problems for Schr\"odinger operators with weights \cite{Sh05}.

Additional results on determining the potential uniquely from spectra associated to
three intervals of the type $[0,1]$, $[0,a]$, and $[a,1]$ for some $a\in(0,1)$ (and similarly for whole-line problems with purely discrete spectra) can be found in \cite{GS99}. This has been inspired by work of Pivovrachik \cite{Pi99}, who also addressed the reconstruction algorithm from three spectra in the symmetric case $a=1/2$ (see also 
\cite{Pi99b}, \cite{Pi01}, \cite{Pi05}). He also considered the analogous Sturm--Liouville problem applicable to a smooth inhomogeneous partially damped string in \cite{Pi99a} and extended some of these results to Sturm--Liouville equations on graphs in \cite{Pi00}, \cite{Pi06}. Uniqueness and characterization problems for a class of singular Sturm--Liouville problems associated with three spectra were studied by Hryniv and Mykytyuk \cite{HM03}.  The reconstruction of a finite Jacobi matrix from three of its spectra was presented by Michor and Teschl \cite{MT04}.

\smallskip
\begin{center}
$\bigstar$  \qquad $\bigstar$ \qquad $\bigstar$
\end{center}
\smallskip

These results are related to two other papers: In \cite{GS97a}, we
considered, among other topics, analogs of Theorems \ref{t6.1.1} and
\ref{t6.1.3} for finite tri-diagonal (Jacobi) matrices, extending a result
in \cite{Ho79}. The approach there is very similar
to the current one except that the somewhat subtle theorems on zeros of 
entire functions in this paper are replaced by the elementary fact that a
polynomial of degree at most $N$ with $N+1$ zeros must be identically
zero. In \cite{GS97}, we consider results related to Theorem \ref{t6.1.4}
in that for Schr\"odinger operators on
$(-\infty, \infty)$, ``spectral'' information plus the potential on one of
the half-lines determine the potential on all of $(-\infty,\infty)$. In
that paper, we considered situations where there are scattering states for
some set of energies and the ``spectral'' data are given by a reflection
coefficient on a set of positive Lebesgue
measure in the a.c.\ spectrum of $H$. The approach is not as close to this
paper as is \cite{GS97a}, but $m$-function techniques (see also
\cite{GS96a}) are critical in all three papers.

\medskip

{\it More recent references:} For additional results on inverse scattering with partial information on the potential we refer to Aktosun and Papanicolaou \cite{AP03}, Aktosun and Sacks \cite{AS98}, Aktosun and Weder \cite{AW02},  and the references therein.

\smallskip
\begin{center}
$\bigstar$  \qquad $\bigstar$ \qquad $\bigstar$
\end{center}
\smallskip

We conclude this section by briefly describing some of the results in
\cite{GS97a}, where inverse spectral analysis for finite and
semi-infinite Jacobi operators $H$ was studied. While discussing a variety of topics
(including trace formulas), we also provided a new proof of a result of
Hochstadt \cite{Ho79} and its extension, which can be viewed as the
discrete analog of the Hochstadt and Lieberman result in \cite{HL78}.
Moreover, we solved the inverse spectral problem for $(\delta_n, 
(H-z)^{-1}\delta_n)$ in the case of finite Jacobi matrices. As mentioned
earlier, the tools we apply are grounded in $m$-function techniques.
 
Explicitly, \cite{GS97a} studied finite $N\times N$ matrices of the form:
$$
H=\begin{pmatrix}
b_1 & a_1 & 0    & 0 & \cdot & 
\cdot & \cdot \\
a_1 & b_2 & a_2 & 0 & \cdot & 
\cdot & \cdot \\
0    & a_2 & b_3 & a_3 & \cdot & 
\cdot & \cdot \\
\cdot & \cdot & \cdot & \cdot & 
\cdot & \cdot & \cdot \\
\cdot & \cdot & \cdot & \cdot & 
\cdot & \cdot & \cdot \\
\cdot & \cdot & \cdot & \cdot & 
0 & a_{N-1} & b_N
\end{pmatrix}  
$$
and the semi-infinite analog $H$ defined on
$$
\ell^2 (\bbN) =  \bigg\{ u = (u(1), u(2), \dots) \,\bigg|\,
\sum^\infty_{n=1} |u(n)|^2 < \infty \bigg\}
$$
given by
$$
(Hu)(n) = \begin{cases} a_n u(n+1) + b_n u(n) + a_{n-1} 
u(n-1), & n \geq 2,  \\
a_1 u(2) + b_1 u(1), &  n=1. \end{cases}
$$
In both cases, we assume $a_n, b_n \in\bbR$   
with $a_n >0$. To avoid inessential technical complications, 
we will only consider the case where $\sup_n [|a_n| + |b_n|] <\infty$ in
which  case $H$ is a map from $\ell^2$ to $\ell^2$ and defines a
bounded self-adjoint operator. In the semi-infinite case, we will set
$N=\infty$. It will also be useful to consider 
the $b$'s and $a$'s as a single sequence
 $b_1, a_1, b_2, a_2, \dots = c_1, c_2, \dots$, that is,
$$
c_{2n-1} = b_n, \quad c_{2n} = a_n, \quad n\in\bbN.  
$$

Concerning the recovery of a finite Jacobi matrix from parts of the matrix
and additional spectral information (i.e., mixed data), Hochstadt
\cite{Ho79} proved the following remarkable theorem.

\begin{theorem} \lb{t6.4.1a} Let $N \in \bbN$. Suppose that
$c_{N+1}, \dots, c_{2N-1}$ are known, as well as the eigenvalues
$\lambda_1, \dots,\lambda_N$ of $H$. Then $c_1, \dots, c_N$ are uniquely
determined.
\end{theorem}

The discrete Hochstadt--Lieberman-type theorem proved in \cite{GS97a}
reads as follows.

\begin{theorem} \lb{t6.4.2a} Suppose that 
$1 \leq j \leq N$ and $c_{j+1}, \dots, c_{2N-1}$ are known, as well as $j$
of the eigenvalues. Then $c_1, \dots, c_j$ are uniquely determined.
\end{theorem}

We emphasize that one need {\it{not}} know which of the $j$ eigenvalues
one has.
 
Borg \cite{Bo46} proved the celebrated theorem that the spectra for two
boundary  conditions of a bounded interval regular Schr\"odinger operator 
uniquely determine the  potential. Later refinements (see, e.g.,
\cite{Bo52}, \cite{Ho67a}, \cite{Le49}, \cite{Le68}, \cite{LG64}, \cite{Ma73}) 
imply that they even determine the two boundary conditions.
 
Next, we consider analogs of this result for a finite Jacobi
matrix. Such analogs were first considered by Hochstadt
\cite{Ho67}, \cite{Ho74} (see also \cite{BG78}, \cite{FH74}, \cite{GH84},
\cite{GW76}, \cite{Ha76}, \cite{Ho79}).  The results below are
adaptations of known results for the continuum or the semi-infinite case,
but the ability to determine parameters by counting sheds light on facts
like the one that the lowest eigenvalue in the Borg result is not needed
under certain circumstances.
 
Given $H$, an $N\times N$ Jacobi matrix, one defines $H(b)$ to be the
Jacobi matrix where all $a$'s and $b$'s are the same as $H$, except $b_1$
is replaced by $b_1 + b$, that is,
$$
H(b) = H + b(\delta_1, \, \cdot\, 
)\delta_1.  
$$

\begin{theorem} \lb{t6.5.1a} The eigenvalues $\lambda_1, \dots,
\lambda_N$ of $H$, together with $b$ and $N-1$ eigenvalues $\lambda 
(b)_1, \dots, \lambda (b)_{N-1}$ of $H(b)$, determine $H$ uniquely.
\end{theorem}
 
Again it is irrelevant which $N-1$ eigenvalues of the $N$ eigenvalues of
$H(b)$ are known.

\begin{theorem} \lb{t6.5.2a} The eigenvalues 
$\lambda_1, \dots, \lambda_N$ of $H$, together with the $N$ eigenvalues
$\lambda  (b)_1, \dots, \lambda (b)_N$ of some $H(b)$
$($with $b$ unknown$)$, determine $H$ and $b$.
\end{theorem}

\begin{remark} Since
$$
b= \text{Tr}(H(b)-H) = \sum^N_{j=1} 
(\lambda (b)_j - \lambda_j),
$$
we can a priori deduce $b$ from the $\lambda (b)$'s and $\lambda$'s and
so deduce Theorem \ref{t6.5.2a} from Theorem \ref{t6.5.1a}. We note that
the parameter counting  works out. In Theorem \ref{t6.5.1a}, $2n-1$
eigenvalues determine $2n-1$ parameters; and in Theorem \ref{t6.5.2a}, 
$2n$ eigenvalues determine $2n$ parameters.
\end{remark}
 
The basic inverse spectral theorem for finite Jacobi matrices shows that
 $(\delta_1, (H-z)^{-1} \delta_1)$ determines $H$ uniquely. In
\cite{GS97a} we considered $N \in \bbN$, $1\leq n \leq N$, and asked
whether $(\delta_n, (H-z)^{-1} \delta_n)$ determines $H$ uniquely. For
 notational convenience, we occasionally allude to $G(z,n,n)$ as the $n,n$
Green's function in the remainder of this section. The $n=1$ result can 
be summarized via:

\begin{theorem} \lb{6.6.1a} $(\delta_1, (H-z)^{-1}\delta_1)$ has the form
 $\sum^N_{j=1}\alpha_j (\lambda_j -z)^{-1}$ with $\lambda_1 <\cdots <
\lambda_N, \, \sum^N_{j=1} \alpha_j =1$ and each $\alpha_j >0$. Every
such  sum arises as the $1,1$ Green's function of an $H$ and of exactly
one such $H$.
\end{theorem}
 
For general $n$, define $\tilde n = \min (n, N+1 -n)$. Then the following
theorems were proved in \cite{GS97a}:
 
\begin{theorem} \lb{t6.6.2a} $(\delta_n, (H-z)^{-1}\delta_n)$ has the form
 $\sum^k_{j=1} \alpha_j (\lambda_j -z)^{-1}$ with $k$ one of $N, N-1,
\dots, N-\tilde n + 1$ and $\lambda_1< \cdots < \lambda_k, \,\sum^k_{j=1} 
\alpha_j =1$ and each $\alpha_j >0$. Every such sum arises
as the $n,n$ Green's function of at least one $H$.
\end{theorem}

\begin{theorem} \lb{t6.6.3a} If $k=N$, then precisely $\binom{N-1}{n-1}$
operators $H$ yield the given $n,n$ Green's function.
\end{theorem}

\begin{theorem} \lb{t6.6.4a} If $k<N$, then infinitely many Jacobi
matrices $H$ yield the given $n,n$ Green's function. Indeed, the inverse
spectral  family is then a collection of
$\binom{k-1}{N-k} \binom{k-1-N+k}{n-1-N+k}$ disjoint  manifolds, each of
dimension $N-k$ and diffeomorphic to an $(N-k)$-dimensional open ball.
\end{theorem}

\medskip

{\it More recent references:} Additional geometric information in connection with Theorem \ref{t6.6.4a} and a version for off-diagonal Green's functions were studied by Gibson \cite{Gi02}. 
Borg- and discrete Hochstadt--Lieberman-type results for generalized (i.e., certain 
tri-diagonal block) Jacobi matrices were studied by Derevyagin \cite{De06} (see also 
Shieh \cite{Sh04}). The case of non-self-adjoint Jacobi matrices with a rank-one imaginary part, and an extension of Hochstadt--Lieberman-type results to this situation was recently discussed by Arlinski\u i and Tsekanovski\u i \cite{AT06}.  An extension of results of Hochstadt \cite{Ho74} to the case of normal matrices was found by 
S.\ Malamud \cite{Ma04}. A detailed treatment of two-spectra inverse problems of semi-infinite Jacobi operators, including reconstruction, has recently been presented by Silva and Weder \cite{SW06}. 

\section{The Crown Jewel: Simon's New Approach to Inverse Spectral
Theory} 
\label{s6}
 
In this section we summarize some of the principal results of
the following papers:  \\
\cite{Si99} B.~Simon, {\it A new approach to inverse spectral theory, I. 
Fundamental formalism}, Ann. Math. {\bf 150}, 1029--1057 (1999). \\
\cite{GS00a} F.\ Gesztesy and B.\ Simon, {\it A new approach to 
inverse spectral theory, II.  General real potentials and the connection to
the spectral measure}, Ann. Math. {\bf 152}, 593--643 (2000). \\
\cite{RS00} A.\ Ramm and B.\ Simon, {\it A new approach to inverse spectral
theory, III. Short range potentials}, J. Analyse Math. {\bf 80},
319--334 (2000). \\
\cite{GS00b} F.~Gesztesy and B.~Simon, {\it On local Borg--Marchenko
uniqueness results}, Commun. Math. Phys. {\bf 211}, 273--287 (2000).

\medskip\smallskip

As the heading of this section suggests, we are approaching the pinnacle of
Barry Simon's contributions to inverse scattering theory thus far: In his
spectacular paper \cite{Si99}, he single-handedly developed a new approach
to inverse spectral theory for Schr\"odinger operators on a half-line, by
starting from a particular representation of the Weyl--Titchmarsh
$m$-function as a finite Laplace-type transform with control over the error term. In addition to establishing this feat,  it also led to a completely unexpected 
uniqueness result for Weyl--Titchmarsh functions, what is now called the local
Borg--Marchenko uniqueness theorem, but which really should have been
named Simon's local uniqueness theorem. The inverse spectral approach for
Schr\"odinger operators on a half-line (including a reconstruction algorithm for the potential) originated with the celebrated
paper \cite{GL55} by Gelfand and Levitan in 1951 and an independent approach by Krein \cite{Kr51} in the same year, followed by a seminal
contribution \cite{Ma73} by Marchenko in 1952. The
Borg--Marchenko uniqueness result was first published by Marchenko \cite{Ma50} in 1950 but Borg apparently had it in 1949 and it was independently published by Borg
\cite {Bo52} and again by Marchenko \cite{Ma73} in 1952. 
Both results, the uniqueness theorem and the Gel'fand--Levitan (reconstruction) formalism remained pillars of the inverse spectral theory that withstood any reformulation or improvement for nearly fifty years. Hence it was an incredible
achievement by Barry Simon to have changed the inverse spectral landscape
by offering such a reformulation of inverse spectral theory and in the
very same paper \cite{Si99} to have been able to substantially improve the
Borg--Marchenko uniqueness theorem from a global to a local version.   

We start by highlighting the approach in Simon's paper \cite{Si99} and then switch 
to a more detailed treatment of some aspects of the theory by borrowing from 
\cite{GS00a}.  

Inverse spectral methods have been actively studied in the past 
years both via their relevance in a variety of applications and due to  
their connection with integrable evolution equations such as the KdV equation. In this section, however, we will not deal with the full-line inverse spectral approach relevant to integrable equations but exclusively focus on inverse spectral theory for half-line 
Schr\"odinger operators. In this particular context, a major role is played by 
the Gel'fand--Levitan equations \cite{GL55} (see also, 
\cite[Chs.\ 3, 4]{CCPR97}, \cite{CS89}, \cite{Kr53}, \cite{Kr53a}, \cite{Kr54}, 
\cite[Ch.\ 2]{Le87}, \cite{Ma73}, \cite[Ch.\ 2]{Ma86}, \cite[Ch.\ VIII]{Na68}, \cite{Re02},
\cite{Sy79}, \cite{Th79}). The goal in Barry Simon's paper \cite{Si99} was
to  present a new approach to their basic results. In particular, he 
introduced a  new basic object, the $A$-function (see \eqref{7.1.24}
below), the remarkable equation  \eqref{7.1.28} it satisfies, and
illustrated its fundamental importance with several new results including
improved asymptotic expansions of the Weyl--Titchmarsh $m$-function in
the high-energy regime and the local uniqueness result.

To present some of these new results, we will first describe the 
major players in this game. One is concerned with self-adjoint differential
operators on  either $L^2 ([0,b])$ with $b<\infty$, or $L^2 ([0,\infty))$
associated with differential expressions of the  form
\begin{equation}
-\frac{d^2}{dx^2} + V(x), \quad x\in (0,b). \lb{7.1.1}
\end{equation}
If $b$ is finite, we suppose
$$
\int_0^b dx \, |V(x)| < \infty  
$$
and place a boundary condition
\begin{equation}
u'(b) + hu(b) =0 \lb{7.1.3}
\end{equation}
at $b$, where $h\in\bbR \cup \{\infty\}$ with $h=\infty$ shorthand 
for the Dirichlet boundary condition $u(b)=0$. If $b=\infty$, we suppose 
$$
\int_y^{y+1} dx \, |V(x)| < \infty \, \text{ for all } \, y\geq 0 
$$
and
\begin{equation}
 \sup_{y>0} \int_y^{y+1} dx \, \max(V(x), 0) < 
\infty . \lb{7.1.5}
\end{equation}
Under condition \eqref{7.1.5}, it is known that \eqref{7.1.1} is limit
point at  infinity \cite[App.\ to Sect.\ X.1]{RS75}. In addition, a fixed self-adjoint 
boundary condition at $x=0$ is assumed when talking about the self-adjoint operator associated with \eqref{7.1.1}. 

In either case, for each $z\in\bbC \backslash [\beta, \infty)$ 
with $-\beta$ sufficiently large, there is a unique solution (up 
to an overall constant), $u(z,x)$, of $-u''+Vu =zu$ which satisfies 
\eqref{7.1.3} at $b$ if $b<\infty$ or which is $L^2$ at $\infty$ if 
$b=\infty$. The principal $m$-function $m(z)$ is defined by
$$ 
m(z) =\frac{u'(z,0)}{u(z,0)}.  
$$

If we replace $b$ by $b_1 = b-x_0$ with $x_0\in (0,b)$ and let 
$V(s) = V(x_0+s)$ for $s\in (0,b_1)$, we get a new $m$-function 
we will denote by $m(z,x_0)$. It is given by
$$
m(z,x) = \frac{u'(z,x)}{u(z,x)}.  
$$
$m(z,x)$ satisfies the Riccati-type equation
\begin{equation}
\f{d}{dx}m(z,x) = V(x) - z - m^2 (z,x). \lb{7.1.8}
\end{equation}

Obviously, $m(z,x)$ depends only on $V$ on $(x,b)$ (and on $h$ if 
$b<\infty)$. A basic result of the inverse spectral theory says that the 
converse is true as was shown independently by Borg \cite{Bo52} and 
Marchenko \cite{Ma73} in 1952:

\begin{theorem} \lb{t7.1.1} $m$ uniquely determines $V$. Explicitly, if $V_j$ 
are potentials with corresponding $m$-functions $m_j$, $j=1,2$, and
$m_1 =m_2$, then $V_1=V_2$ a.e.\ $($including $h_1 =  h_2$$)$.
\end{theorem}

In 1999, Simon \cite{Si99} spectacularly improved this to obtain a local 
version of the Borg--Marchenko uniqueness result as follows:

\begin{theorem} \lb{t7.1.2} If $(V_1, b_1, h_1)$, $(V_2, b_2, h_2)$ 
are two potentials and $a <\min(b_1, b_2)$ and if
\begin{equation}
V_1 (x) = V_2 (x) \, \text{ on } \, (0,a), \lb{7.1.9}
\end{equation}
then as $\kappa\to\infty$,
\begin{equation}
m_1 (-\kappa^2) - m_2 (-\kappa^2) = \tilde O(e^{-2\kappa a}). 
\lb{7.1.10}
\end{equation}
Conversely, if \eqref{7.1.10} holds, then \eqref{7.1.9} holds. 
\end{theorem}

In \eqref{7.1.10}, we use the symbol $\tilde O$ defined by 
\begin{align*}
&f=\tilde O(g) \, \text{  
as $x\to x_0$ (where $\lim_{x\to x_0} g(x)=0$)} \\
& \quad \text{if and only if 
$\lim_{x\to x_0} \frac{|f(x)|}{|g(x)|^{1-\varepsilon}} =0$ for 
all $\varepsilon >0$.}
\end{align*}

From a results point of view, this local version of the 
Borg--Marchenko uniqueness theorem was the most significant new 
result in Simon's paper \cite{Si99}, but a major thrust of this paper was
the new set of methods introduced which led to a new approach of the inverse spectral problem. Theorem \ref{t7.1.2} implies that $V$ is
determined by the asymptotics of $m(-\kappa^2)$ as $\kappa \to \infty$. One
can also read off  differences of the boundary condition from these
asymptotics.  Moreover, the following result is proved in \cite{Si99}:

\begin{theorem} \lb{t7.1.3} Let $(V_1, b_1, h_1)$, $(V_2, b_2, h_2)$ 
be two potentials and suppose that
\begin{equation}
b_1 = b_2 \equiv b <\infty, \quad |h_1| + |h_2| < \infty,  
\quad V_1(x) = V_2 (x) \, \text{ on } \, (0,b). \lb{7.1.11}
\end{equation}
Then
\begin{equation}
\lim_{\kappa\to\infty} e^{2b\kappa} 
|m_1 (-\kappa^2) - m_2 (-\kappa^2)| = 4(h_1 - h_2). \lb{7.1.12}
\end{equation}
Conversely, if \eqref{7.1.12} holds for some $b <\infty$ with a 
limit in $(0,\infty)$, then \eqref{7.1.11} holds. 
\end{theorem}

\medskip

To understand Simon's new approach, it is useful to recall briefly 
the two approaches to the inverse problem for Jacobi matrices on 
$\ell^2 (\bbN_0)$ \cite[Ch.\ VII]{Be68}, \cite{GS97a},
\cite{St95}:
$$
A = \begin{pmatrix} b_0 & a_0 & 0 & 0 & \dots \\
a_0 & b_1 & a_1 & 0 & \dots \\
0 & a_1 & b_2 & a_2 & \dots \\
\dots & \dots & \dots & \dots & \dots \end{pmatrix}
$$
with $a_j >0, b_j\in\bbR$. Here the $m$-function is just $(\delta_0, 
(A-z)^{-1}\delta_0) = m(z)$ and, more generally, $m_n(z)=
(\delta_n, (A^{(n)}-z)^{-1}\delta_n)$ with $A^{(n)}$ on 
$\ell^2 (\{n, n+1, \dots\})$ obtained by truncating the 
first $n$ rows and $n$ columns of $A$. Here $\delta_n$ is the 
Kronecker vector, that is, the vector with $1$ in slot $n$ and 
$0$ in other slots. The fundamental theorem in this case is that 
$m(z) \equiv m_0(z)$ determines the $b_n$'s and $a_n$'s.

$m_n(z)$ satisfies an analog of the Riccati equation \eqref{7.1.8}:
\begin{equation}
a^2_n m_{n+1}(z) = b_n - z - \frac{1}{m_n(z)}\, . \lb{7.1.13}
\end{equation}

One solution of the inverse problem is to turn \eqref{7.1.13} around to 
see that
\begin{equation}
m_n (z)^{-1} = -z + b_n - a^2_n m_{n+1} (z) \lb{7.1.14}
\end{equation}
which, first of all, implies that as $z\to\infty$, $m_n(z) = 
-z^{-1} + O(z^{-2})$, so \eqref{7.1.14} implies
\begin{equation}
m_n (z)^{-1} = -z + b_n + a^2_n z^{-1} + O(z^{-2}). \lb{7.1.15}
\end{equation}
Thus, \eqref{7.1.15} for $n=0$ yields $b_0$ and $a^2_0$ and so $m_1(z)$ 
by \eqref{7.1.13}, and then an obvious induction yields successive 
$b_k$, $a^2_k$, and $m_{k+1}(z)$.

A second solution involves orthogonal polynomials. Let $P_n(z)$ 
be the eigensolutions of the formal $(A-z)P_n =0$ with boundary 
conditions $P_{-1}(z)=0$, $P_0(z)=1$. Explicitly,
\begin{equation}
P_{n+1}(z) = a^{-1}_n [(z-b_n)P_n (z)] -a_{n-1} P_{n-1}. 
\lb{7.1.16}
\end{equation}
Let $d\rho$ be the spectral measure for $A$ and vector $\delta_0$ so
that
\begin{equation}
m(z) = \int_{\bbR} \frac{d\rho (\lambda)}{\lambda-z}. \lb{7.1.17}
\end{equation}
Then one can show that
\begin{equation}
\int_{\bbR} d\mu(\lambda) \, P_n(\lambda) P_m(\lambda) =\delta_{n,m}, \quad 
n,m\in\bbN_0. \lb{7.1.18}
\end{equation}

Thus, $P_n(z)$ is a polynomial of degree $n$ with positive  
leading coefficients determined by \eqref{7.1.18}. These orthonormal 
polynomials are determined via Gram--Schmidt from $\rho$ and by 
\eqref{7.1.17} from $m$. Once one has the polynomials $P_n$, one can 
determine the $a$'s and $b$'s from equation \eqref{7.1.16}.

Of course, these approaches via Riccati equation and orthogonal 
polynomials are not completely disjoint. The Riccati solution 
gives the $a_n$'s and $b_n$'s as continued fractions and the 
connection between continued fractions and orthogonal polynomials 
played a fundamental role in Stieltjes' work \cite{St95} on the moment problem 
in 1895.

The Gel'fand--Levitan approach  to the continuum case (cf.\ \cite{GL55},
\cite[Ch.\ 2]{Le87}, \cite{Ma73}, \cite[Ch.\ 2]{Ma86}) is a direct analog
of this orthogonal  polynomial case. One looks at solutions $U(k,x)$ of
\begin{equation}
-U''(k,x) + V(x)U(k,x) = k^2 U(k,x) \lb{7.1.19}
\end{equation}
satisfying $U(k,0)=1$, $U'(k,0)=ik$, and proves that they satisfy a 
representation
\begin{equation}
U(k,x)=e^{ikx} + \int_{-x}^x dy \, K(x,y) e^{iky}, \lb{7.1.20}
\end{equation}
the analog of $P_n(z)=cz^n +$ lower order. One defines $s(k,x) 
=(2ik)^{-1} [U(k,x) - U(-k,x)]$ which satisfies \eqref{7.1.19} with
$s(k,0_+)  =0$, $s'(k,0_+)=1$.

The spectral measure $d\rho$ associated to $m(z)$ by 
$$
d\rho (\lambda) = (2\pi)^{-1} \lim_{\varepsilon\downarrow 0} [ 
\Im (m(\lambda + i\varepsilon))\, d\lambda]
$$ 
satisfies
\begin{equation}
\int_{\bbR}  d\rho(k^2) \, s(k,x) s(k,y) = \delta (x-y), \lb{7.1.21}
\end{equation}
at least formally. \eqref{7.1.20} and \eqref{7.1.21} yield an integral
equation  for $K$ depending only on $d\rho$ and then once one has $K$, 
one can find $U$ and hence $V$ via \eqref{7.1.19} (or via another relation 
between $K$ and $V$).

The principal goal in \cite{Si99} was to present a new approach to the 
continuum case, that is, an analog of the Riccati equation 
approach to the discrete inverse problem. The simple idea for 
this is attractive but has a difficulty to overcome. $m(z,x)$ 
determines $V(x)$, at least if $V$ is continuous by the known 
asymptotics (\cite{DL91}, \cite{Ry02}):
\begin{equation}
m(-\kappa^2, x) = -\kappa - \frac{V(x)}{2\kappa} + o(\kappa^{-1}). 
\lb{7.1.22}
\end{equation}
We can therefore think of \eqref{7.1.8} with $V$ defined by \eqref{7.1.22}
as an  evolution equation for $m$. The idea is that using a suitable 
underlying space and uniqueness theorem for solutions of 
differential equations, \eqref{7.1.8} should uniquely determine $m$ for 
all positive $x$, and hence $V(x)$ by \eqref{7.1.22}.

To understand the difficulty, consider a potential $V(x)$ on 
the whole real line. There are then functions $u_\pm (z,x)$ 
defined for $z\in\bbC\backslash [\beta, \infty)$ which are 
$L^2$ at $\pm\infty$ and two $m$-functions $m_\pm (z,x) = 
u'_\pm (z,x)/u_\pm (z,x)$. Both satisfy \eqref{7.1.8}, yet $m_+ 
(z,0)$ determines and is determined by $V$ on $(0,\infty)$ 
while $m_- (z,0)$ has the same relation to $V$ on $(-\infty, 0)$. 
Put differently, $m_+ (z,0)$ determines $m_+ (z,x)$ for $x>0$ 
but not at all for $x<0$. $m_-$ is the reverse. So uniqueness 
for \eqref{7.1.8} is one-sided and either side is possible! That this 
does not make the scheme hopeless is connected with the fact 
that $m_-$ does not satisfy \eqref{7.1.22}, but rather
\begin{equation}
m_- (-\kappa^2, x)=\kappa + \frac{V(x)}{2\kappa} + 
o(\kappa^{-1}). \lb{7.1.23}
\end{equation}
We will see the one-sidedness of the solvability is intimately 
connected with the sign of the leading $\pm\kappa$ term in 
\eqref{7.1.22}, \eqref{7.1.23}.

The key object in this new approach is a function $A(\alpha)$ 
defined for $\alpha \in (0,b)$ related to $m$ by
\begin{equation}
m(-\kappa^2) = -\kappa - \int_0^a d\alpha\,A(\alpha) e^{-2\alpha\kappa} 
 +\tilde O(e^{-2a\kappa}) \lb{7.1.24}
\end{equation}
as $\kappa\to\infty$. We have written $A(\alpha)$ as a function 
of a single variable but we will allow similar dependence on 
other variables. Since $m(-\kappa^2, x)$ is also an $m$-function, 
\eqref{7.1.24} has an analog with a function $A(\alpha, x)$. 

By uniqueness of inverse Laplace transforms (see  
\cite[Appendix 2, Theorem A.2.2]{Si99}), \eqref{7.1.24} and $m$ near
$-\infty$ uniquely determine $A(\alpha)$.

Not only will \eqref{7.1.24} hold but, in a sense, $A(\alpha)$ is close 
to $V(\alpha)$. Explicitly, one can prove the following result:

\begin{theorem} \lb{t7.1.4} Let $m$ be the $m$-function of the 
potential $V$. Then there is a function $A \in L^1 
([0,b])$ if $b< \infty$ and $A \in L^1 ([0,a])$ for all 
$a<\infty$ if $b=\infty$ so that \eqref{7.1.24} holds for any 
$a\leq b$ with $a<\infty$. $A(\alpha)$ only depends on $V(y)$ 
for $y\in [0,\alpha]$. Moreover, $A(\alpha)=V(\alpha) + 
E(\alpha)$ where $E(\alpha)$ is continuous and satisfies
$$
|E(\alpha)|\leq \biggl( \int_0^\alpha dy \, |V(y)|\biggr)^2
\exp\biggl(\alpha \int_0^\alpha dy \,|V(y)|\biggr). 
$$
\end{theorem}

Restoring the $x$-dependence, we see that $A(\alpha, x) =
V(\alpha + x) + E(\alpha, x)$ where 
$$
\lim_{\alpha\downarrow 0}\, \sup_{0\leq x\leq a} |E(\alpha, x)|=0
$$
for any $a>0$, so
\begin{equation}
\lim_{\alpha\downarrow 0} A(\alpha, x) = V(x), \lb{7.1.26}
\end{equation}
where this holds in general in the $L^1$-sense. If $V$ is 
continuous, \eqref{7.1.26} holds pointwise. In general, \eqref{7.1.26} will
hold  at any point of right Lebesgue continuity of $V$.

Because $E$ is continuous, $A$ determines any discontinuities 
or singularities of $V$. More is true. If $V$ is $C^k$, then
$E$ is $C^{k+2}$ in $\alpha$, and so $A$ determines 
$k^{\text{th}}$ order kinks in $V$. Much more is true and one can also
prove the following result:  

\begin{theorem} \lb{t7.1.5} $V$ on $[0,a]$ is only a function of 
$A$ on $[0,a]$. Explicitly, if $V_1, V_2$ are two potentials,  
let $A_1, A_2$ be their $A$-functions. If $a<b_1$, $a<b_2$, and 
$A_1(\alpha)=A_2(\alpha)$ for $\alpha \in [0,a]$, then 
$V_1 (x)=V_2 (x)$ for $x\in[0,a]$.
\end{theorem}

Theorems \ref{t7.1.4} and \ref{t7.1.5} imply Theorem
\ref{t7.1.2}. 

As noted, the singularities of $V$ come from singularities of 
$A$. A boundary condition is a kind of singularity, so one 
might hope that boundary conditions correspond to very 
singular $A$. In essence, we will see that this is the case -- there are
delta-function and delta-prime singularities  at $\alpha =b$. Explicitly,
one can prove the following result: 

\begin{theorem} \lb{t7.1.6} Let $m$ be the $m$-function for a 
potential $V$ with $b<\infty$. Then for $a<2b$,
\begin{equation}
m(-\kappa^2)=-\kappa -\int_0^a d\alpha \, A(\alpha) e^{-2\alpha\kappa} 
 - A_1 \kappa e^{-2\kappa b} - B_1 e^{-2\kappa b} 
+\tilde O(e^{-2a\kappa}), \lb{7.1.27}
\end{equation}
where the following facts hold: \\
$(a)$ If $h=\infty$, then $A_1=2, \quad B_1 = -2
\int_0^b V(y)\, dy$. \\
$(b)$ If $|h|<\infty$, then $A_1=-2, \quad B_1 =
2[2h+\int_0^b V(y)\, dy]$.
\end{theorem}
This implies Theorem \ref{t7.1.3}.

The reconstruction theorem, Theorem \ref{t7.1.5}, depends on the 
differential equation that $A(\alpha, x)$ satisfies. Remarkably, $V$ 
drops out of the translation of \eqref{7.1.8} to the equation for $A$:
\begin{equation}
\frac{\partial A(\alpha,x)}{\partial x} = 
\frac{\partial A(\alpha,x)}{\partial\alpha} + 
\int_0^\alpha d\beta \, A(\beta,x) A(\alpha-\beta,x). \lb{7.1.28}
\end{equation}

If $V$ is $C^1$, the equation holds in classical sense. For 
general $V$, it holds in a variety of weaker senses. Either way, 
$A(\alpha,0)$ for $\alpha\in [0,a]$ determines $A(\alpha,x)$ for 
all $x,\alpha$ with $\alpha >0$ and $0<x+\alpha <a$. \eqref{7.1.26} 
then determines $V(x)$ for $x\in [0,a)$. That is the essence 
from which uniqueness comes. We will return to this circle of ideas later on when discussing Simon's approach to the inverse spectral problem in detail.

\smallskip
\begin{center}
$\bigstar$  \qquad $\bigstar$ \qquad $\bigstar$
\end{center}
\smallskip

Now we switch to \cite{GS00a} and take a closer look at some of the concepts introduced in \cite{Si99}. In particular, we continue the
study of the $A$-amplitude associated to half-line Schr\"odinger
operators, $-\f{d^2}{dx^2}+ V$ in $L^2 ([0,b))$, $b\leq \infty$. $A$ is
related to the Weyl--Titchmarsh $m$-function via 
$m(-\kappa^2) =-\kappa - \int_0^a d\alpha \, A(\alpha) e^{-2\alpha\kappa}
+O(e^{-(2a -\varepsilon)\kappa})$ for all $\veps >0$. Three main issues will be 
discussed: 

\medskip

\noindent 
$\bullet$ First, we describe how to extend the theory to general $V$
in $L^1 ([0,a])$  for all $a>0$, including $V$'s which are limit circle at
infinity. \\
$\bullet$ Second, the following relation between the
$A$-amplitude and the spectral  measure $\rho$: 
$$
A(\alpha) = -2\int_{-\infty}^\infty d\rho(\lambda) \, \lambda^{-\frac12} 
\sin (2\alpha \sqrt{\lambda}), 
$$ 
will be discussed. (Since the integral is 
divergent, this formula has to be properly interpreted.) \\
$\bullet$ Third, a Laplace transform representation for $m$ without error term 
in the case $b<\infty$ will be presented. 
\medskip

We consider Schr\"odinger operators
\begin{equation} \label{1.1}
-\frac{d^2}{dx^2} + V
\end{equation}
in $L^2 ([0,b))$ for $0<b<\infty$ or $b=\infty$ and real-valued locally  
integrable $V$. There are essentially four distinct cases. 

\medskip
\noindent{\bf Case 1.} $b<\infty$. We suppose $V\in L^1 ([0,b])$. We then 
pick $h\in \bbR \cup \{\infty\}$ and add the boundary condition at $b$
\begin{equation}
\label{1.2}
u'(b_-) + hu(b_-)=0, 
\end{equation}
where $h=\infty$ is shorthand for the Dirichlet boundary condition 
$u(b_-)=0$. 

For Cases 2--4, $b=\infty$ and
\begin{equation}
\label{1.3}
\int_0^a dx \, |V(x)| < \infty \quad\text{for all }a<\infty. 
\end{equation}

\medskip
\noindent{\bf Case 2.} $V$ is ``essentially'' bounded from below in the 
sense that 
\begin{equation} \lb{1.3a}
\sup_{a>0}  \left(\int_a^{a+1} dx\, \max (-V(x),0)\right) < \infty.
\end{equation} 
Examples include $V(x) = c(x+1)^\beta$ for $c>0$ and all $\beta\in\bbR$ 
or $V(x) = -c (x+1)^\beta$ for all $c>0$ and $\beta \leq 0$.

\medskip
\noindent{\bf Case 3.} \eqref{1.3a} fails but \eqref{1.1} is limit point at 
$\infty$ (see, e.g., \cite[Ch.\ 9]{CL55}, \cite[Sect.\ X.1]{RS75} for a
discussion of  limit point/limit circle), that is, for each $z\in\bbC_+ =
\{z\in\bbC \,|\, \Im (z)>0\}$,
\begin{equation}
\label{1.4}
-u'' + Vu = zu 
\end{equation}
has a unique solution, up to a multiplicative constant, which is $L^2$ at 
$\infty$. An example is $V(x) = -c(x+1)^\beta$ for $c>0$ and 
$0< \beta\leq 2$.

\medskip
\noindent{\bf Case 4.} \eqref{1.1} is limit circle at infinity, that is, 
every solution of \eqref{1.4} is $L^2 ([0,\infty))$ at infinity if 
$z\in\bbC_+$. We then pick a boundary condition by choosing a nonzero
solution $u_0$ of \eqref{1.4} for $z=i$. Other functions $u$ satisfying
the associated boundary  condition at infinity then are supposed to satisfy
\begin{equation}
\label{1.5}
\lim_{x\to\infty} W(u_0,u)(x)=\lim_{x\to\infty} [u_0 (x) u'(x) - u'_0 (x) u(x)] = 0. 
\end{equation}
Examples include $V(x) = -c(x+1)^\beta$ for $c>0$ and $\beta >2$.

\medskip
The Weyl--Titchmarsh $m$-function, $m(z)$, is defined for $z\in\bbC_+$ as 
follows. Fix $z\in\bbC_+$. Let $u(x,z)$ be a nonzero solution of 
\eqref{1.4} which satisfies the boundary condition at $b$. In Case 1, that
means $u$ satisfies \eqref{1.2}; in Case 4, it satisfies \eqref{1.5}; and
in Cases 2--3, it satisfies  
$\int_R^\infty |u (z,x)|^2\, dx<\infty$ for some (and hence for all) 
$R\geq 0$. Then,
\begin{equation} 
\label{1.6}
m(z) = \f{u'(z,0_+)}{u(z,0_+)} 
\end{equation}
and, more generally,
\begin{equation}\label{1.7}
m(z,x) = \f{u'(z,x)}{u(z,x)}\,. 
\end{equation}
$m(z,x)$ satisfies the Riccati equation 
(with $m'=\partial m/\partial x$), 
\begin{equation}\label{1.8a}
m'(z,x) = V(x) - z - m(z,x)^2. 
\end{equation}

$m$ is an analytic function of $z$ for $z\in\bbC_+$, and the following properties hold:
\medskip

\noindent{\bf Case 1.} $m$ is meromorphic in $\bbC$ with a discrete set 
$\lambda_1 < \lambda_2 < \cdots$ of poles on $\bbR$ (and none on $(-\infty, 
\lambda_1)$).

\smallskip
\noindent{\bf Case 2.} For some $\beta\in\bbR$, $m$ has an analytic continuation 
to $\bbC\backslash[\beta, \infty)$ with $m$ real on $(-\infty, \beta)$.

\smallskip
\noindent{\bf Case 3.} In general, $m$ cannot be continued beyond $\bbC_+$ 
(there exist $V$'s where $m$ has a dense set of polar singularities on
$\bbR$).

\smallskip
\noindent{\bf Case 4.} $m$ is meromorphic in $\bbC$ with a discrete set of 
poles (and zeros) on $\bbR$ with limit points at both $+\infty$ and $-\infty$.

\smallskip

Moreover,
\[
\text{if } z\in\bbC_+ \text{ then } m(z,x)\in\bbC_+ ,
\]
so $m$ admits the Herglotz representation, 
\begin{equation}
\label{1.8b}
m(z) = \Re(m(i)) + \int_\bbR d\rho(\lambda)  
\left[ \f{1}{\lambda-z} -\f{\lambda}{1+\lambda^2} \right], \quad z\in\bbC\backslash\bbR,
\end{equation}
where $\rho$ is a positive measure called the spectral measure, which satisfies 
\begin{gather}
\int_\bbR \f{d\rho(\lambda)}{1+|\lambda|^2} < \infty , \label{1.9} \\
d\rho(\lambda) = \wlim_{\veps\downarrow 0} 
\frac{1}{\pi} \Im (m(\lambda + i\veps))\, d\lambda , \label{1.10}
\end{gather}
where $\wlim$ is meant in distributional sense.

\smallskip
All these properties of $m$ are well known (see, e.g.\ 
\cite[Ch.\ 2]{LS75}).

\medskip

In \eqref{1.8b}, the constant $\Re (m(i))$ is determined 
by the result of Everitt \cite{Ev72} that for each $\veps>0$, 
\begin{equation}
\label{1.11}
m(-\kappa^2) = -\kappa +o(1) \quad\text{as } |\kappa| \to \infty 
\text{ with } -\f{\pi}{2} + \veps < \arg (\kappa) <  - \veps < 0.
\end{equation}

Atkinson \cite{At81} improved \eqref{1.11} to read, 
\begin{equation}
\label{1.12}
m(-\kappa^2) = -\kappa + \int_0^{a_0} d\alpha \, V(\alpha)
e^{-2\alpha\kappa}  + o(\kappa^{-1}) 
\end{equation}
again as $|\kappa|\to\infty$ with $- \frac{\pi}{2} +\veps < \arg(\kappa) < 
- \veps <0$ (actually, he allows $\arg(\kappa) \to 0$ as 
$|\kappa| \to \infty$ 
as long as $\Re (\kappa) >0$ and $\Im (\kappa) > -\exp(-D|\kappa|)$ for 
suitable $D$). In \eqref{1.12}, $a_0$ is any fixed $a_0 >0$.

One of the main results in \cite{GS00a} was to go way beyond the two 
leading orders in \eqref{1.12}. 

\begin{theorem} \label{t1} There exists a function $A(\alpha)$ for $\alpha
\in [0,b)$ so that $A\in L^1 ([0,a])$ for all $a<b$ and
\begin{equation}
\label{1.13}
m(-\kappa^2) = -\kappa - \int_0^a d\alpha \,A(\alpha) e^{-2\alpha\kappa}
+ \tilde O(e^{-2\alpha\kappa})
\end{equation}
as $|\kappa|\to\infty$ with $-\f{\pi}{2} + \veps<\arg(\kappa)< - \veps <0$. 
Here we say $f=\tilde O(g)$ if $g\to 0$ and for all
$\veps >0$, $\big(\frac{f}{g}\big)|g|^\veps\to 0$ as 
$|\kappa| \to \infty$. Moreover, $A-q$ is continuous and
\begin{equation}
\label{1.14}
|(A-q)(\alpha)|\leq \left[\int_0^\alpha dx \,|V(x) \right]^2 \exp 
\left(\alpha \int_0^\alpha dx \, |V(x)|\right).
\end{equation}
\end{theorem}

This result was proved in Cases 1 and 2 in \cite{Si99}. The proof of this
result if one only assumes \eqref{1.3}  (i.e., in Cases~3 and 4) has been
provided in \cite{GS00a}. 

Actually, in \cite{Si99}, \eqref{1.13} was proved in Cases~1 and 2 for 
$\kappa$ real with $|\kappa| \to\infty$. The proof in \cite{GS00a}, assuming 
only condition \eqref{1.3},  
includes Case~2 in the general $\kappa$-region $\arg(\kappa)\in (-\f{\pi}{2}
+ \veps, -\veps)$ and, as can be shown, the proof also holds in this region 
for Case~1.

\begin{remark} At first sight, it may appear that Theorem \ref{t1}, 
as stated, does not imply the $\kappa$ real result of \cite{Si99}, but 
if the spectral measure $\rho$ of \eqref{1.8b} has $\supp(\rho)\subseteq [a,\infty)$ 
for some $a\in \bbR$, \eqref{1.13} extends to all $\kappa$ in $|\arg(\kappa)| 
< \f{\pi}{2} - \veps$, $|\kappa| \geq a +1$. To see this, one notes by 
\eqref{1.8b} that $m'(z)$ is bounded away from $[a,\infty)$ so one has 
the a priori bound $|m(z)| \leq C|z|$ in the region $\Re (z) < a -1$. 
This bound and a Phragm\'en--Lindel\"of argument let one extend \eqref{1.13} 
to the real $\kappa$ axis.
\end{remark}

The following is a basic result from \cite{Si99}:

\begin{theorem} \lb{t2} \mbox{\rm (Theorem~2.1 of \cite{Si99})} Let 
$V\in L^1 ([0,\infty))$. Then there exists a function $A$ on
$(0,\infty)$ so  that $A-V$ is continuous and satisfies \eqref{1.14} such
that for $\Re (\kappa)  > \|V\|_1/2$,
\begin{equation} \lb{1.15}
m(-\kappa^2) = -\kappa - \int_0^\infty d\alpha\, A(\alpha)
e^{-2\alpha\kappa}.
\end{equation}
\end{theorem}

\begin{remark} In \cite{Si99}, this is only stated for $\kappa$ real with 
$\kappa > \|V\|_1/2$, but \eqref{1.14} implies that 
$|A(\alpha) -V(\alpha)| 
\leq \|V\|_1^2 \exp (\alpha \|V\|_1)$ so the right-hand side of 
\eqref{1.15} converges to an analytic function in $\Re (\kappa) > 
\|V\|_1/2$. Since $m(z)$ is analytic in $\bbC\backslash [\alpha, \infty)$
for suitable $\alpha$,  we have equality in $\{ \kappa\in\bbC \,|\, \Re
(\kappa) > \|V\|_1/2\}$  by analyticity.
\end{remark}

Theorem \ref{t1} in all cases follows from Theorem \ref{t2} and the following 
result which was proved in \cite{GS00a}:

\begin{theorem} \lb{t3} Let $V_1, V_2$ be potentials defined on $(0, b_j)$ 
with $b_j >a$ for $j=1,2$. Suppose that $V_1 =V_2$ on $[0,a]$. Then in the 
region $\arg (\kappa) \in (-\f{\pi}{2}+\veps, -\veps)$, $|\kappa| \geq 
K_0$, we have that
\begin{equation} \lb{1.16}
|m_1 (-\kappa^2) - m_2 (-\kappa^2)| \leq C_{\veps,\delta} \exp 
(-2a \Re (\kappa)), 
\end{equation}
where $C_{\veps,\delta}$ depends only on $\veps$, $\delta$, and 
$\sup_{0\leq x\leq a} \big(\int_x^{x+\delta} dy\, |V_j (y)|\big)$, where $\delta>0$
is  any number so that $a+\delta \leq b_j$, $j=1,2$.
\end{theorem}

\begin{remark} $(i)$ An important consequence of Theorem \ref{t3} is that 
if  $V_1(x) = V_2 (x)$ for $x\in [0,a]$, then $A_1(\alpha) = A_2(\alpha)$ 
for $\alpha\in [0,a]$. Thus, $A(\alpha)$ is only a function of $V$ on
$[0,\alpha]$.  (We emphasize that, conversely, one can show that 
also $V(x)$ is only a function of $A$ on $[0,x]$.) \\
$(ii)$ This implies Theorem \ref{t1} by taking $V_1 = V$ and $V_2 =
V\chi_{[0,a]}$  and using Theorem \ref{t2} on $V_2$.  \\
$(iii)$ The actual proof implies \eqref{1.16} on a larger region than $\arg
(\kappa) \in   (-\f{\pi}{2}+\veps, -\veps)$. Basically, one needs $\Im
(\kappa) \geq  -C_1 \exp(-C_2|\kappa|)$ as $\Re (\kappa) \to\infty$.
\end{remark}

The basic connection between the spectral 
measure $d\rho$ and the $A$-amplitude established in \cite{GS00a} says 
\begin{equation} \lb{1.18}
A(\alpha) = -2 \int_{-\infty}^\infty d\rho(\lambda) \,\lambda^{-\f12}
\sin(2\alpha \sqrt{\lambda}\,).
\end{equation}
However, the integral in \eqref{1.18} is not convergent. Indeed, the asymptotics  
\eqref{1.11} imply that $\int_0^R d\rho(\lambda)\sim\f{2}{3\pi} R^{\f32}$ 
so \eqref{1.18} is never absolutely convergent. Thus, \eqref{1.18} has to be suitably
interpreted.

We will indicate how to demonstrate \eqref{1.18} as a distributional
relation,  smeared in $\alpha$ on both sides by a function $f\in
C_0^\infty ((0,\infty))$.  This holds for all $V$'s in Cases~1--4.
Finally, we will discuss an  Abelianized version of \eqref{1.18}, namely,
\begin{equation} \lb{1.19}
A(\alpha) = -2\lim_{\veps\downarrow 0} \int_{-\infty}^\infty 
d\rho(\lambda) \, e^{-\veps\lambda} 
\lambda^{-\f12} \sin(2\alpha\sqrt{\lambda}\,) 
\end{equation}
at any point, $\alpha$, of Lebesgue continuity for $V$. \eqref{1.19} is 
proved only for a restricted class of $V$'s including Case~1, 2 and those
$V$'s  satisfying
$$
V(x) \geq -Cx^2, \quad x\geq R
$$
for some $R>0$, $C>0$, which are always in the limit point case at infinity. 
Subsequently, we will use \eqref{1.19} as the point of departure for
relating $A(\alpha)$ to scattering data.

In order to prove \eqref{1.18} for finite $b$, one needs to analyze the
finite $b$  case extending \eqref{1.13} to all $a$ including $a=\infty$
(by allowing $A$ to  have $\delta$ and $\delta'$ singularities at
multiples of $b$). This was originally done in 
\cite{Si99} for $\kappa$ real and positive and $a<\infty$. We now need
results  in the entire region $\Re (\kappa)\geq K_0$. Explicitly, the following 
was proved in \cite{GS00a}: 

\begin{theorem} \lb{t6} In Case~1, there are $A_n, B_n$ for $n=1,2,\dots$, and a 
function $A(\alpha)$ on $(0,\infty)$ with 
\begin{align*}
& |A_n|\leq C, \quad |B_n| \leq Cn, \\
& \int_0^a  d\alpha \, |A(\alpha)| \leq C\exp (K_0 |a|) \, 
\text{ so that for $\Re (\kappa) > \f12 K_0$:} \\
& m(-\kappa^2) = -\kappa - \sum_{n=1}^\infty A_n \kappa e^{-2\kappa bn} 
-\sum_{n=1}^\infty B_n e^{-2\kappa bn} - \int_0^\infty 
 d\alpha \, A(\alpha) 
e^{-2\alpha\kappa}.
\end{align*}
\end{theorem}

\eqref{1.18} can be used to obtain a priori bounds on $\int_{-R}^0
d\rho(\lambda)$ as $R\to\infty$.

\medskip

Now we turn to more details and start by illustrating how to use the
Riccati equation and a priori control on $m_j$ to obtain exponentially
small estimates on $m_1 - m_2$.

\begin{lemma} \lb{p2.1} Let $m_1, m_2$ be two absolutely
continuous  functions on $[a,b]$ so that for some $Q\in L^1 ([a,b])$, 
\begin{equation} \lb{2.1}
m'_j (x) = Q(x) - m_j (x)^2, \quad j=1,2, \ x\in (a,b).
\end{equation}
Then
\[
[m_1(a) - m_2(a)] = [m_1(b) - m_2(b)] \exp 
\left(\int_a^b dy \, [m_1(y) + m_2(y)] \right).
\]
\end{lemma}

As an immediate corollary, one obtains the following result (which implies 
Theorem~\ref{t3}):

\begin{theorem} \lb{t2.2} Let $m_j (-\kappa^2,x)$ be functions defined for 
$x\in [a,b]$ and $\kappa\in K$ some region of $\bbC$. Suppose that for each 
$\kappa$ in $K$\!, $m_j$ is absolutely continuous in $x$ and satisfies
{\rm (}note that $V$ is the same for $m_1$ and $m_2${\rm)},
\[
m'_j(-\kappa^2,x) = V(x) + \kappa^2 - m_j (-\kappa^2,x)^2, \quad j=1,2. 
\]
Suppose $C$ is such that for each $x\in [a,b]$ and $\kappa\in K$\!,
\begin{equation} \lb{2.2}
|m_j (-\kappa^2,x) + \kappa| \leq C, \quad j=1,2,
\end{equation}
then
\begin{equation} \lb{2.3}
|m_1 (-\kappa^2,a) - m_2 (-\kappa^2,a)| \leq 2C \exp [-2(b-a) 
[\Re (\kappa) -C]].
\end{equation}
\end{theorem}

Theorem~\ref{t2.2} places importance on a priori bounds of the form 
\eqref{2.2}. Fortunately, by modifying ideas of Atkinson \cite{At81}, we
can  obtain estimates of this form as long as $\Im (\kappa)$ is bounded 
away from zero.

Atkinson's method allows one to estimate $|m(-\kappa^2)+\kappa|$ in two
steps.  We will fix some $a<b$ finite and define $m_0(-\kappa^2)$ by
solving
\begin{subequations} \lb{3.3}
\begin{align}
m'_0 (-\kappa^2, x) &= V(x) + \kappa^2 - m_0 (-\kappa^2, x)^2, \lb{3.3a} \\
m_0(-\kappa^2, a) &= - \kappa \lb{3.3b} \\
\intertext{and then setting}
m_0 (-\kappa^2) & := m_0 (-\kappa^2, 0_+). \lb{3.3c}  
\end{align}
\end{subequations}

One then proves the following result.

\begin{lemma} \lb{p3.1} There is a $C>0$ depending only on $V$ and a 
universal constant $E>0$ so that if $\Re (\kappa) \geq C$ and $\Im
(\kappa) 
\neq 0$, then 
\begin{equation} \lb{3.4}
|m(-\kappa^2) -m_0(-\kappa^2)| \leq E\, \f{|\kappa|^2}{|\Im (\kappa)|}\, 
e^{-2a\Re (\kappa)}.
\end{equation}
In fact, one can take
\[
C =\max \left(a^{-1} \ln (6), \, 4 \int_0^a dx\, |V(x)| \right), \quad
E = \f{3 \cdot 2 \cdot 12^2}{5}\, .
\]
\end{lemma}

\begin{lemma} \lb{p3.2} There exist constants $D_1$ and $D_2$  
{\rm(}depending only on $a$ and $V${\rm)}, so that for $\Re (\kappa)
>D_1$, 
\[
|m_0 (-\kappa^2) + \kappa| \leq D_2.
\]
Indeed, one can take
\[
D_1 = D_2 = 2\int_0^a dx \, |V(x)|.
\]
\end{lemma}

These Lemmas together with Theorem \ref{t2} yield the following
explicit  form of Theorem \ref{t3}. 

\begin{theorem} \lb{t3.3} Let $V_1, V_2$ be defined on $(0, b_j)$ with
$b_j  >a$ for $j=1,2$. Suppose that $V_1 = V_2$ on $[0,a]$. Pick $\delta$
so that $a + \delta \leq \min(b_1, b_2)$ and let $\eta=\sup_{0\leq x \leq
a; j=1,2}  (\int_x^{x+\delta} dy\, |V_j (y)|)$. Then if $\Re (\kappa) \geq
\max(4\eta, 
\delta^{-1}\ln (6))$ and $\Im (\kappa)\neq 0$, one obtains
\[ 
|m_1 (-\kappa^2) - m_2(-\kappa^2)|\leq 2 g(\kappa) \exp (-2a[\Re (\kappa) 
- g(\kappa)]),
\]
where
\[
g(\kappa) = 2\eta + \f{864}{5}\, \f{|\kappa|^2}{|\Im (\kappa)|}\, 
e^{-2\delta\Re (\kappa)}\, .
\]
\end{theorem}

\begin{remark} $(i)$ To obtain Theorem~\ref{t3}, we need only note that 
in the region $\arg (\kappa) \in (-\f{\pi}{2}+\veps, -\veps)$, 
$|\kappa|\geq K_0$, $g$ is bounded. \\
$(ii)$ We need not require that $\arg(\kappa)< -\veps$ to obtain $g$
bounded.  It suffices, for example, that $\Re (\kappa) \geq |\Im
(\kappa)| \geq  e^{-\alpha\Re (\kappa)}$ for some $\alpha < 2\delta$. \\
$(iii)$ For $g$ to be bounded, we need not require that $\arg(\kappa) > -
\f{\pi}{2} +\veps$. It suffices that $|\Im (\kappa)| \geq \Re (\kappa)
\geq \alpha\ln  [|\Im (\kappa)|]$ for some $\alpha >(2\delta)^{-1}$.
Unfortunately, this does  not include the region $\Im (-\kappa^2)=c$,
$\Re (-\kappa^2)\to\infty$, where 
$\Re (\kappa)$ goes to zero as $|\kappa|^{-1}$. However, as
$\Re (-\kappa^2)
\to\infty$, we only need that $|\Im(-\kappa^2)|\geq 2\alpha |\kappa| \ln 
(|\kappa|)$.
\end{remark}

Next, we turn to finite $b$ representations with no
errors: Theorem \ref{t2} implies that if $b=\infty$ and $V\in L^1
([0,\infty))$, then 
\eqref{1.15} holds, a Laplace transform representation for $m$ without 
errors. It is, of course, of direct interest that such a formula holds,
but we are  especially interested in a particular consequence of it -- 
namely, that it implies  that the formula \eqref{1.13} with error holds in
the region $\Re (\kappa) > K_0$  with error uniformly bounded in
$\Im(\kappa)$; that is, one proves the following result:

\begin{theorem}\lb{t4.1} If $V\in L^1 ([0, \infty))$ and $\Re (\kappa) > 
 \|V\|_1/2$, then for all $a>0$:
\begin{equation} \lb{4.1}
\left| m(-\kappa^2) + \kappa + \int_0^a d\alpha \, 
A(\alpha) e^{-2\alpha\kappa} \right| \leq \left[ \|V\|_1 + \f{\|V\|^2_1
e^{a\|V\|_1}} {2\Re (\kappa) - \|V\|_1}\right] e^{-2a\Re (\kappa)}.
\end{equation}
\end{theorem}

The principal goal is to prove an analog of this result in the case 
$b<\infty$. To do so, we will need to first prove an analog of 
\eqref{1.15} in case $b<\infty$ -- something of interest in its own
right. The idea will be to  mimic the proof of Theorem~2 from \cite{Si99}
but use the finite $b$, $V^{(0)} (x)=0$, $x\geq0$ Green's function where
\cite{Si99} used the infinite $b$ Green's  function. The basic idea is
simple, but the arithmetic is a bit involved.

We will start with the $h=\infty$ case. Three functions for 
$V^{(0)}(x)=0$, $x \geq 0$ are significant. First, the kernel of the
resolvent $(-\f{d^2}{dx^2} + \kappa^2)^{-1}$ with $u(0_+)=u(b_-)=0$
boundary conditions. By an elementary  calculation (see, e.g.,
\cite[Sect.\ 5]{Si99}), it has the form
\begin{equation} \lb{4.2}
G^{(0)}_{h=\infty} (-\kappa^2,x,y) = \f{\sinh(\kappa x_<)}{\kappa} 
\left[ \f{e^{-\kappa x_>} - e^{-\kappa (2b-x_>)}}{1-e^{-2\kappa b}}\right],
\end{equation}
with $x_<=\min(x,y)$, $x_> = \max(x,y)$.

The second function is
\begin{equation} \lb{4.3}
\psi^{(0)}_{h=\infty} (-\kappa^2,x) = \lim_{y\downarrow 0} 
\f{\partial G^{(0)}_{h=\infty}}{\partial y}\, (-\kappa^2,x,y) = 
\f{e^{-\kappa x} - e^{-\kappa (2b-x)}}{1-e^{-2\kappa b}}
\end{equation}
and finally (notice that $\psi^{(0)}_{h=\infty}(-\kappa^2,0_+)=1$ and 
$\psi^{(0)}_{h=\infty}$ satisfies the equations 
$-\psi'' =-\kappa^2 \psi$ and $\psi(-\kappa^2,b_-)=0$):
\begin{equation} \lb{4.4}
m^{(0)}_{h=\infty} (-\kappa^2) 
= \psi^{(0)\prime}_{h=\infty} (-\kappa^2,0_+) 
= -\f{\kappa + \kappa e^{-2\kappa b}}{1-e^{-2\kappa b}}\, .
\end{equation}
In \eqref{4.4}, prime means $d/dx$.

Now fix $V\in C^\infty_0 ((0,b))$. The pair of formulas
\begin{equation*}
\begin{split}
&\left(-\f{d^2}{dx^2} + V + \kappa^2 \right)^{-1} \\
&\qquad \qquad = \sum_{n=0}^\infty 
(-1)^n \left( -\f{d^2}{dx^2} + \kappa^2\right)^{-1} \left[ V
\left( -\f{d^2}{dx^2} + \kappa^2\right)^{-1}\right]^n 
\end{split}
\end{equation*}
and
\[
m(-\kappa^2) = \lim_{x<y;\,  y\downarrow 0} 
\f{\partial^2 G (-\kappa^2,x,y)}{\partial x\partial y}
\]
yields the following expansion for the $m$-function of $-\f{d^2}{dx^2} +
V$  with $u(b_-)=0$ boundary conditions. 

\begin{lemma} \lb{p4.2} Let $V\in C^\infty_0 ((0, b))$, $b<\infty$. Then 
\begin{equation} \lb{4.5}
m(-\kappa^2) = \sum_{n=0}^\infty M_n (-\kappa^2; V),
\end{equation}
where
\begin{align}
M_0 (-\kappa^2; V) &= m^{(0)}_{h=\infty} (-\kappa^2), \lb{4.6} \\
M_1 (-\kappa^2; V) &= -\int_0^b V(x) \psi^{(0)}_{h=\infty}
(-\kappa^2,x)^2 \, dx , \lb{4.7}
\end{align}
and for $n\geq 2$, 
\begin{align} 
\begin{split}
M_n (-\kappa^2; V) &= (-1)^n \int_0^b dx_1 \dots \int_0^b 
dx_n \, V(x_1) \dots V(x_n)  \\
& \quad \times \psi^{(0)}_{h=\infty}(-\kappa^2,x_1) \psi^{(0)}_{h=\infty}
(-\kappa^2,x_n) \prod_{j=1}^{n-1} 
G^{(0)}_{h=\infty} (-\kappa^2, x_j, x_{j+1}) .  \lb{4.8}
\end{split}
\end{align}
\end{lemma}

The precise region of convergence is unimportant since one can  
expand regions by analytic continuation. For now, we note it certainly converges 
in the region $\kappa$ real with $\kappa^2 > \|V\|_\infty$.

Writing each term in \eqref{4.5} as a Laplace transform then yields the
following result:

\begin{theorem} \lb{t4.3} {\rm{(Theorem~\ref{t6} for $h =\infty$)}} 
Let $b<\infty$, $h =\infty$, and $V\in L^1 ([0,b])$. Then for $\Re
(\kappa) >  \|V\|_1/2$, we have that
\begin{equation} \lb{4.22}
m(-\kappa^2) = -\kappa - \sum_{j=1}^\infty A_j \kappa e^{-2\kappa bj} 
-\sum_{j=1}^\infty B_j e^{-2\kappa bj} - \int_0^\infty d\alpha\, A(\alpha) 
e^{-2\alpha\kappa} ,
\end{equation}
where
\begin{align*}
& A_j =2, \quad B_j = -2j\int_0^b dx \, V(x), \quad j\in\bbN, \\
& |A(\alpha) - A_1 (\alpha)| \leq \f{(2\alpha 
+ b)(2\alpha +2b)}
{2b^2} \|V\|^2_1 \exp (\alpha \|V\|_1)\, \text{ with $A_1$ given by} \\
& A_1(\alpha) = \begin{cases}
V(\alpha), & 0 \leq\alpha < b, \\
(n+1) V(\alpha - nb) + nq ((n+1)b-\alpha), & nb\leq\alpha < (n+1)b, \;
 n\in\bbN .\end{cases}
\end{align*}
In particular, for all $a\in(0,b)$,
\[
\int_0^a d\alpha \,|A(\alpha)| \leq C (b, \|V\|_1) (1+a^2) \exp(a\|V\|_1).
\]
\end{theorem}

This implies the following estimate:

\begin{corollary} \lb{c4.4} If $V\in L^1 ([0,\infty))$ and $\Re (\kappa) 
\geq\f12 \|V\|_1 + \veps$, then for all $a\in(0,b)$, $b<\infty$, we have 
that  
\[
\left| m(-\kappa^2) + \kappa + \int_0^a d\alpha \, 
A(\alpha) e^{-2\alpha\kappa} \right| \leq C(a,\veps) e^{-2a\Re (\kappa)},
\]
where $C(a,\veps)$ depends only on $a$ and $\veps$ {\rm(}and
$\|V\|_1$\rm{)}  but not on $\Im(\kappa)$.
\end{corollary}

Next, we turn to the case $h\in\bbR$. Then \eqref{4.2}--\eqref{4.4} become
\begin{align}
G^{(0)}_{h}(-\kappa^2,x,y) &= \f{\sinh (\kappa x_<)}{\kappa}\, 
\psi^{(0)}_{h} (-\kappa^2,x_>) , \lb{4.26} \\
\psi^{(0)}_{h}(-\kappa^2,x) &= \left[ \f{e^{-\kappa x} + \zeta(h,\kappa) 
e^{-\kappa (2b-x)}}{1+\zeta (h,\kappa) e^{-2b\kappa}}\right], \lb{4.27} \\
m^{(0)}_{h}(-\kappa^2) &= -\kappa + 2\kappa\, \f{\zeta(h,\kappa) e^{-2\kappa b}}
{1+\zeta(h,\kappa) e^{-2\kappa b}}\, , \lb{4.28}
\end{align}
where
\begin{equation} \lb{4.29}
\zeta(h,\kappa) = \f{\kappa-h}{\kappa+h}\,. 
\end{equation}

This then leads to the following result:

\begin{theorem} \lb{t4.7} {\rm{(Theorem~\ref{t6} for general
$h\in\bbR$)}}  Let $b<\infty$, $|h|<\infty$, and $V\in L^1 ([0,b])$. Then
for
$\Re (\kappa)  > \f12 D_1 [\|V\|_1 + |h|+b^{-1} +1]$ for a suitable
universal constant $D_1$, 
\eqref{4.22} holds, where
\begin{align}
& A_j = 2(-1)^j, \quad B_j = 2(-1)^{j+1} 
j\bigg[2h + \int_0^b dx \, V(x) \bigg], \\
& |A(\alpha) - V(\alpha)| \leq \|V\|^2_1\exp(\alpha \|V\|_1) \, 
\text{ if $|\alpha| <b$, and for any $a>0$}, \\ 
& \quad \int_0^a d\alpha \,|A(\alpha)| \leq D_2 (b, \|V\|_1,h) \exp(D_1 a
(\|V\|_1  + |h|+b^{-1}+1)).
\end{align}
\end{theorem}

Hence, one obtains the following estimate:

\begin{corollary} \lb{c4.8} Fix $b<\infty$, $V\in L^1 ([0,b])$, and $|h| 
<\infty$. Fix $a<b$. Then there exist positive constants $C$ and $K_0$ so 
that for all complex $\kappa$ with $\Re (\kappa) > K_0$, 
\[
\left| m(-\kappa^2) + \kappa + \int_0^a d\alpha \, 
A(\alpha) e^{-2\alpha\kappa} \right| \leq Ce^{-2a\kappa}.
\]
\end{corollary}

Next we return to the relation between $A$ and $\rho$ and
discuss a first distributional form of this relation: Our primary goal in
the following is to discuss a formula which formally 
says that
\begin{equation} \lb{5.1}
A(\alpha) = -2\int_{-\infty}^\infty d\rho(\lambda) \, 
\lambda^{-\f12} \sin(2\alpha 
\sqrt{\lambda}\,),
\end{equation}
where for $\lambda \leq 0$, we define
\[
\lambda^{-\f12}\sin(2\alpha\sqrt{\lambda}\,) = \begin{cases} 
2\alpha &\text{if } \lambda =0,  \\
(-\lambda)^{-\f12} \sinh(2\alpha \sqrt{-\lambda}\, ) &\text{if } \lambda <0.
\end{cases}
\]
In a certain sense which will become clear, the left-hand side of 
\eqref{5.1} should be $A(\alpha) - A(-\alpha) + \delta'(\alpha)$.

To understand \eqref{5.1} at a formal level, note the basic formulas,
\begin{align}
m(-\kappa^2) &= -\kappa - \int_0^\infty d\alpha \, 
A(\alpha) e^{-2\alpha\kappa},
\lb{5.2} \\ m(-\kappa^2) &= \Re  (m(i)) + \int_{-\infty}^\infty 
d\rho(\lambda) \left[
\f{1} {\lambda+\kappa^2} - \f{\lambda}{1+\lambda^2}\right], \lb{5.3}
\end{align}
and
\begin{equation} \lb{5.4}
(\lambda + \kappa^2)^{-1} = 2\int_0^\infty d\alpha \,\lambda^{-\f12} 
\sin(2\alpha \sqrt{\lambda}\,) e^{-2\alpha\kappa}, 
\end{equation}
which is an elementary integral if $\kappa >0$ and $\lambda >0$. Plug 
\eqref{5.4} into \eqref{5.3}, formally interchange the order of integrations, 
and \eqref{5.2} should only hold if \eqref{5.1} does. However, a closer 
examination of this procedure reveals that the interchange of the order of 
integrations is not justified and indeed \eqref{5.1} is not true as a 
simple integral since, 
$\int_0^R d\rho (\lambda) \simlim\limits_{R\to\infty}\f{2}{3\pi} R^{\f32}$, 
which implies that \eqref{5.1}  is not absolutely convergent. We will even 
see later that the integral sometimes fails to be conditionally convergent.

Our primary method for understanding \eqref{5.1} is as a distributional 
statement, that is, it will hold when smeared in $\alpha$ for $\alpha$ in 
$(0,b)$. We discuss this next if $V\in L^1 ([0,\infty))$ or if 
$b<\infty$. Later it will be extended to all $V$ (i.e., 
all Cases 1--4) by a limiting argument. Subsequently, we will study 
\eqref{5.1} as a pointwise statement, where  the integral is defined as an
Abelian limit. 

Suppose $b<\infty$ or $b=\infty$ and $V\in L^1 ([0,b))$. Fix $a<b$ and $f
\in C^\infty_0 ((0,a))$. Define
\begin{equation} \lb{5.5}
m_a (-\kappa^2) := -\kappa - \int_0^a d\alpha \, A(\alpha)
e^{-2\alpha\kappa} 
\end{equation}
for $\Re (\kappa) \geq 0$. Fix $\kappa_0$ real and let
\[
g(y,\kappa_0, a) := m_a (-(\kappa_0 + iy)^2),
\]
with $\kappa_0, a$ as real parameters and $y\in\bbR$ a variable. As usual, 
define the Fourier transform by (initially for smooth functions and then
by duality for  tempered distributions \cite[Ch.\ IX]{RS75})
\begin{equation} \lb{5.6}
\hat F(k) = \f{1}{\sqrt{2\pi}} \int_\bbR dy \, e^{-iky} F(y), \quad
\check F(k) = \f{1}{\sqrt{2\pi}} \int_\bbR dy \, e^{iky} F(y).
\end{equation}
Then by \eqref{5.5}, 
\begin{equation}
\widehat{\bar g} (k, \kappa_0, a) = -\sqrt{2\pi}\, \kappa_0 \delta(k) - 
\sqrt{2\pi}\, \delta'(k) - \f{\sqrt{2\pi}}{2}\, e^{-k\kappa_0} 
A\left(\f{k}{2}\right) \chi_{(0,2a)} (k).
\end{equation}
Thus, since $f(0_+)=f'(0_+)=0$, in fact, $f$ has support away from $0$ and 
$a$, 
\begin{align}
\int_0^a d\alpha \,A(\alpha) f(\alpha) &= -\f{2}{\sqrt{2\pi}} \int_0^a  
 d\alpha \, 
\ol{\widehat{\bar g}}\, (2\alpha, \kappa_0, a) e^{2\alpha\kappa_0} 
f(\alpha) \notag \\
&=-\f{1}{\sqrt{2\pi}} \int_0^{2a} d\alpha \, 
\ol{\widehat{\bar g}} (\alpha, \kappa_0, a) 
e^{\alpha\kappa_0} f\left(\f{\alpha}{2}\right) \notag \\
&= -\f{1}{\sqrt{2\pi}} \int_\bbR  dy \,
g(y, \kappa_0, a) \check F (y, \kappa_0),
\lb{5.8}
\end{align}
where we have used the unitarity of $\,\, \widehat{} \,\,$ and
\begin{equation}
\check F (y, \kappa_0) = \f{1}{\sqrt{2\pi}} \int_0^{2a} d\alpha \, 
e^{\alpha (\kappa_0 +iy)} f\left(\f{\alpha}{2}\right)  
= \f{2}{\sqrt{2\pi}} \int_0^a d\alpha \, 
e^{2\alpha (\kappa_0 +iy)} f(\alpha). 
\lb{5.9}
\end{equation}
Notice that
\begin{equation} \lb{5.10}
|\check F (y, \kappa_0)| \leq Ce^{2(a-\veps)\kappa_0} (1+|y|^2)^{-1}
\end{equation}
since $f$ is smooth and supported in $(0, a-\veps)$ for some $\veps >0$.

By Theorem~\ref{t4.1} and Corollary~\ref{c4.8},
\begin{equation} \lb{5.11}
|m_a (-(\kappa_0 + iy)^2) - m(-(\kappa_0 + iy)^2)| \leq Ce^{-2a\kappa_0}
\end{equation}
for large $\kappa_0$, uniformly in $y$. From \eqref{5.8}, \eqref{5.10}, and 
\eqref{5.11}, one concludes the following fact:

\begin{lemma} \lb{l5.1} Let $f\in C_0^\infty ((0,a))$ with $0<a<b$ and 
$V\in  L^1 ([0,b))$. Then
\begin{equation} \lb{5.12}
\int_0^a  d\alpha \, A(\alpha) f(\alpha) 
= \lim_{\kappa_0 \uparrow \infty} \left[-\f{1}{\pi} \int_\bbR  dy  \, 
m(-(\kappa_0 + iy)^2)\int_0^a d\alpha \,e^{2\alpha(\kappa_0 +
iy)}  f(\alpha) \right].
\end{equation}
\end{lemma}

As a function of $y$, for $\kappa_0$ fixed, the alpha integral is 
$O((1+y^2)^{-N})$ for all $N$ because $f$ is $C^\infty$. Now define
\begin{equation} \lb{5.13}
\tilde m_R (-\kappa^2)=\left[ c_R + \int_{\lambda\leq R} 
\f{d\rho(\lambda)}{\lambda + \kappa^2}\right],
\end{equation}
where $c_R$ is chosen so that $\tilde m_R\arrow\limits_{R\to\infty} m$. Because 
$\int_\bbR \f{d\rho(\lambda)}{1+\lambda^2} <\infty$, the convergence is uniform 
in $y$ for $\kappa_0$ fixed and sufficiently large. Thus in \eqref{5.12} we can 
replace $m$ by $m_R$ and take a limit (first $R\to\infty$ and then $\kappa_0
\uparrow\infty$). Since $f(0_+)=0$, the $\int dy\, c_R\, d\alpha$-integrand is 
zero. Moreover, we can now interchange the $dy\,d\alpha$ and $d\rho(\lambda)$ 
integrals. The result is that 
\begin{align} 
\begin{split}
\int_0^a d\alpha\,A(\alpha) f(\alpha) &= \lim_{\kappa_0\uparrow\infty}\, 
\lim_{R\to\infty} \int_{\lambda\leq R} d\rho(\lambda)  \\
& \quad \times \left[ \int_0^a d\alpha\, e^{2\alpha\kappa_0} f(\alpha) 
\left[ -\f{1}{\pi} \int_\bbR 
\f{dy \,e^{2\alpha iy}}{(\kappa_0 + iy)^2 +\lambda}
\right]\right].  \lb{5.14}
\end{split}
\end{align}
In the case at hand, $d\rho$ is bounded below, say $\lambda \geq -K_0$. 
As long as we take $\kappa_0 > K_0$, the poles of $(\kappa_0 + iy)^2 + \lambda$ 
occur in the upper half-plane
\[
y_\pm = i\kappa_0 \pm\sqrt{\lambda}\, .
\]

Closing the contour in the upper plane, we find that if $\lambda \geq -K_0$,
\[
-\f{1}{\pi} \int_\bbR \f{dy\,e^{2\alpha iy}}{(\kappa_0 + iy)^2 + \lambda} 
= 
-2e^{-2\alpha\kappa_0}\, \f{\sin(2\alpha\sqrt{\lambda}\,)}{\sqrt\lambda}\,.
\]
Thus \eqref{5.14} becomes
$$
\int_0^a d\alpha \, A(\alpha)  f(\alpha) 
= -2 \lim_{\kappa_0\uparrow\infty}\, \lim_{R\to\infty} \int_{\lambda\leq
R}  d\rho(\lambda)  
\left[\int_0^a d\alpha \, f(\alpha)\, 
\f{\sin(2\alpha\sqrt\lambda\,)}{\sqrt\lambda}\right].
$$
$\kappa_0$ has dropped out and the $\alpha$ integral is bounded by $C
(1+\lambda^2)^{-1}$, so one can take the limit as $R\to\infty$ since $\int_\bbR  
\f{d\rho(\lambda)}{1+\lambda^2} <\infty$. One is therefore led to the 
following result.

\begin{theorem} \lb{t5.2} Let $f\in C^\infty_0 ((0,a))$ with $a<b$ and either 
$b<\infty$ or $V\in L^1 ([0, \infty))$ with $b=\infty$. Then
\begin{equation} \lb{5.15}
\int_0^a d\alpha \, A(\alpha) f(\alpha) = -2 \int_\bbR 
d\rho(\lambda) \left[ \int_0^a 
 d\alpha \, 
f(\alpha)\, \f{\sin(2\alpha\sqrt\lambda\,)}{\sqrt\lambda}\right].
\end{equation}
\end{theorem}

One can strengthen this in two ways. First, one wants to allow $a>b$ 
if $b<\infty$. As long as $A$ is interpreted as a distribution with 
$\delta$ and $\delta'$ functions at $\alpha=nb$, this is easy. One also
wants to allow $f$ to have a nonzero derivative at $\alpha =0$. The net
result is described in the next theorem:

\begin{theorem} \lb{t5.3} Let $f\in C^\infty_0 (\bbR)$ with $f(-\alpha) = 
-f(\alpha)$, $\alpha\in\bbR$ and either $b<\infty$ or 
$V\in L^1 ([0,\infty))$  with $b=\infty$. Then
\begin{equation} \lb{5.16}
-2\int_\bbR d\rho(\lambda) \left[ \int_{-\infty}^\infty 
d\alpha \, f(\alpha) \, 
\f{\sin(2\alpha\sqrt\lambda\, )}{\sqrt\lambda} \right] 
= \int_{-\infty}^\infty d\alpha \, \tilde A(\alpha) f(\alpha),
\end{equation}
where $\tilde A$ is the distribution
\begin{subequations} \lb{5.17}
\begin{align}
\tilde A(\alpha) &= \chi_{(0,\infty)} (\alpha) A(\alpha) - 
\chi_{(-\infty, 0)} (\alpha) A(-\alpha) + \delta'(\alpha) \lb{5.17a} \\
\intertext{if $b=\infty$ and}
\tilde A(\alpha) &= \chi_{(0,\infty)} (\alpha) A(\alpha) - \chi_{(-\infty,
0)}  (\alpha) A(-\alpha) + \delta' (\alpha) \notag \\
&\quad + \sum_{j=1}^\infty B_j [\delta (\alpha - 2bj) - \delta (\alpha +
2bj)] \notag \\
&\quad + \sum_{j=1}^\infty \tfrac12\, A_j [\delta' (\alpha - 2bj) +
\delta'  (\alpha + 2bj)] \lb{5.17b}
\end{align}
\end{subequations}
if $b<\infty$, where $A_j, B_j$ are $h$ dependent and given in 
Theorems \ref{t4.3} and \ref{t4.7}.
\end{theorem}

\medskip

Next we change the subject temporarily and turn to bounds on $\int_0^{\pm R} d\rho(\lambda)$ which are of independent interest: 
As we will see, \eqref{1.11} implies asymptotic results on $\int_{-R}^R 
d\rho(\lambda)$, and \eqref{5.1} will show that $\int_{-\infty}^0 
e^{b\sqrt{-\lambda}}d\rho(\lambda)<\infty$ for all $b>0$ and more. It 
follows from \eqref{5.3} that
\[
\Im (m(ia)) = a\int_\bbR \f{d\rho(\lambda)}{\lambda^2 + a^2}\,, 
\quad a>0.
\]
Thus, Everitt's result \eqref{1.11} implies that
\[
\lim_{a\to\infty} a^{\f12} \int_\bbR \f{d\rho(\lambda)}{\lambda^2 + a^2} = 
2^{-\f12}.
\]
Standard Tauberian arguments (see, e.g., in \cite[Sect.\ III.10]{Si79}, 
which in this case shows that on even functions $R^{\f32} d\rho
(\lambda/R)\arrow\limits_{R\to\infty} (2 \pi)^{-1}
|\lambda|^{\f12}\, d\lambda$) then  imply the following result: 

\begin{theorem} \lb{t6.1}
\begin{equation} \lb{6.1}
\lim_{R\to\infty} R^{-\f32} \int_{-R}^R d\rho(\lambda) = \f{2}{3\pi}\, .
\end{equation}
\end{theorem}

\begin{remark} $(i)$ This holds in all cases (1--4) we consider here, 
including some with $\supp(d\rho)$ unbounded below. \\
$(ii)$ Since one can show that $\int_{-\infty}^0 d\rho$ is bounded, one can
replace $\int_{-R}^R$ by $\int_0^R$ in \eqref{6.1}.
\end{remark}

Next, we recall the following a priori bound that follows from 
Lemmas \ref{p3.1} and \ref{p3.2}:

\begin{lemma} \lb{p6.2} Let $d\rho$ be the spectral measure for a 
Schr\"odinger operator in Cases 1--4. Fix $a <b$. Then there is a constant 
$C_a$ depending only on $a$ and $\int_0^a dy \, |V(y)|$ so that
\begin{equation} \lb{6.2}
\int_\bbR \f{d\rho(\lambda)}{1+\lambda^2} \leq C_a.
\end{equation}
\end{lemma}

The goal is to bound $\int_{-\infty}^0 e^{2\alpha\sqrt{-\lambda}}$ $d\rho
(\lambda)$ for any $\alpha < b$ and to find  an explicit bound in terms of
$\sup_{0\leq x \leq\alpha + 1}$ $[-V(y)]$ when   that $\sup$ is finite. As
a preliminary, we need the following result from the  standard limit
circle theory \cite[Sect.\ 9.4]{CL55}.

\begin{lemma}\lb{p6.3} Let $b=\infty$ and let $d\rho$ be the spectral 
measure for some Schr\"odinger operator in Cases 2--4. Let $d\rho_{R,h}$ be the spectral measure 
for the problem with $b=R<\infty, h$ and potential equal to $V(x)$ for
$x\leq  R$. Then there exists $h(R)$ so that
\[
d\rho_{R, h(R)}\arrow_{R\to\infty} d\rho,
\]
when smeared with any function $f$ of compact support. In particular, if 
$f\geq 0$, then
\[
\int_\bbR d\rho(\lambda) \, f(\lambda) \leq 
\operatornamewithlimits{\varlimsup}_{R\to\infty} 
\int_\bbR d\rho_{R, h(R)} (\lambda) \, f(\lambda).
\]
\end{lemma}

This result implies that we need only obtain bounds for $b<\infty$ (where 
we have already proved \eqref{5.15}).

\begin{lemma} \lb{l6.4} If $\rho_1$ has support in $[-E_0, \infty)$, $E_0 >0$, 
then
\begin{equation} \lb{6.3}
\int_{-\infty}^0 e^{\gamma\sqrt{-\lambda}} \, d\rho_1 (\lambda) \leq 
e^{\gamma\sqrt{E_0}} (1+ E_0^2) \int_{-\infty}^0 \f{d\rho_1 (\lambda)}
{1+\lambda^2}\,.
\end{equation}
\end{lemma}

Lemmas {\ref{p6.2}, \ref{p6.3} and Lemma \ref{l6.4} imply the following result.

\begin{theorem} \lb{t6.5} Let $\rho$ be the spectral measure for some Schr\"odinger operator in Cases 2--4. Let
\begin{align*} 
&E(\alpha_0) 
:= -\inf \bigg\{ \int_0^{\alpha_0 +1} dx \, (|\varphi'_n(x)|^2 + V(x)
|\varphi(x)|^2) \biggm| \varphi\in 
C^\infty_0 ((0,\alpha_0 +1)),  \\
&\hspace*{2.6cm}  \int_0^{\alpha_0 +1} dx\,|\varphi(x)|^2 \leq 1 \bigg\}.
\end{align*}
Then for all $\delta >0$ and $\alpha_0 >0$,
\begin{equation} 
 \alpha_0 \delta \int_{-\infty}^0 e^{2(1-\delta)\alpha_0\sqrt{-\lambda}} 
 \, d\rho (\lambda) 
\leq \biggl[ C_1 (1+\alpha_0)  + C_2 (1+E(\alpha_0)^2) 
e^{2(\alpha_0 +1)\sqrt{E(\alpha_0)}} \, \biggr],  \lb{"6.7"}
\end{equation}
where $C_1, C_2$ only depend on $\int_0^1 dx\,|V(x)|$. In particular,
\begin{equation} \lb{"6.8"}
\int_{-\infty}^0 e^{B\sqrt{-\lambda}} \, d\rho(\lambda) <\infty 
\end{equation}
for all $B<\infty$.
\end{theorem}

As a special case, suppose $V(x) \geq -C(x+1)^2$. Then $E(\alpha_0) \geq 
-C(\alpha_0 + 2)^2$ and we see that
\begin{equation} \lb{6.7}
\int_{-\infty}^0 e^{B \sqrt{-\lambda}}  \, d\rho(\lambda) \leq D_1 
e^{D_2 B^2}.
\end{equation}
This implies the next result.

\begin{theorem} \lb{t6.6} If $d\rho$ is the spectral measure for a potential 
which satisfies
\begin{equation} \lb{6.8}
V(x) \geq - Cx^2, \quad x\geq R
\end{equation}
for some $R>0$, $C>0$, then for $\veps>0$ sufficiently small, 
\begin{equation} \lb{6.9}
\int_{-\infty}^0  e^{-\veps\lambda} \, d\rho(\lambda) < \infty.
\end{equation}
\end{theorem}

If in addition $V\in L^1([0,\infty))$, then the 
corresponding Schr\"odinger operator is bounded from below and hence 
$d\rho$ has compact support on $(-\infty, 0]$. This fact will be useful
later in the  scattering-theoretic context.

The estimate \eqref{"6.8"}, in the case of non-Dirichlet boundary 
conditions at $x=0_+$, appears to be due to Marchenko \cite{Ma73}. Since
it is a  fundamental ingredient in the inverse spectral problem, it
generated  considerable attention; see, for instance, \cite{GL55},
\cite{Le73}, \cite{Le73a}, \cite{Le77}, \cite{LG64}, \cite{Ma73},
\cite[Sect.\ 2.4]{Ma86}. The case of a  Dirichlet boundary at $x=0_+$
was studied in detail by Levitan \cite{Le77}. These  authors, in addition
to studying the spectral asymptotics of $\rho(\lambda)$ as 
$\lambda\downarrow -\infty$, were also particularly interested in the 
asymptotics of $\rho(\lambda)$ as $\lambda\uparrow \infty$ and
established Theorem \ref{t6.1}. In the latter context,
we also refer to Bennewitz \cite{Be89},  Harris \cite{Ha97}, and the
literature cited therein. In contrast to these activities,  we were not
able to find estimates of the type \eqref{"6.7"} (which implies 
\eqref{"6.8"}) and \eqref{6.9} in the literature.

At this point one can return to the relation between $A$ and $\rho$ and
discuss a second distributional form of this relation which extends
Theorem \ref{t5.2} to all four cases.

\begin{theorem} \lb{t7.1} Let $f\in C^\infty_0 ((0, \infty))$ and suppose  
$b=\infty$. Assume $V$ satisfies \eqref{1.3} and let $d\rho$ be the 
associated spectral measure and $A$ the associated $A$-function. Then
\eqref{5.16} and \eqref{5.17} hold.
\end{theorem}

Next we establish a third relation between $A$ and $\rho$ and turn to
Abelian limits: 

For $f\in C^\infty_0 (\bbR)$, define for $\lambda\in\bbR$, 
\begin{equation} \lb{8.1}
Q(f)(\lambda) = \int_{-\infty}^\infty d\alpha \, f(\alpha) 
\f{\sin(2\alpha \sqrt\lambda\,)}{\sqrt\lambda} 
\end{equation}
and then
\begin{align}
T(f) &= -2 \int_\bbR d\rho(\lambda) \,Q(f)(\lambda) \lb{8.2} \\
&= \int_{-\infty}^\infty d\alpha \, \tilde A(\alpha) f(\alpha). \lb{8.3}
\end{align}
Relations \eqref{5.16}, \eqref{5.17} show that for 
$f\in C^\infty_0 (\bbR)$, the two expressions \eqref{8.2}, \eqref{8.3}
define the same $T(f)$. This was proved for odd $f$'s but both
integrals vanish for even $f$'s. Now one wants to use \eqref{8.2} to extend to a
large class of $f$, but needs to exercise some care not  to use \eqref{8.3},
except for $f\in C^\infty_0 (\bbR)$.

$Q(f)$ can be defined as long as $f$ satisfies
\begin{equation} \lb{8.4} 
|f(\alpha)| \leq C_k e^{-k|\alpha|}, \quad \alpha\in\bbR
\end{equation}
for all $k >0$. In particular, a simple calculation shows that
\begin{equation} \lb{8.5}
f(\alpha) = (\pi\veps)^{-\f12} \left[ e^{-(\alpha-\alpha_0)^2/ \veps}  
\right] \, \text{ implies } \, 
 Q(f)(\lambda) = \f{\sin(2\alpha_0 \sqrt\lambda\, )}
{\sqrt\lambda}\, e^{-\veps\lambda}.
\end{equation}
We use $f(\alpha, \alpha_0, \veps)$ for the function $f$ in \eqref{8.5}.

For $\lambda \geq 0$, repeated integrations by parts show that
\begin{equation} \lb{8.6}
|Q(f)(\lambda)| \leq C(1+\lambda^2)^{-1} \left[ \|f\|_1 + 
\left\| \f{d^3f}{d\alpha^3} \right\|_1 \right], 
\end{equation}
where $\| \dott\|_1$ represents the $L^1 (\bbR)$-norm. Moreover, 
essentially by repeating the calculation that led to \eqref{8.5}, one sees
that for $\lambda 
\leq 0$,
\begin{equation} \lb{8.7}
|Q(f)(\lambda)| \leq Ce^{\veps |\lambda|} \bigl\| e^{+\alpha^2/ \veps} 
f\bigr\|_\infty.
\end{equation}

One then concludes the following result. 

\begin{lemma} \lb{p8.1} If $\int_\bbR (1+\lambda^2)^{-1} \, d\rho(\lambda) < 
\infty$ {\rm(}always true{\rm!)} and $\int_{-\infty}^0 
e^{-\veps_0 \lambda} \, d\rho(\lambda) <\infty$ {\rm(}see Theorem~\ref{t6.6} 
and the remark following its proof{\rm)}, then using \eqref{8.2},
$T(\dott)$ can be  extended to functions $f\in C^3(\bbR)$ that satisfy
$e^{\alpha^2/\veps_0}f\in L^\infty  (\bbR)$ for some $\veps_0 >0$ and
$\f{d^3 f}{d\alpha^3}\in L^1 (\bbR)$, and  moreover,
\begin{equation}
|T(f)| \leq C \left[ \, \biggl\| \f{d^3 f}{d\alpha^3}\biggr\|_1 + 
\bigl\| e^{\alpha^2/ \veps_0} f\bigr\|_\infty \right] 
 := C ||| f |||_{\veps_0}. \lb{8.8} 
\end{equation}
\end{lemma}

Next, fix $\alpha_0$ and $\veps_0>0$ so that 
$\int_{-\infty}^0 e^{-\veps_0 \lambda} \, 
d\rho(\lambda) < \infty$. If $0<\veps < \veps_0$, 
$f(\alpha,  \alpha_0, \veps)$ 
satisfies $||| f|||_{\veps_0} <\infty$ so we can define $T(f)$. Fix $g\in 
C^\infty_0 (\bbR)$ with $g:= 1$ on $(-2\alpha_0, 2\alpha_0)$. Then 
$||| f(\dott, \alpha_0, \veps)(1-g)|||_{\veps_0} \to 0$ as 
$\veps \downarrow 0$. So
\[
\lim_{\veps\downarrow 0} T(f(\dott , \alpha_0, \veps)) 
= \lim_{\veps \downarrow 0} 
T(gf(\dott, \alpha_0, \veps)).
\]
For $gf$, we can use the expression \eqref{8.3}. $f$ is approximately 
$\delta (\alpha - \alpha_0)$ so standard estimates show if $\alpha_0$ is a
point of  Lebesgue continuity of $\tilde A(\alpha)$, then 
\[
\int_{-\infty}^\infty  d\alpha \, f(\alpha, \alpha_0, \veps) g(\alpha) 
\tilde A(\alpha) \arrow_{\veps\downarrow 0} \tilde A (\alpha_0).
\]
Since $A-q$ is continuous, points of Lebesgue continuity of $A$ exactly are 
points of Lebesgue continuity of $V$. Thus, one obtains the following
theorem.

\begin{theorem} \lb{t8.2} Suppose either $b<\infty$ and $V\in L^1 ([0,b])$ 
or $b=\infty$, and then either $V\in L^1 ([0,\infty))$ or 
$V\in L^1 ([0,a])$ for all $a>0$ and
\[
V(x) \geq -Cx^2, \quad x\geq R
\]
for some $R>0$, $C>0$. Let $\alpha_0 \in (0,b)$ and be a point of Lebesgue 
continuity of $V$. Then
\begin{equation} \lb{8.9}
A(\alpha_0) = - 2\lim_{\veps\downarrow 0} \int_\bbR  d\rho(\lambda) \, 
e^{-\veps\lambda} 
\f{\sin(2\alpha_0 \sqrt\lambda \,)}{\sqrt\lambda}.
\end{equation}
\end{theorem}

\medskip

Finally, we specialize \eqref{8.9} to the scattering-theoretic setting. 
Assuming $V\in L^1 ([0,\infty); (1+x)\, dx)$, the corresponding Jost
solution $f(z,x)$  is defined by 
\begin{equation} \lb{8.11}
f(z,x) = e^{i\sqrt{z}\,x} - \int_x^\infty dx' \, \f{\sin(\sqrt{z}\, (x-x'))}
{\sqrt z}\, V(x') f(z,x'), \quad \Im(\sqrt z\, )\geq 0,
\end{equation}
and the corresponding Jost function, $F(\sqrt z\,)$, and scattering matrix, 
$S(\lambda)$, $\lambda \geq 0$, then read
\begin{align}
F(\sqrt z\, ) &= f(z,0_+), \lb{8.12} \\
S(\lambda) &= \overline{F(\sqrt\lambda\, )}/F(\sqrt\lambda\,), \quad 
\lambda \geq 0. \lb{8.13}
\end{align}
The spectrum of the Schr\"odinger operator $H$ in $L^2 ([0,\infty))$ 
associated with the differential expression $-\f{d^2}{dx^2} + V(x)$ and a
Dirichlet boundary  condition at $x=0_+$ (cf.\ \eqref{8.1.5} for precise details) is simple and of the type
\[
\sigma(H) = \{ -\kappa^2_j <0\}_{j\in J} \cup [0,\infty).
\]
Here $J$ is a finite (possibly empty) index set, $\kappa_j >0$, $j\in J$, 
and the essential spectrum is purely absolutely continuous. The
corresponding spectral  measure explicitly reads
\begin{equation} \lb{8.14}
d\rho(\lambda) = \begin{cases}
\pi^{-1} |F(\sqrt\lambda\,)|^{-2}\sqrt\lambda\, d\lambda, 
&\lambda \geq 0, \\
\sum_{j\in J} c_j \delta (\lambda + \kappa^2_j)\, d\lambda, &\lambda <0,
\end{cases}
\end{equation}
where
\begin{equation} \lb{8.15}
c_j = \|\varphi (-\kappa^2_j,\dott)\|^{-2}_2, \quad j\in J
\end{equation}
are the norming constants associated with the eigenvalues $\lambda_j = 
-\kappa^2_j <0$. Here the regular solution $\varphi(z,x)$ of 
$-\psi''(z,x)+[V(x)-z]\psi(z,x)=0$ (defined by
$\varphi(z,0_+)=0$, $\varphi'(z,0_+)=1$) and $f(z,x)$ in \eqref{8.11} are
linearly dependent precisely for $z=-\kappa^2_j$, $j\in J$.

Since
\[
|F(\sqrt\lambda\, )| = \prod_{j\in J} 
\left( 1+\f{\kappa^2_j}{\lambda}\right) 
\exp\left(\f{1}{\pi}\, P\int_0^\infty \f{d\lambda' \, \delta(\lambda')}
{\lambda - \lambda'}\right), \quad \lambda \geq 0,
\]
where $P\int_0^\infty$ denotes the principal value symbol and 
$\delta(\lambda)$ 
the corresponding scattering phase shift, that is, 
$S(\lambda) =\exp(2i \delta
(\lambda))$, $\delta(\lambda)\arrow\limits_{\lambda\uparrow\infty} 0$, the 
scattering data
\[
\{-\kappa^2_j, c_j\}_{j\in J} \cup \{S(\lambda)\}_{\lambda \geq 0}
\]
uniquely determine the spectral measure \eqref{8.14} and hence $A(\alpha)$. 
Inserting \eqref{8.14} into \eqref{8.9} then yields the following expression for 
$A(\alpha)$ in terms of scattering data.

\begin{theorem} \lb{t8new} Suppose that $V\in L^1 ([0,\infty);
(1+x)dx)$.  Then
\begin{equation} \lb{8.16}
\begin{split}
A(\alpha) &= -2\sum_{j\in J} c_j \kappa^{-1}_j \sinh(2\alpha\kappa_j) \\
&\quad  -2\pi^{-1} \lim_{\veps\downarrow 0}\int_0^\infty 
d\lambda \, 
e^{-\veps\lambda}  |F(\sqrt\lambda)|^{-2} \sin(2\alpha \sqrt\lambda\,)
\end{split}
\end{equation}
at points $\alpha \geq 0$ of Lebesgue continuity of $V$.
\end{theorem}

\begin{remark} In great generality $|F(k)|\to 1$ as $k\to\infty$, 
so one cannot take the limit in $\veps$ inside the integral in \eqref{8.16}. 
In general, though, one can can replace $|F(\sqrt\lambda\,)|^{-2}$ by 
$(|F(\sqrt\lambda\,)|^{-2} -1) \equiv X(\lambda)$ and ask if one can take 
a limit there. As long as $V$ is $C^2 ((0,\infty))$ with $V'' \in L^1 
([0,\infty))$, it is not hard to see that as $\lambda\to\infty$
\[
X(\lambda) = -\f{V(0)}{2\lambda} + O(\lambda^{-2}).
\]
Thus, if $V(0) = 0$, then
\begin{equation} \lb{8.17}
\begin{split} 
A(\alpha) &= -2 \sum_{j\in J} c_j \kappa^{-1}_j \sinh (2\alpha \kappa_j)
\\ &\quad -2\pi^{-1} \int_0^\infty  d\lambda \, (|F(\sqrt\lambda\,)|^{-2}-1) 
\sin(2\alpha \sqrt\lambda\,).
\end{split}
\end{equation}
The integral in \eqref{8.17} is only conditionally convergent if 
$V(0)\neq 0$.
\end{remark}
 
We note that in the present case, where $V\in L^1 ([0,\infty);
(1+x)\,dx)$,  the representation \eqref{1.15} of the $m$-function in terms
of the $A(\alpha)$-amplitude was considered in a paper by Ramm \cite{Ra87} 
(see also \cite[p.\ 288--291]{Ra92}).

We add a few more remarks in the scattering-theoretic setting. Assuming
$V\in L^1 ([0,\infty); (1+x)\, dx)$, one sees that
\begin{equation} \lb{"9.17"}
|F(k)| \eqlim_{k\uparrow \infty} 1 + o(k^{-1}) 
\end{equation}
(cf.\ \cite[eq.\ II.4.13]{CS89} and apply the Riemann-Lebesgue lemma; 
actually,  one only needs $V\in L^1 ([0,\infty))$ for the asymptotic
results on $F(k)$ as 
$k\uparrow\infty$ but we will ignore this refinement in the following). A 
comparison of \eqref{"9.17"} and \eqref{8.16} then clearly demonstrates the 
necessity of an Abelian limit in \eqref{8.16}. Even replacing $d\rho$ in 
\eqref{8.9} by $d\sigma = d\rho - d\rho^{(0)}$, that
is,  effectively replacing $|F(\sqrt\lambda\,)|^{-2}$ by
$[|F(\sqrt\lambda\,)|^{-2}  -1]$ in \eqref{8.16}, still does not
necessarily produce an absolutely convergent  integral in \eqref{8.16}. 

The latter situation changes upon increasing the smoothness properties of 
$V$ since, for example, assuming $V\in L^1 ([0,\infty); (1+x)\, dx)$,
$V'\in L^1  ([0,\infty))$, yields
\[
|F(k)|^{-2} - 1 \eqlim_{k\uparrow \infty} O(k^{-2})
\]
as detailed high-energy considerations (cf.~\cite{GPT80}) reveal. Indeed,
if $V'' \in L^1 ([0,\infty))$, then the  integral one gets is absolutely
convergent if and only if
$V(0)=0$.

\medskip

As a final issue related to the representation \eqref{5.1}, we discuss the
issue  of bounds on $A$ when $|V(x)| \leq Cx^2$. One has two general
bounds on $A$: the  estimate of \cite{Si99} (see \eqref{1.14}), 
\begin{equation} \lb{9.16}
|A(\alpha) - V(\alpha)| \leq \left[ \int_0^\alpha dy \, |V(y)| \right]^2 
\exp \left[ \alpha\int_0^\alpha dy \, |V(y)| \right] ,
\end{equation}
and the estimate in Theorem~\ref{t10.2}, 
\begin{equation} \lb{9.17}
|A(\alpha)| \leq \f{\gamma(\alpha)}{\alpha} \, 
I_1 (2\alpha \gamma(\alpha)),
\end{equation}
where $|\gamma(\alpha)| = \sup_{0\leq x\leq\alpha} |V(x)|^{1/2}$ and 
$I_1 (\dott)$ is the modified Bessel function of order one (cf., e.g.,
\cite{AS}, Ch.~9). Since  (\cite{AS}, p.~375) 
\begin{equation} \lb{9.18}
0\leq I_1 (x) \leq e^x , \quad x\geq 0,
\end{equation}
one concludes that
\begin{equation} \lb{9.19}
|A(\alpha)| \leq \sqrt C\, e^{2\sqrt C\, \alpha^2}
\end{equation}
if $|V(x)|\leq Cx^2$. 

We continue with a discussion of the case of constant $V$: 

\begin{example} If $b=\infty$ and $V(x)=V_0$, $x\geq 0$, then if $V_0 >0$, 
\begin{equation} \lb{10.1}
A(\alpha) = \f{V^{1/2}_0}{\alpha}\, J_1 (2\alpha V^{1/2}_0), 
\end{equation}
where $J_1 (\dott)$ is the Bessel function of order one {\rm(}cf., e.g., 
\cite{AS}, Ch.~9{\rm )}; and if $V_0 < 0$, 
\begin{equation} \lb{10.2}
A(\alpha) = \f{(-V_0)^{1/2}}{\alpha}\, I_1 (2\alpha (-V_0)^{1/2}),
\end{equation}
with $I_1 (\dott) $ the corresponding modified Bessel function.
\end{example}

This example is important because of the following monotonicity
property:

\begin{theorem} \lb{t10.2} Let $|V_1(x)| \leq - V_2(x)$ on $[0,a]$ with 
$a\leq \min(b_1, b_2)$. Then,
$$
|A_1 (\alpha)| \leq -A_2 (\alpha) \, \text{ on } \, [0,a].
$$
In particular, for any $V$ satisfying $\sup_{0\leq x \leq\alpha} |V(x)| < 
\infty$, one obtains  
\begin{equation} \lb{10.5a}
|A(\alpha)| \leq \f{\gamma(\alpha)}{\alpha}\, I_1 (2\alpha \gamma(\alpha)),  
\end{equation}
where
\begin{equation}\lb{10.5b}
\gamma(\alpha) = \sup_{0\leq x\leq\alpha} (|V(x)|^{1/2}).  
\end{equation}
In particular, \eqref{9.18} implies
\begin{equation} \lb{10.6}
|A(\alpha)| \leq \alpha^{-1} \gamma(\alpha) e^{2\alpha\gamma(\alpha)},
\end{equation}
and if $V$ is bounded,
\begin{equation} \lb{10.7}
|A(\alpha)| \leq \alpha^{-1} \|V\|^{1/2}_\infty 
\exp(2\alpha \|V\|^{1/2}_\infty).
\end{equation}
\end{theorem}

For $\alpha$ small, \eqref{10.6} is a poor estimate and one should use 
\eqref{9.16} which implies that $|A(\alpha) \leq \|V\|_\infty + \alpha^2 
\|V\|^2_\infty e^{\alpha^2\|V\|_\infty}$.

This lets one prove the following result:

\begin{theorem} \lb{t10.3} Let $h=\infty$ and $V\in L^\infty
([0,\infty))$.  Suppose 
$\kappa^2 > \|V\|_\infty$. Then
\begin{equation} \lb{10.9}
m(-\kappa^2) = -\kappa - \int_0^\infty d\alpha \, A(\alpha) e^{-2\alpha\kappa} 
\end{equation}
{\rm (}with an absolutely  convergent integral and no error term{\rm )}.
\end{theorem} 

\begin{remark} 
We recall (cf.\ \eqref{1.15}) that the representation \eqref{10.9} also
holds with $A\in L^1([0,a])$ for all $a>0$ and as an absolutely convergent
integral for $\Re(\kappa) > \|V\|_1/2$ if $V\in L^1([0,\infty))$. This
fact will be used below.
\end{remark}

The case of Bargmann potentials has
been worked out in
\cite[Sect.\ 11]{GS00a} and explicit formulas for the $A$-function have
been obtained. 

We end this survey of \cite{GS00a} and \cite{Si99} by recalling the major thrust
of \cite{Si99} -- the connection  between $A$ and the inverse spectral
theory. Namely, there is an $A(\alpha, x)$ function associated to $m(z,x)$ by
\begin{equation} \lb{1.20}
m(-\kappa^2, x) = -\kappa - \int_0^a d\alpha \, 
A(\alpha, x) e^{-2\alpha\kappa}  +
\tilde O(e^{-2\alpha\kappa})
\end{equation}
for $a<b -x$. This, of course, follows from Theorem \ref{t1} by
translating  the origin. The point is that $A$ satisfies the simple
differential equation  in distributional sense
\begin{equation} \lb{1.21}
\f{\partial A}{\partial x}\, (\alpha, x) 
= \f{\partial A}{\partial\alpha} \, 
(\alpha, x) + \int_0^a d\beta \, A(\alpha-\beta, x) A(\beta,x).
\end{equation}
This is proved in \cite{Si99} for $V\in L^1 ([0, a])$ (and some other 
$V$'s) and so holds in the generality of \cite{GS00a} since Theorem \ref{t3}
implies $A(\alpha,x)$ for $\alpha+x \leq a$ is only a function of $V(y)$
for $y\in[0,a]$.

Moreover, by \eqref{1.14}, one has
\begin{equation} \lb{1.22}
\lim_{\alpha\downarrow 0} |A(\alpha, x) - V(\alpha + x)| =0
\end{equation}
uniformly in $x$ on compact subsets of the real line, so by the uniqueness 
theorem for solutions of \eqref{1.21} (proved in \cite{Si99}), $A$ on 
$[0,a]$ determines $V$ on $[0,a]$. 

In the limit circle case, there is an additional issue to discuss. Namely, 
that $m(z, x=0)$ determines the boundary condition at $\infty$. This is 
because, as we just discussed, $m$ determines $A$ which determines $V$ on 
$[0,\infty)$. $m(z, 0_+)$ and $V$ determine $m(z,x)$ by the Riccati 
equation. Once we know $m$, we can recover $u(i, x) = \exp 
\big(\int_0^x m(i, y)\, dy\big)$, and so the particular solution that
defined  the boundary condition at $\infty$.

Thus, the inverse spectral theory aspects of the framework easily extend to 
the general case of potentials considered in \cite{GS00a}.

To turn this into an inverse spectral approach alternative to and fully equivalent to that of Gel'fand and Levitan, one needs to settle necessary and sufficient conditions for solvability of the differential equation \eqref{1.21} in terms of an initial condition 
$A(\alpha, 0_+)=A_0(\alpha)$, that is, in terms of properties of $A_0$. This final step was accomplished by Remling \cite{Re03} and we briefly describe its major elements next.

Remling's first result is of local nature and determines a necessary and sufficient condition on $A$ to be the $A$-function of a potential $V$. Assuming $V\in L^1([0,b])$ for all $b>0$, he introduces the set
\begin{equation}
\cA_b=\{A\in L^1([0,b])\,|\, \text{$A$ real-valued}, \, I+\cK_A >0 \},   \lb{1.23}
\end{equation}
where
\begin{align*}
& (\cK_A f)(\alpha)=\int_0^b d\beta \, K(\alpha,\beta) f(\beta), \quad 
\alpha \in [0,b], \; f\in L^2([0,b]), \\
& K(\alpha, \beta)=[\phi(\alpha-\beta)-\phi(\alpha+\beta)]/2, \quad 
\phi(\alpha)=\int_0^{|\alpha|/2} d\gamma \, A(\gamma), \; \alpha, \beta \in [0,b].
\end{align*}

Based on his reformulation of the Gel'fand--Levitan approach in terms of de Branges spaces in \cite{Re02}, Remling obtained the following characterization of $A$-functions:

\begin{theorem}
$\cA_b$ is precisely the set of $A$-functions in 
$$
m(-\kappa^2) = -\kappa - \int_0^a d\alpha \,  A(\alpha) e^{-2\alpha\kappa}  +
\tilde O(e^{-2\alpha\kappa}) \, \text{ for all $a<b$.}
$$
Equivalently, given $A_0\in L^1([0,b])$, there exists a potential 
$V\in L^1([0,b])$ such that $A_0$ is the $A$-function of $V$ if and only if 
$A_0\in\cA_b$.
\end{theorem}
(We recall that all potentials $V$ in this survey are assumed to be real-valued.)

As a second result, Remling also proved in \cite{Re03} that the positivity condition in \eqref{1.23} is necessary and sufficient to solve \eqref{1.21} on 
$\Delta_b=\{(\alpha,x)\in\bbR^2\,|\, \alpha\in [0,b-x], \, x\in [0,b] \}$ given an initial condition $A(\cdot,0_+)=A_0 \in L^1([0,b])$. The potential $V$ can then be read off from 
\begin{equation}
V(x)=A(0_+,x) \, \text{ for $x\in [0,b]$.}  \lb{1.24} 
\end{equation}
Necessity of this positivity condition had been established independently by Keel and Simon (unpublished). To make this precise, it pays to slightly rewrite 
\eqref{1.21} as follows: Let
\begin{equation}
B(\alpha,x)=A(\alpha-x,x)-A_0(\alpha), \quad (\alpha, x) \in\wti \Delta_b,  \lb{1.25}
\end{equation}
where
$$
\wti \Delta_b = \{(\alpha,x)\in\bbR^2\,|\, 0\leq x\leq \alpha \leq b\}. 
$$
Then \eqref{1.21} together with the initial condition $A(\cdot,0_+)=A_0 \in L^1([0,b])$, becomes
\begin{align}
& B(\alpha,x)=\int_0^x dy \int_0^{\alpha -y} d\beta \, [B(y+\beta,y)+A_0(y+\beta)] 
[B(\alpha-\beta,y)+A_0(\alpha-\beta)]    \no \\
& B(\alpha,0_+)=0, \quad (\alpha,x)\in \wti \Delta_b.  \lb{1.26}
\end{align}

If $A$ is actually the $A$-function of a potential, then $B\in C(\wti \Delta_b)$ by \cite[Theorem\ 2.1]{Si99}. Remling \cite{Re03} then proves the following result:

\begin{theorem}
Suppose $A_0\in L^1([0,b])$. Then \eqref{1.26} has a solution $B\in C(\wti \Delta_b)$ if and only if $A_0\in \cA_b$.
\end{theorem}

This brings Simon's inverse approach to full circle and one can envision the following two scenarios. First, Simon's inverse $A$-function approach, as complemented by Remling \cite{Re03}:

\vspace{5pt} 
\fbox{\begin{minipage}{330pt}
\begin{align}
& A_0\in \cA_b 
\xrightarrow{\text{by \eqref{1.26}}} B(\alpha,x), \; (\alpha,x)\in\wti\Delta_b 
\xrightarrow{\text{by \eqref{1.25}}} A(\alpha,x), \; (\alpha,x)\in\Delta_b   \no \\
& \quad \xrightarrow{\text{by \eqref{1.24}}} V=A(0_+,\cdot)\in L^1([0,b]).  \lb{1.27}
\end{align}
\end{minipage}} 
\vspace{10pt}

Second, denote by $\cR$ the set of spectral functions $\rho$ associated with self-adjoint half-line Schr\"odinger operators with a Dirichlet boundary condition at $x=0$ and a self-adjoint boundary condition \eqref{1.5} at infinity (if any, i.e., if \eqref{1.1} is in the limit circle case at $\infty$). For characterizations of $\cR$ we refer, for instance, to 
\cite{LG64}, \cite[Ch.\ 2]{Le87}, \cite[Ch.\ 2]{Ma86}. Then combining \eqref{1.27} with \eqref{5.16} yields Simon's inverse spectral approach as an alternative to that by Gel'fand and Levitan:

\vspace{5pt} 
\fbox{\begin{minipage}{330pt}
\begin{align}
\begin{split}
& \rho \in \cR \xrightarrow{\text{by \eqref{5.16}}} A_0\in \cA_b \, \text{ for all $b>0$ } \\
& \quad \xrightarrow{\text{by \eqref{1.26}}} B(\alpha,x), \; (\alpha,x)\in\wti\Delta_b  
\, \text{ for all $b>0$}   \\
& \qquad \xrightarrow{\text{by \eqref{1.25}}} A(\alpha,x), \; (\alpha,x)\in\Delta_b  
\, \text{ for all $b>0$} \\
& \qquad \quad \xrightarrow{\text{by \eqref{1.24}}} V=A(0_+,\cdot)\in L^1([0,b]) 
\, \text{ for all $b>0$}.  
\end{split}
\end{align} 
\end{minipage}} 
\vspace{10pt}

\medskip

{\it More recent references:} Local solvability and a necessary condition for global solvability of the $A$-equation \eqref{1.21} were recently discussed by Zhang \cite{Zh03}, \cite{Zh06}. Connections between the $A$-amplitude and the scattering transform for Schr\"odinger operators on the real line have been discussed by Hitrik \cite{Hi01}.

\smallskip
\begin{center}
$\bigstar$  \qquad $\bigstar$ \qquad $\bigstar$
\end{center}
\smallskip

Next we briefly quote the main results by Ramm and Simon \cite{RS00}. The
primary goal in this paper was to study $A$ as an interesting  object in
its own right and, in particular, using ideas implicit in  Ramm
\cite{Ra87} to obtain detailed information on the behavior of $A(\alpha)$ 
as $\alpha\to\infty$ when $V$ decays sufficiently fast as $x\to\infty$. 
Indeed, for potentials decaying rapidly enough, Ramm \cite{Ra87}  stated
the representation \eqref{10.9} (actually, \eqref{1.15}), but no proof was
given (nor was there any connection of the function $A$ to the inverse
problem for $V$). In \cite{Ra87} the inverse problem of finding the
potential from the knowledge of the $m$-function has been solved for
short-range  potentials. A more detailed discussion of the result in
\cite{Ra87}  can be found in \cite{Ra00}, \cite{Ra02}.

Throughout \cite{RS00} it is assumed that
\begin{equation} \label{8.1.4}
\int_0^\infty (1+x) \, dx \, |V(x)| < \infty
\end{equation}
and the Dirichlet-type Schr\"odinger operator $H$ in
$L^2([0,\infty))$ defined by
\begin{align}
\begin{split}
&H f=-f''+Vf, \quad 
 f\in \dom(H)=\{u\in L^2([0,\infty))\,|\, u, u' \in AC_{\loc}([0,b]) \\ 
& \hspace*{3.4cm}
\text{for all $b>0$}; \, u(0_+)=0; \, (-g''+Vg)\in L^2([0,\infty))\}   \lb{8.1.5}
\end{split}
\end{align}
is considered.

More generally, for $n\in\bbN_0$, $B\leq 0$ and $\ell 
\geq 0$, the space $C^{B,\ell}_n$ of all functions $q$ with $n-1$
classical  derivatives and $q^{(n)}\in L^1 ([0,\infty))$ so that
$$ 
\int_0^\infty (1 +x)^\ell\,  e^{-Bx} \,  dx \, |q^{(j)}(x)| < \infty
$$
for $j=0,1,\dots, n$. Thus, \eqref{8.1.4} says $V\in C^{B=0,
\ell=1}_{n=0}$.

Under condition \eqref{8.1.4}, general principles (see, e.g.,
\cite[Ch.\ 3]{Ma86})  imply that for all $\kappa\in\bbC$ with $\Re(\kappa)
\geq 0$, there is a unique  solution $F(\kappa,x)$ of $-f'' +
Vf=-\kappa^2 f$ normalized so that $F(\kappa,x) =  e^{-\kappa x} (1
+ o(1))$ as $x\to
\infty$.  We set $F(\kappa):= F(\kappa,0_+)$. Except 
for the change of variables $\kappa = -ik$, $F(\kappa,x)$ and $F(\kappa)$ 
are the standard Jost solution and Jost function. Both $F(\kappa,x)$ and 
$F(\kappa)$ are analytic with respect to $\kappa$ in $\{\kappa\in\bbC\,|\,
\Re(\kappa) > 0\}$.  If $V\in C^{B,\ell}_n$ for any $n,\ell$ and $B<0$,
then $F(\kappa,x)$ and $F(\kappa)$ have analytic continuations into the
region $\Re (\kappa) > B/2$.

The following is easy to see and well known (cf.\ \cite[Ch.\ 3]{Ma86}): 
\begin{enumerate}
\item[(1)] The zeros of $F$ in $\{\kappa\in \bbC\,|\,\Re(\kappa) >0\}$
occur precisely  at those points $\kappa_j$ with $-\kappa^2_j$ an
eigenvalue of the operator $H$ and each such zero is simple.
\item[(2)] $F$ has no zeros in $\{\kappa\in\bbC\,|\,\Re(\kappa)
=0,\,\kappa\neq 0\}$.
\item[(3)] If $F(0)=0$ and $V\in C^{B=0, \ell=2}_{n=0}$, then $F$ is $C^1$ 
and $F'(0)\neq 0$. If $F(0)=0$, we say that $H$ has a {\it zero energy 
resonance}.
\end{enumerate}

If $F$ can be analytically continued to $\{\kappa\in\bbC\,|\,\Re(\kappa)
>B/2\}$ for $B<0$, then zeros of $F$ in $\{\kappa\in\bbC\,|\,\Re(\kappa)
<0\}$ are called  {\it resonances} of $H$. They occur in complex
conjugate pairs (since $F$ is real  on the real axis). If
$F'(\kappa_0)\neq 0$ at a zero
$\kappa_0$, we say that 
$\kappa_0$ is a simple resonance. Resonances need not be simple if $\Re
(\kappa_0) <0$ although they are generically simple.

The result stated in \cite{Ra87} can be phrased as follows:

\begin{theorem} \label{t8.8.1} Suppose that $V$ satisfies \eqref{8.1.4}
{\rm(}i.e., it lies in 
$C^{B=0, \ell=1}_{n=0}${\rm)} and that $H$ does not have a zero energy
resonance.  Let $\{-\kappa^2_j\}^J_{j=1}$ be the negative eigenvalues of
$H$ with $\kappa_j >0$. Then
\begin{equation} \label{8.1.6}
A(\alpha) = \sum_{j=1}^J B_j \, e^{2\alpha\kappa_j} + g(\alpha),
\end{equation}
where $g\in L^1 ([0,\infty))$. In particular, if $H$ has no eigenvalues
and no zero energy resonance {\rm(}e.g., if $V\geq 0${\rm)}, then $A\in
L^1([0,\infty))$.
\end{theorem}

\begin{remark} $(i)$ The result stated in \cite{Ra87} assumes 
implicitly that there is no zero energy resonance.
Details can be found in \cite{Ra00}.  \\
$(ii)$ If $A\in L^1([0,\infty))$, then the representation \eqref{10.9}
(resp., \eqref{1.15}) can be analytically continued to  the entire region
$\Re(\kappa)
\geq 0$.  \\
$(iii)$ If $u_j$ is the eigenfunction of $H$ corresponding to the
eigenvalue $-\kappa^2_j$,  normalized by $\|u_j\|_2 = 1$, then
$$
B_j = -\f{|u'_j (0_+)|^2}{\kappa_j}\, .
$$
This follows from \eqref{8.16} 
and the fact that 
$$
d\rho(\lambda)\restriction (-\infty, 0) = \sum_{j=1}^J 
|u'_j(0_+)|^2 \delta (\lambda + \kappa^2_j)\, d\lambda.
$$
\end{remark}

To handle zero energy resonances of $H$, 
one needs an extra two powers of decay 
(just as \eqref{1.4} says more or less that $|V(x)|$ is bounded by 
$O(x^{-2-\veps})$, the condition in the next theorem says that 
$|V(x)|$ is more or less $O(x^{-4-\veps})$):

\begin{theorem} \label{t8.8.2} Let $V\in C^{B=0, \ell=3}_{n=0}$. Suppose
that $H$ has a zero energy 
resonance and negative eigenvalues at $\{-\kappa^2_j\}^J_{j=1}$ with 
$\kappa_j >0$. Then
\begin{equation} \label{8.1.7b}
A(\alpha) = B_0 + \sum_{j=1}^J B_j \, e^{2\alpha\kappa_j} + g(\alpha),
\end{equation}
where $g\in L^1 ([0,\infty))$.
\end{theorem}

These results are special cases of the following theorem:

\begin{theorem} \label{t8.8.3} Let $V\in C^{B=0, \ell}_n$ where $\ell \geq
1$  and, if $H$ has a zero energy resonance, then $\ell\geq 3$. Then
\eqref{8.1.6}  {\rm(}resp., \eqref{8.1.7b} if there is a zero energy
resonance{\rm)}  holds, where 
$g\in C^{B=0, \ell-1}_n$ {\rm(}resp., $C^{B=0, \ell-3}_n${\rm)}.
\end{theorem}

Finally, for $B<0$, the following result was proved in \cite{RS00}:

\begin{theorem} \label{t8.8.4} Let $V\in C^{B, \ell=0}_n$ with $B<0$. Let 
$\tilde B\in (B, 0)$ and let $\{-\kappa^2_j\}_{j=1}^J$ with $\kappa_j >0$ 
be the negative eigenvalues of $H$, $\{\lambda_j\}_{j=1}^M$ with
$\lambda_j\leq 0$  the real resonances {\rm(}a.k.a.~anti-bound
states{\rm)} of $H$, and $\{\mu_j 
\pm i\nu_j\}_{j=1}^N$ the complex resonances of $H$ with $\tilde B \leq
\mu_j <0$  and $\nu_j > 0$. Suppose each resonance is  simple. Then for
suitable 
$\{B_j\}_{j=1}^J$, $\{C_j\}_{j=1}^M$, $\{D_j\}_{j=1}^N$, 
$\{\theta_j\}_{j=1}^N$, one obtains
\[
A(\alpha) = \sum_{j=1}^J B_j \, e^{2\alpha\kappa_j} + \sum_{j=1}^M 
C_j \, e^{2\alpha\lambda_j} + \sum_{j=1}^N D_j \, e^{2\mu_j\alpha} 
\cos(2\nu_j\alpha + \theta_j) + \tilde g(\alpha), 
\]
where $\tilde g \in C^{\tilde B, \ell=0}_n$. In particular, if 
$H$ has no negative eigenvalues, the rate of decay of $A$ is 
determined by the resonance with the least negative value of $\lambda$ or 
$\mu$.
\end{theorem}

\smallskip
\begin{center}
$\bigstar$  \qquad $\bigstar$ \qquad $\bigstar$
\end{center}
\smallskip

We conclude this section with a brief look at the principal results in \cite{GS00b}. 

Let $H_j = -\f{d^2}{dx^2} + V_j$, $V_j\in L^1 ([0,b])$ for 
all $b>0$, $V_j$ real-valued, $j=1,2$, be two self-adjoint operators in 
$L^2 ([0,\infty))$ with a Dirichlet boundary condition at $x=0_+$. Let $m_j(z)$, 
$z\in\bbC\backslash\bbR$ be the Weyl--Titchmarsh $m$-functions associated with 
$H_j$, $j=1,2$. The main purpose of \cite{GS00b} was to provide a short proof of the following local uniqueness theorem in the spectral theory of one-dimensional 
Schr\"odinger operators, originally obtained by Simon \cite{Si99}, but under slightly more general assumptions than in \cite{Si99}. 

We summarize the principal results of \cite{GS00b} as follows:

\begin{theorem}\lb{t9.1.1} $(i)$ Let $a>0$, $0<\veps<\pi/2$ and 
suppose that
$$
|m_1 (z) - m_2(z)| \eqlim_{|z|\to\infty} 
O(e^{-2\Im (z^{1/2})a}) 
$$ 
along the ray $\arg(z) = \pi-\veps$. Then
$$
V_1(x) = V_2 (x) \, \text{ for a.e.\ } \, x \in [0,a].
$$
$(ii)$ Conversely, let $\arg(z)\in (\veps, \pi -\veps)$ 
for some $0<\veps<\pi$ and suppose $a>0$. If
$$
V_1(x) = V_2(x) \, \text{ for a.e.\ } \, x \in [0,a],
$$
then
\begin{equation} \lb{9.2.13}
|m_1 (z) - m_2(z)| \eqlim_{|z|\to\infty} 
O(e^{-2\Im (z^{1/2})a}).
\end{equation}
$(iii)$ In addition, suppose that $H_j$, $j=1,2$, are bounded 
from below. Then \eqref{9.2.13} extends to all $\arg(z) \in (\veps,\pi]$.
\end{theorem}

\begin{corollary}\lb{c9.2.6} Let $0<\veps <\pi/2$ and suppose 
that for all $a>0$,
$$
|m_1(z) - m_2(z)| \eqlim_{|z|\to\infty} 
O(e^{-2\Im (z^{1/2})a})
$$
along the ray $\arg(z) = \pi -\veps$. Then
$$
V_1(x) = V_2(x) \, \text{ for a.e.\ } \, x\in [0,\infty).
$$
\end{corollary}

Theorem \ref{t9.1.1} and Corollary \ref{c9.2.6} follow by combining some of the Riccati equation methods in \cite{GS00a} with properties of transformation operators (cf.\ 
\cite[Sect.\ 3.1]{Ma86}) and a uniqueness theorem for finite Laplace transforms 
\cite[Lemma\ A.2.1]{Si99}.  

In particular, Corollary \ref{c9.2.6} represents a considerable strengthening of the original Borg--Marchenko uniqueness result \cite{Bo52}, \cite{Ma50}, \cite{Ma73}:

\begin{theorem} \lb{t9.1.2}  
Suppose 
$$
m_1(z) = m_2(z), \quad z\in\bbC\backslash\bbR,
$$
then
$$
V_1(x) = V_2(x) \, \text{ for a.e.\ }  \, x\in [0,\infty).
$$
\end{theorem}
 
\begin{remark} \lb{r9.2.7} 
$(i)$ Marchenko \cite{Ma50} first published Theorem \ref{t9.1.2} in 1950. His extensive treatise on spectral theory of one-dimensional Schr\"odinger operators 
\cite{Ma73}, repeating the proof of his uniqueness theorem, then appeared in 1952, which also marked the appearance of Borg's proof of the uniqueness theorem \cite{Bo52} (apparently, based on his lecture at the 11th Scandinavian Congress of Mathematicians held at Trondheim, Norway, in 1949). \\
We emphasize that Borg and Marchenko also treat the general case of non-Dirichlet 
boundary conditions at $x=0_+$ (see also item $(vi)$ below). Moreover, Marchenko simultaneously discussed the half-line and finite interval case (cf.\ item $(vii)$ below).  \\
$(ii)$ As pointed out by Levitan \cite{Le87} in his Notes to Chapter\ 2, Borg and
Marchenko  were actually preceded by Tikhonov \cite{Ti49} in 1949, who
proved a special case of  Theorem~\ref{t9.1.2} in connection with the string equation (and hence under certain  additional hypotheses on $V_j$). \\
$(iii)$ Since Weyl--Titchmarsh functions $m$ are uniquely  related to the spectral measure $d\rho$ of $H$ by the standard Herglotz  representation theorem, 
\eqref{1.8b}, Theorem~\ref{t9.1.2} is equivalent to the following statement: Denote by 
$d\rho_j$ the spectral measures of $H_j$, $j=1,2$. Then
$$
d\rho_1 = d\rho_2 \, \text{ implies } \, V_1 = V_2 \, \text{ a.e.~on } \, [0,\infty).
$$
In fact, Marchenko took the spectral measures $d\rho_j$ as his point of 
departure while Borg focused on the Weyl--Titchmarsh functions $m_j$.  \\
$(iv)$ The Borg--Marchenko uniqueness result, Theorem~\ref{t9.1.2} (but not the  strengthened version, Corollary \ref{c9.2.6}), under the additional  condition of 
short-range potentials $V_j$ satisfying $V_j\in L^1 ([0,\infty); (1+x)\, dx)$, $j=1,2$, can also be proved using Property~C, a device used  by Ramm \cite{Ra99}, \cite{Ra00} in a variety of uniqueness results. \\
$(v)$  The ray $\arg(z) = \pi -\veps$, $0<\veps < \pi/2$ chosen in Theorem \ref{t9.1.1}\,$(i)$ and Corollary \ref{c9.2.6} is of no particular importance. A limit taken along any non-self-intersecting curve $\cC$ going to infinity in the sector 
$\arg(z)\in ((\pi/2)+\veps, \pi -\veps)$ will do since we can apply the 
Phragm\'en--Lindel\"of principle (\cite[Part~III, Sect.~6.5]{PS72}) to the region enclosed by $\cC$ and its complex conjugate $\ol\cC$. \\
$(vi)$ For simplicity of exposition, we only discussed the Dirichlet boundary condition
$$
u(0_+)=0
$$
in the Schr\"odinger operator $H$. Everything extends to the the general boundary condition
$$
u'(0_+) + h u(0_+) = 0, \quad h\in\bbR, 
$$
and we refer to \cite[Remark\ 2.9]{GS00b} for details.  \\
$(vii)$ Similarly, the case of a finite interval problem on 
$[0,b]$, $b\in (0,\infty)$, instead  of the half-line $[0,\infty)$ in Theorem~\ref{t9.1.1}\,$(i)$, with $0<a<b$, and a self-adjoint boundary condition at $x=b_-$ of
the type
$$
u'(b_-) + h_b u(b_-) = 0,  \quad h_b \in \bbR, 
$$
can be discussed (cf.\ \cite[Remark\ 2.10]{GS00b}).
\end{remark}

While we have separately described a few extensions in 
Remarks \ref{r9.2.7}\,$(v)$--$(vii)$,  it is clear that they can all be combined at once.

Without going into further details, we also mention that \cite{GS00b} contains the analog of the local Borg--Marchenko uniqueness result, Theorem \ref{t9.1.1}\,$(i)$ for Schr\"odinger operators on the real line. In addition, the case of half-line Jacobi operators and half-line matrix-valued Schr\"odinger operators was dealt with in \cite{GS00b}.

\medskip

{\it More recent references:} An even shorter proof of Theorem \ref{t9.1.1}\,$(i)$, close in spirit to Borg's original paper \cite{Bo52}, was found by Bennewitz \cite{Be01}. Still other proofs were presented by Horv\'ath 
\cite{Ho01} and Knudsen \cite{Kn01}. Various local and global uniqueness results for matrix-valued Schr\"odinger, Dirac-type, and Jacobi operators were considered in \cite{GKM02}. The analog of the local 
Borg--Marchenko theorem for certain Dirac-type systems was also studied by A.\ Sakhnovich \cite{Sa02}. The matrix-valued weighted Sturm--Liouville case has further been studied by Andersson \cite{An03}. He also studied uniqueness questions for certain scalar higher-order differential operators in \cite{An05}. A local Borg--Marchenko theorem for complex-valued potentials has been proved by Brown, Peacock, and Weikard \cite{BPW02}. The case of semi-infinite Jacobi operators with complex-valued coefficients was studied by Weikard \cite{We04}. A (global) uniqueness result for trees in terms of the (generalized) Dirichlet-to-Neumann map was found by Brown and Weikard \cite{BW05}. 
 
 \medskip

\noindent {\bf Acknowledgments.}
I'm grateful to Mark Asbaugh and the anonymous referee for a critical reading of this manuscript. I am also indebted to  T.\ Tombrello for the kind invitation to visit the
California Institute of Technology, Pasadena, CA, during the fall
semester of 2005, where parts of this paper were written. The great
hospitality of the Department of Mathematics at Caltech is gratefully
acknowledged. Moreover, I thank the Research Council and the
Office of Research of the University of Missouri-Columbia for granting me 
a research leave for the academic year 2005--2006.
 

\end{document}